\documentclass{article}

\usepackage{amsfonts}
\usepackage{amsmath}
\usepackage{amssymb}
\usepackage{xcolor}
\usepackage{dsfont}
\usepackage[a4paper, total={6.5in, 9in}]{geometry}
\usepackage{graphicx}

\usepackage{subcaption}

\usepackage{enumitem}

\usepackage{adjustbox,graphicx} 

\usepackage{multirow}
\usepackage{booktabs}
\usepackage{comment}
\usepackage{algorithm}
\usepackage{algpseudocode}

\usepackage{stmaryrd}
\usepackage{bm}

\usepackage{array} 

\usepackage{tikz}
\usetikzlibrary{shapes,arrows}
\usetikzlibrary{patterns,hobby}
\usetikzlibrary{calc}
\usetikzlibrary{shapes.geometric, arrows}
\usetikzlibrary{arrows}
\usepackage{pgfplots}
\pgfplotsset{compat=1.17}
\usetikzlibrary{shapes,arrows,positioning}
\usetikzlibrary{shapes.symbols,shapes.callouts,patterns}
\usetikzlibrary{calc}
\tikzset{
  cloud/.style = {draw, ellipse, fill=red!20, node distance=0.87cm, minimum height=2em},
  line/.style  = {draw, -latex'}
}
\usetikzlibrary{shapes.geometric}

\renewcommand{\t}{^{\mbox{\tiny\sf T}}}

\usepackage{aliascnt}

\usepackage{amsthm}
\usepackage{hyperref}
\usepackage{cleveref}

\newtheorem{theorem}{Theorem}[section]
\newaliascnt{corollary}{theorem}
\newaliascnt{lemma}{theorem}

\newtheorem{lemma}[lemma]{Lemma}
\newtheorem{corollary}[corollary]{Corollary}
\newtheorem{definition}{Definition}[section]
\newtheorem{remark}{Remark}[section]

\aliascntresetthe{corollary}

\aliascntresetthe{lemma}

\hypersetup{
    colorlinks=true,
    linkcolor=blue,
    filecolor=magenta,      
    urlcolor=cyan,
    citecolor=red
}

\newcommand{\timetot}{\tau\big(x_+(a),0\big)}

\newcommand{\nbtot}{\Sigma_s\big(x_+(a),0\big)}

\renewcommand{\t}{^{\mbox{\tiny\sf T}}}

\definecolor{vert}{rgb}{0.00, 0.55, 0.45}  
\definecolor{marronfonce}{RGB}{120,72,48}

\definecolor{bleu}{rgb}{0, 0.45, 1}

\definecolor{bleu_fonce}{RGB}{0,70,140}

\makeatletter
\renewcommand\paragraph{\@startsection{paragraph}{4}{\z@}%
            {-2.5ex\@plus -1ex \@minus -.25ex}%
            {1.25ex \@plus .25ex}%
            {\normalfont\normalsize\bfseries}}
\makeatother
\title{Sterile Insect Technique in a $n$-patch system with Allee effect and mass trapping:  modelling, analysis and simulations}

\author{Pierre-Alexandre Bliman${}^1$, Manon de la Tousche${}^{1,3}$, Yves Dumont${}^{2,3,4}$ \\
\small ${}^1$Sorbonne Université, Université Paris Cité, CNRS, Inria, Laboratoire Jacques-Louis Lions, LJLL, \\ 
\small EPC MUSCLEES, F-75005 Paris, France\\
\small ${}^2$CIRAD, UMR AMAP, Pôle de Protection des Plantes, F-97410 Saint Pierre, France,\\
\small ${}^3$AMAP, Univ Montpellier, CIRAD, CNRS, INRAe, IRD, Montpellier, France\\
\small ${}^4$University of Pretoria, Department of Mathematics and Applied Mathematics, Pretoria, South Africa\\
}

\def\diag{{\mathrm {\tt diag}}}

\date{\today}

\begin{document}

\maketitle

\begin{abstract}
    \newcommand\blfootnote[1]{%
        \begingroup
        \renewcommand\thefootnote{}\footnote{#1}%
        \addtocounter{footnote}{-1}%
        \endgroup
    }
    \blfootnote{The authors of this paper are listed in alphabetical order.}

     The sterile insect technique (SIT) is a biological control method aimed at reducing or eliminating populations of pests and disease vectors. This technique involves releasing sterilised insects which, by mating with wild individuals, will reduce the target population.
        In this study, spatial aspects are incorporated through an explicit metapopulation model, in which both wild and sterile insects disperse across multiple patches. We derive a general sufficient condition ensuring elimination of the wild population in all patches through SIT.
        When a strong Allee effect is naturally present  in each patch, the release programme can be terminated after a finite time, and an upper bound on the required release duration is provided. We then study an optimisation problem aimed at minimising the  daily number of sterile insects released to ensure 
        population elimination, under the constraint that the releases occur only in a prescribed subset of patches.  For illustrative purposes, we explore through numerical simulations different dispersal network configurations. We  also consider additional control by mass trapping (MT), which can affect the sterile insects entering trapped areas. 
        The main contribution of this work is to enable  the identification of the best spatial strategies for releasing sterile insects  under practical constraints, in order to achieve population elimination. 
\end{abstract}

\vspace{1cm}

\newpage

\tableofcontents

\newpage

\section{Introduction}

Pests and disease vectors have become a critical global concern, impacting both food security and public health.    In particular, fruit flies (Diptera: Tephritidae) cause substantial damage to vegetable crops and orchards. Among them, \textit{Bactrocera dorsalis}, commonly known as the oriental fruit fly, is a major pest primarily found in tropical and subtropical regions \cite{efsa2019bactrocera} and is considered one of the world's most invasive species \cite{clarke2005invasive,mutamiswa2021overview}. Females lay eggs beneath the fruit's surface, promoting bacterial and fungal growth that renders the fruit unmarketable.  In African countries where \textit{Bactrocera dorsalis} is established, high losses have been reported on mango crops.
For example, in Mozambique, the percentage of damaged mango fruits ranged from $21\%$ to $78\%$ between September 2014 and August 2015 \cite{cugala2020economic}. In Réunion Island, a French  territory near Madagascar, mango infestation by fruit flies  has greatly increased since the invasion of \textit{Bactrocera dorsalis}  in 2017 \cite{grechi2022bactrocera}.

To address these  issues, a  wide range of  control methods have been developed, with varying degrees of success. Among the large variety of control techniques  employed in the field, we focus on   the Sterile Insect Technique (SIT). This  biological control tool,  developed by Knipling in 1955 \cite{knipling1955possibilities},  involves the mass production of the target insect, followed by the sterilisation of  males (or both males and females, depending on the species) at the larval or  pupal stage  using ionising radiation. These sterilised insects are then released in large numbers into the field, where they mate with wild populations, leading to a gradual decline in the pest/vector population over time \cite{dyck2021sterile}.

 Motivated by Integrated Pest Management  and  Integrated Vector Management  strategies, we  also consider  additional treatments that increase mortality of the population.  Mass trapping (MT) methods are  widely used  in the field and are specifically designed to target and increase adult mortality,  through the combination of a semiochemical attractant, a killing agent---which is usually an insecticide \cite{barclay1991combining}---and a trapping device that is designed to attract and capture a large number of adult insects.  While MT is useful for reducing wild  populations, traps do not distinguish between wild and sterile individuals. As a result, they may also remove sterile  insects from the system, which can reduce the overall efficiency of SIT if not properly accounted for. The success of the combined strategy thus depends both on the spatial deployment of control measures and on their potential interactions.

Since Knipling’s early work, a wide variety of SIT models---simple, sex-structured, and/or stage-structured---have been developed and extensively studied to improve the effectiveness of SIT, depending on the target pest or vector. Some examples include probabilistic models \cite{berryman1967mathematical}, computational models \cite{markhof1973computer, diouf2022agent}, discrete models \cite{van1985some,li2015modelling,barclay2021mathematical}, semi-discrete models \cite{Huang,Strugarek2019,aronna2020nonlinear,bliman2019implementation,Dumont-Oliva2024}, and continuous models \cite{barclay1980sterile,esteva2005mathematical,bliman2019implementation,aronna2020nonlinear,anguelov2012mathematical,anguelov2020sustainable}, with some incorporating tools from control theory \cite{thome2010optimal,bliman2019implementation,almeida2022optimal,bliman2022,bliman2024optimal} to optimise release protocols  or to design feedback control. However, except in \cite{bliman2022}, these models typically assume that the targeted area is isolated from  surrounding areas—an idealised scenario that rarely holds in reality. In practice, the environment is heterogeneous, and both wild and sterile insects can disperse between neighbouring locations.

This spatial structure has important implications for  release  strategies.  For example, the application of SIT often faces spatial limitations when certain areas are inaccessible for direct releases.  Such constraints may arise from environmental conditions, such as dense forests or narrow mountainous terrain requiring aerial releases by helicopter for better manoeuvrability \cite{vargas1995aerial}, or from local stakeholders forbidding releases in their areas. A key advantage of SIT, however, is that sterile insects released in one location can disperse into neighbouring areas, potentially enabling suppression over a broader region even with localised releases.

Spatial dynamics of the SIT have been  mainly studied using partial differential equations (PDEs). Early work, such as \cite{lewis1993waves}, examined wave extinction phenomena resulting from  sterile insect releases, while later studies, including \cite{jiang2014reaction}, expanded on this framework. Other works, like \cite{DUFOURD2013},  also used PDEs to include environmental factors, such as wind and vegetation patterns.  However, the analysis of such models is challenging, especially in heterogeneous environments and when the spatial structure is complex.

     An alternative  to  PDE models is the  explicit metapopulation  approach.  In ecology, a metapopulation consists of a group of populations of the same species occupying
        distinct, homogeneous, and connected patches of suitable habitat within a larger landscape. Different modelling approaches have been developed to characterise metapopulation dynamics. In the classical, spatially implicit, approach introduced by Levins in 1969 \cite{levins1969some}, local population dynamics are not described explicitly; instead, the model tracks the fraction of occupied patches, assuming that individuals have the same probability of moving to any habitat patch. By contrast, spatially explicit metapopulation models keep track of the evolution of the  magnitude of the population in each patch and explicitly account for dispersal between patches, where movement typically depends on spatial proximity between habitat patches. This second  approach is adopted in this work.

     To our knowledge, all  existing SIT metapopulation models \cite{bliman2024efficacy,yang2020dynamics,dumont2025sterile_patch} have been restricted to the two-patch case. In particular,  a recent study \cite{dumont2025sterile_patch} established a sufficient condition for population elimination and subsequently optimised the  release rates. This study also formulated and solved an optimal control problem aimed at reducing the population   in a single patch below a prescribed threshold within finite time, while minimising the total number of sterile insects released. 

By contrast, in the modelling of mass trapping (MT),  the paper \cite{bliman2024feasibility} addressed the general case of $n$ patches, deriving optimisation results for complete elimination that depend on the spatial network structure and the number of  patches available for control.

\vspace{0.5cm}

 The present paper develops a general $n$-patch model combining MT and SIT, in which both wild and sterile insects can disperse between neighbouring patches, with potentially distinct dispersal networks for wild and sterile populations.  Unlike the two-patch case, the model accommodates arbitrary landscape configurations, including complex connectivity patterns and spatial heterogeneity that arise in realistic situations. In such a general framework, explicit analytical characterisations available in two-dimensional settings are no longer possible, making the development of new mathematical approaches necessary.

    A general  elimination condition  is thus derived for this $n$-patch system.
    Based on this result, we provide a framework for determining where sterile insects should be released under practical constraints, such as a limited number of feasible release sites, in order to efficiently eliminate the wild population throughout the landscape.

    Moreover,  unlike most SIT models, the  proposed framework  can incorporate a  natural  strong Allee effect  in each patch---which guarantees that the extinction equilibrium of the  system  without release is locally asymptotically stable---by adjusting the mating probability in each patch.
This allows to study elimination in finite time and to analyse how long  releases must be maintained. Specifically, once the system’s state enters the basin of attraction of the extinction equilibrium in the absence of  release, the population will not recover even if the   release efforts are subsequently relaxed.
     Below, we call release duration the time necessary to enter this basin.

The paper is organised as follows. In the next section, we introduce some notations and definitions that are relevant to our study. In \Cref{Pest/vector population model and preliminary results},  the  uncontrolled  model describing the natural evolution of the  wild population is  studied, first  in a homogeneous isolated area, and then in  a patchy environment. In particular, a   numerically tractable condition that ensures population elimination is provided.  Then  in \Cref{The n-patch sterile insect model}, we study the $n$-patch model describing   the evolution of sterile insects with constant releases. In \Cref{introduction_ control_chap3},  the models for wild and sterile insects are coupled, and the feasibility of population elimination,  under the constraint that the releases of sterile insects are limited to a prescribed subset of  patches, is studied. When the (sufficient) condition for  elimination feasibility is satisfied and  when  the origin of the uncontrolled system is locally asymptotically stable, we give an upper bound for the  release duration. Then,  in \Cref{minimization_section}, an optimisation problem consisting  of achieving  elimination while minimising a certain cost function, chosen as a non-decreasing and convex function of  the sterile insect release rates,  is studied.  Finally in \Cref{simulations}, we use numerical simulations on \textit{Bactrocera dorsalis} to study how different  network configurations  and practical constraints, including some where the releases are limited to a specific number of patches, influence  release strategies.   We also study SIT in combination with MT, and  the numerical explorations suggest that, although MT increases the critical  release rates required to eliminate the population, for fixed release rates it reduces the release duration, and thus decreases the total number of sterile insects released   over the entire
        programme. For reader’s convenience, most of the  proofs of the results are put in the Appendix.

\section{Notations and definitions}

\label{notations and definitions}

Let us introduce some general notations and definitions that will be used in this paper.

\begin{definition}
    For any $z \in \mathbb{R}$,  its \emph{positive part} is defined as
    $$z^+ := \max\{0,z\}.$$
\end{definition}

The inequalities between vectors  are considered in their usual coordinate-wise sense, that is, for any $x = (x_1,\ldots,x_n)$ and $y = (y_1,\ldots,y_n)$,
\begin{itemize}
    \item \(x \le y \Leftrightarrow x_i \le y_i, ~i = 1, \ldots, n\),
    \item \(x < y \Leftrightarrow x \le y, ~ x \neq y\),
    \item \(x \ll y \Leftrightarrow x_i < y_i, ~i = 1, \ldots, n\).
\end{itemize}
A vector $x$ is said to be   \textit{non-negative} if $x \geq 0_n$, \textit{positive} if $x > 0_n$, and \textit{strictly positive} if $x \gg 0_n$. 

These definitions are extended to matrices.

\vspace{0.5cm}

Now, let $f : \mathbb{R}^n \rightarrow \mathbb{R}^n$. The function $f$ is said:
\begin{itemize}
    \item \textit{non-decreasing} if $x \le y \Rightarrow f(x) \le f(y)$,
    \item \textit{increasing} if $x < y \Rightarrow f(x) < f(y)$,
\end{itemize}
Similarly, $f$ is said \textit{non-increasing} if $-f$ is non-decreasing, and \textit{decreasing} if $-f$ is increasing.

\begin{definition}
    \label{positive_real_line}
    The extended positive real number  line is denoted by
    $$\overline{\mathbb{R}}_+ := \mathbb{R}_+ \cup \{+ \infty\}.$$
\end{definition}

Now, we recall the definition of a Metzler  matrix.

\begin{definition}
    A matrix $A \in \mathbb{R}^{n \times n}$ is called \emph{Metzler} if all its off-diagonal components are non-negative. Moreover, $A$ is said to be \emph{irreducible} if it is not similar via a permutation to a block  triangular matrix;  otherwise, $A$ is said to be \emph{reducible}.
\end{definition}

\begin{definition}
    The \emph{stability modulus} of a  matrix $A \in \mathbb{R}^{n \times n}$ is  defined as
    $$s(A) := \{\max \{{  \Re(\lambda)} \}: \lambda \text{ is an eigenvalue of } A\}.$$
\end{definition}

\begin{definition}
    A   matrix $A \in  \mathbb{R}^{n \times n}$ is said \emph{Hurwitz} if $s(A) < 0$.
\end{definition}

We now introduce the following theorem. Its original form, known as the \emph{Perron--Frobenius theorem} (see, for instance,  \cite[Theorem 1.4, Chapter 2]{nonnegative}), is stated for  non-negative matrices. Its extension to the class of Metzler matrices follows directly from the spectral shift property of eigenvalues.

\begin{theorem}[Extension of the Perron-Frobenius theorem]
     Let $A \in \mathbb{R}^{n \times n}$  be a Metzler matrix.
        Then, $A$ has a positive right (resp. left) eigenvector. This right (resp. left) eigenvector is associated with a real eigenvalue, called the Perron-Frobenius eigenvalue, that is equal to $s(A)$.

    If moreover $A$ is irreducible, then $s(A)$ is a simple eigenvalue of $A$. Furthermore, $A$ has a strictly positive right (resp. left) eigenvector associated with the Perron-Frobenius eigenvalue, and any positive right (resp. left) eigenvector of $A$ is a scalar multiple of the previous one.
\end{theorem}

\begin{definition}
    The cardinal of any subset  $\mathcal{C}$ of $\{1,\ldots,n\}$ is denoted by  $n_\mathcal{C}$.
    For any $x = (x_1,\ldots,x_n)$, the point of $\mathbb{R}^{n_{\mathcal{C}}}$ composed of the $n_{\mathcal{C}}$ components of $x$ with index in $\mathcal{C}$ is denoted by  $x|_{\mathcal{C}}$:
    $$x|_{\mathcal{C}} := \left(x_i \right)_{i \in \mathcal{C}}.$$
    Similarly, for any matrix $A \in \mathbb{R}^{n \times n}$,  define the  submatrix
    $$A|_{\mathcal{C}} := \left(A_{ij} \right)_{i,j \in \mathcal{C}}.$$
\end{definition}

We now recall the definition  of the Kronecker delta.

\begin{definition}
    For any $x,y \in \mathbb{R}$, define the  Kronecker delta
    $\delta_x^y$ as
    $$\delta_x^y :=  \left\{
        \begin{array}{ll}
            1 & \mbox{if } x = y  \\
            0 & \mbox{otherwise.}
        \end{array}
        \right.$$
\end{definition}

 Last,
we   recall the definition of a cooperative and irreducible system of ODEs.

\begin{definition}\cite{smith1995monotone}
    \label{def cooperative}
    Let $\mathcal{D}$ be an open $p$-convex subset of $\mathbb{R}^n$, that is,
    $t x + (1-t) y \in \mathcal{D}$ for all $t \in [0,1]$ whenever $x,y \in \mathcal{D}$ and $x \le y$.
    Consider the system of ODEs $\dot{x} = f(x)$,
    where $f: \mathcal{D}\to \mathbb{R}^n$ is differentiable.
    The system is \emph{cooperative} on $\mathcal{D}$ if the Jacobian matrix of $f$
    is Metzler for all $x \in \mathcal{D}$. It is \emph{irreducible} on $\mathcal{D}$ if the Jacobian matrix
    of $f$ is irreducible for all $x \in \mathcal{D}$.
\end{definition}

\section{Pest/vector population model and preliminary results}

\label{Pest/vector population model and preliminary results}

In \Cref{One-patch model}, the  model describing the  evolution of the  wild population in a homogeneous isolated area is presented. Then in \Cref{n_patch}, we derive a general $n$-patch model describing the evolution of the population in presence of dispersal.

\subsection{One-patch, scalar model}
\label{One-patch model}

 Let us introduce the following dynamical system:
\begin{equation}
    \label{simple_allee}
    \dot{x} = x \big( p(x,a)b  - \mu_1  -\rho - \mu_2 x\big).
\end{equation}
The scalar $x(t)$ represents the number of  wild individuals  capable of reproducing (adults) at time $t$, and $\dot{x}$ is the derivative of $x$ with respect to time.   The parameters $b> 0$ and  $\mu_1> 0$ represent the natural daily birth rate per individual and the natural daily death rate, respectively, while $\rho \in \mathbb{R}_+$ denotes the additional mortality rate induced by mass trapping. 
        In fruit flies and mosquitoes,  the immature phase, gathering eggs, larvae and pupae, is subject to competition for space and food resources 
        \cite{burrack2009intraspecific,clarke2019biology,reinbold2021comparative}, as immature stages are unable to move to other sites when conditions become unfavourable. This intraspecific competition is taken into account by the term $-\mu_2 x^2$ in the model, with $\mu_2>0$. 

For  every $x,a \in \mathbb{R}_+$, the function
\begin{equation}
    \label{chap3_eq1}
    p(x,a):=  \left\{
    \begin{array}{ll}
        \dfrac{x}{x+a} & \mbox{if } a > 0 \\[3pt]
        1              & \mbox{if } a=0
    \end{array}
    \right.
\end{equation}
represents the density-dependent probability
of a successful mating.

When $a=0$, this probability is equal to $1$, and the growth is \textit{logistic}.  However, we know that in  practice, depending on several environmental and/or biological  factors, mating is not  always guaranteed.
More precisely, at low population densities, mating can become difficult, since males and females may  have difficulties to find each other. To account for this effect, we consider a probability of mating, $\dfrac{x}{x+a}$, where $a>0$ is a parameter that models the difficulty for the species  to mate \cite{durrett1994importance,dennis1989allee}. 
With this additional factor when $a>0$, we introduce what is called in ecology a  {\em  strong Allee effect} \cite{fauvergue2013review}, meaning that
     at low density, the per capita growth rate of the population is negative, which implies that  the extinction steady state is locally asymptotically stable. While an  Allee effect is an issue for species that need to be protected \cite{stephens1999consequences}, it  is  useful   to eliminate pests and disease vectors.
     Another common approach  to model a  strong Allee effect is to  introduce a factor $1-e^{-\beta x}$,  where  $\beta$ is
related  to the  mean distance an individual can effectively search for mates \cite{mccarthy1997allee}.   Both formulations lead to similar qualitative behaviour. We choose here the hyperbolic function $\frac{x}{x+a}$ because it has a direct interpretation as a mating probability and allows a unified treatment of both the low-density mating limitations and the impact of sterile insects. As will be seen later, the presence of sterile insects can be interpreted as an effective increase in the parameter $a$. In the remainder of this article, $a$ is referred to as the \emph{Allee parameter}.

When $a = 0$, the mating conditions are perfect and  the growth in model \eqref{simple_allee} is logistic.  We introduce the \emph{basic offspring number}
        \begin{equation}
            \label{N_chap3}
          \mathcal{N}:=  \frac{b}{\mu_1+\rho},
        \end{equation}
        which represents the average number of offspring produced by an individual during its lifetime.

          By denoting
        $$ Q:= \frac{\mu_1+\rho}{\mu_2},$$
        one easily shows  the following result.
\begin{theorem}
    \label{logistique_a=0}
    \mbox{}
    If $a = 0$, then the following holds for \eqref{simple_allee}:
    \begin{enumerate}
        \item If $\mathcal{N}\le 1$, then the origin $0$ is the only non-negative equilibrium and is globally asymptotically stable (GAS) on $\mathbb{R}_+$.
        \item If $\mathcal{N} > 1$,  then there exists  a  unique positive  equilibrium. Its value $x^*$ is given by
              $$x^* = Q \left(\mathcal{N}-1 \right).$$
              Moreover, this equilibrium is GAS on $\mathbb{R}^*_+$.
    \end{enumerate}
\end{theorem}

The following result is analogous to \Cref{logistique_a=0} in the case $a>0$. The equilibrium values were stated in \cite{dumont2025sterile_patch}; here, we also establish the basins of attraction and the stability properties of the equilibria, and a complete proof  is provided in \Cref{proof_of_allee_effet_simple}, which was not included in \cite{dumont2025sterile_patch}.

\begin{theorem}
    \label{allee_effet_simple}
    \mbox{}
    If $a > 0$, then the origin $0$ of \eqref{simple_allee} is locally asymptotically stable (LAS).

    Moreover, let
    $$ a^{\text{crit}} := Q\left(\sqrt{\mathcal{N}}-1 \right)^2.$$
    The following holds for \eqref{simple_allee}:
    \begin{enumerate}
        \item  If $\mathcal{N} \le 1$, then the origin $0$ is the only non-negative equilibrium  and is GAS on $\mathbb{R}_+$.
        \item If $\mathcal{N} > 1$,
              \begin{enumerate}[label = (\alph*)]
                  \item If
                        $$a> a^{\text{crit}},$$
                        then the origin is the only non-negative equilibrium and is GAS on $\mathbb{R}_+$.
                  \item If
                        $$a= a^{\text{crit}},$$
                        then there exists  only one positive equilibrium. Its value $x^*$ is given by
                        $$x^* =Q \left(\sqrt{\mathcal{N}}-1 \right).$$
                        Moreover, this equilibrium is unstable. The basins of attraction of $0$ and $x^*$ in $\mathbb{R}_+$ are  $[0, x^*)$ and $[x^*,+\infty)$, respectively.
                  \item If $$0< a< a^{\text{crit}},$$
                        then there exist  exactly two positive equilibria $x_{\pm}$,  with $x_- < x_+$. Their value is given by
                        $$x^*_{\pm} = \frac{Q}{2}\left(\mathcal{N}- 1 - \frac{a}{Q} \right)\left( 1 \pm \sqrt{1 - \frac{4a}{Q\left(\mathcal{N}-1- \frac{a}{Q} \right)^2}}\right).$$
                        Moreover, $x_-$ is unstable and $x_+$ is stable. The basins of attraction of $0$ and   $x_+$ in $\mathbb{R}_+$ are  $[0, x_-)$ and $(x_-, +\infty)$, respectively.
              \end{enumerate}
    \end{enumerate}
\end{theorem}

Now define the per capita growth rate
$$g(x) :=  p(x,a)b - \mu_1 - \rho - \mu_2 x,$$
for $p(x,a)$ defined in~\eqref{chap3_eq1}.
By \cite[Section 3.2]{dumont2025sterile_patch},
\begin{equation}
    \label{max_equality}
    \max_{x \geq 0} g(x) = - \mu_1 - \rho + \left(\left( \sqrt{b} - \sqrt{a \mu_2}  \right)^{+}\right)^2.
\end{equation}

We then establish the following result for population elimination. It is essentially a reformulation of the previous theorem, but with elimination now depending on the value of $ \max_{x \geq 0} g(x)$. This result will serve as a basis for the generalisation to the case where the species evolves across interconnected areas,  as discussed later in \Cref{n_patch}.

\begin{theorem}
    \label{max}
    The origin of \eqref{simple_allee} is GAS on $\mathbb{R}_+$ whenever
    \begin{equation}
        \label{condition_elim}
        - \mu_1 - \rho + \left(\left( \sqrt{b} - \sqrt{a \mu_2} \right)^+\right)^2 < 0.
    \end{equation}
    Moreover, if $\mathcal{N}>1$, this condition is both necessary and sufficient for global asymptotic stability of the origin on $\mathbb{R}_+$.
\end{theorem}

\begin{proof}
    See \Cref{proof_max}
\end{proof}

When the population evolves over $n$   interconnected patches (see \Cref{n_patch}), an analogous stability condition will be derived, which generalises \Cref{max} (see \Cref{convergence_allee} below).

\subsection{The $n$-patch pest/vector population model}
\label{n_patch}

In \Cref{One-patch model}, the model \eqref{simple_allee}  typically assumes that the  area in which the insects evolve is   homogeneous and is isolated from  surrounding areas—an idealised scenario that rarely holds in reality.  In practice, the environment is heterogeneous, and  insects can disperse between neighbouring locations.    To incorporate this spatial structure, we now extend the model to a metapopulation framework consisting of $n$ patches that are  interconnected.

        We first present the $n$-patch model in \Cref{Presentation of the model} and some local stability properties of
        the extinction equilibrium in \Cref{local_stability_properties_chap3}. Then in \Cref{maximal_equilibrium_chap3}, the existence of a maximal (i.e., componentwise largest) equilibrium is established, together with a stability property of this equilibrium.  In \Cref{global_elimination}, a sufficient condition for global elimination is provided. Finally in \Cref{section_monotonicity_properties}, we establish some monotonicity properties  relative to the Allee parameter. These results will be applied later in \Cref{introduction_ control_chap3}, when incorporating the release of sterile insects into the model.

\subsubsection{Wild population model}
\label{Presentation of the model}

The $n$-patch metapopulation model obtained from the one-patch model \eqref{simple_allee} is the following  system:
\begin{equation}
    \label{log_allee}
    \dot{x}_i = x_i \big(p(x_i,a_i) b_i  - \mu_{1,i}  - \rho_i- \mu_{2,i} x_i\big) + \sum\limits_{\substack{j = 1\\ j \neq i}}^{n} d_{ij} x_j -  \sum\limits_{\substack{j = 1\\ j \neq i}}^{n} d_{ji} x_i, \qquad i = 1,\ldots,n,
\end{equation}
which  extends  the logistic-type growth framework studied in \cite{bliman2024feasibility}  by incorporating an Allee effect   in the patches. In the remaining of the present paper, system \eqref{log_allee} is referred to as the \emph{uncontrolled system}.

For every $i = 1,\ldots,n$,
$x_i(t)$ is a scalar representing the number of individuals in the $i$-th patch at  time $t$ and $\dot{x}_i$ is its time derivative. The parameters $b_i, \mu_{1,i}$, $\mu_{2,i}$ are positive,  $\rho_i \in \mathbb{R}_+$, and the function
\begin{equation}
    \label{g_i}
    g_i(x_i) :=    p(x_i,a_i)b_i  - \mu_{1,i} - \rho_i- \mu_{2,i} x_i
\end{equation}
is  the   per capita growth rate  of the population in the $i$-th patch, for $p$ defined in~\eqref{chap3_eq1}.  In particular, having $a_i>0$ corresponds to a strong Allee effect in  Patch $i$, since as mentioned in  \Cref{One-patch model}, the local per capita growth rate $g_i$ becomes negative when the density in the patch is low. Note however that, contrary to the one-patch model \eqref{simple_allee}, a strong Allee effect in a patch does not automatically imply the extinction of its population when its density is low, as immigration from other patches can potentially rescue it from collapsing.

Finally, $d_{ij} \in \mathbb{R}_+$ is a non-negative dispersal coefficient describing the flow of individuals from  Patch $j$ to Patch $i$ ($i \neq j$). Thus,  define the  connectivity matrix $D$ of the wild population as follows:
$$
    D_{ij} := \begin{cases}
        d_{ij} & \text{ if } i \ne j \\
        -\sum\limits_{\substack{k=1  \\ k \neq i}}^{n} d_{ki} & \text{ if } i=j.
    \end{cases}
$$
This matrix is  Metzler, i.e.,  its non-diagonal components are non-negative. In some results of this paper, we will assume that  $D$ is irreducible, i.e., that the underlying directed graph is irreducible, implying that
the species can  migrate between any patches.

Denote $x:= (x_1,\ldots,x_n)$, $a: = (a_1,\ldots,a_n)$, and  define $f_a(x)$ as the right-hand side of  system \eqref{log_allee}, i.e.,
\begin{equation}
    \label{f_a}
    f_a(x) := \diag \big(
    p(x_i,a_i) b_i   - \mu_{1,i} - \rho_i- \mu_{2,i} x_i
    \big) x + Dx.
\end{equation}
When $n = 2$,  this type of system has been studied in the context of sterile insect releases in \cite{dumont2025sterile_patch}. This section  generalises some  results from that paper to an arbitrary number $n$ of patches.

\subsubsection{Local stability properties}

\label{local_stability_properties_chap3}

This section establishes some local stability properties of the extinction equilibrium of \eqref{log_allee}. To this aim, denote by $J_a$ the Jacobian matrix at the origin  of the function $f_a$ defined in \eqref{f_a}. It writes
\begin{equation}
    \label{J_a}
    J_a = \diag\left(b_i  \delta_{a_i}^0 - \mu_{1,i} - \rho_i\right) + D,
\end{equation}
where $\delta^0_{a_i}$ denotes the Kronecker delta.

\begin{theorem}[Local stability of the extinction equilibrium]
    \label{0LAS}
    For any $a \in \mathbb{R}^n_+$, if $s(J_a) < 0$, then the origin of \eqref{log_allee} is LAS.

    In particular, if $a \gg 0_n$, then the origin is LAS.
\end{theorem}

\begin{proof}
    When $s(J_a) < 0$, the origin of \eqref{log_allee} being LAS is a classical result in stability theory.

    When $a \gg 0_n$, one has $ J_a = -\diag(\mu_{1,i}+\rho_i) + D$. Since $\mathbf{1}_n^T J_a = -\left(\mu_1+\rho_1,\ldots,\mu_n+\rho_n \right)$, the  Metzler matrix $J_a$ is strictly column diagonally dominant, and is therefore Hurwitz, implying the local asymptotic stability of the origin.
\end{proof}

\begin{remark}
    From \Cref{0LAS}, a strong Allee effect in each patch ($a \gg 0_n$) implies the local asymptotic stability of the origin. However, the local asymptotic stability of the origin does not imply the presence of a  strong Allee effect in every patch. Consider for example  a two-patch model ($n=2$) where Patch 1 does not exhibit a (strong) Allee effect ($a_1=0$) while Patch 2 does ($a_2>0$).
    Assuming Patch 1 is intrinsically viable ($b_1 > \mu_{1,1} + \rho_1$), the origin can still become LAS  if the dispersal rate from the viable patch to the patch with strong Allee effect satisfies:
    \begin{equation*}
        d_{21} > (b_1-\mu_{1,1}-\rho_1)\left(1+ \frac{d_{12}}{\mu_{1,2}+\rho_2} \right).
    \end{equation*}
    In this scenario, the extinction of the population when densities are low in both patches is not caused by a strong  Allee effect in each patch,
    but rather by a dispersal-induced phenomenon where the emigration rate $d_{21}$ drains individuals from Patch 1 into Patch 2 faster than Patch 1's intrinsic low-density growth rate can compensate, while Patch 2 acts as a sink due to its strong Allee effect.
\end{remark}

\begin{definition}
    Denote by $\mathcal{B}_a$ the basin of attraction of the origin of system \eqref{log_allee}.
\end{definition}

In the next result, we give an approximation of $\mathcal{B}_a$.

\begin{theorem}[Estimation of the basin of attraction]
    \label{estimation basin of attraction}
    Let  $a > 0_n$. Assume $D$ is irreducible and  $s(J_a) < 0$, and let $c \gg 0_n$ be a corresponding left  (Perron) eigenvector of $J_a$,  normalised so that $\min_i c_i = 1$. Then,
    \begin{equation*}
        \left\{ x \in \mathbb{R}^n_+: \Psi\left(c^T x\right) < -s(J_a)  \right\}  \subset \mathcal{B}_a,
    \end{equation*}
    where,
    for every $z \in \mathbb{R}_+$,
    $$\Psi(z):=   \max_{\substack{1 \le k \le n \\  a_k > 0}}  \phi_k(z), \qquad  \phi_k(z) :=
        \begin{cases}
            b_k \dfrac{z}{z + a_k} - \mu_{2,k} z                               & \text{if } z \le -a_k + \sqrt{\dfrac{a_k b_k}{\mu_{2,k}}} \\
            \left( \left( \sqrt{b_k} - \sqrt{a_k \mu_{2,k}} \right)^+\right)^2 & \text{if } z > -a_k + \sqrt{\dfrac{a_k b_k}{\mu_{2,k}}}
        \end{cases}.$$
    Moreover, the previous set contains the origin.
\end{theorem}

\begin{proof}
    See \Cref{proof_estimation_basin}.
\end{proof}

\subsubsection{Maximal equilibrium}

\label{maximal_equilibrium_chap3}

Using the theory of cooperative systems, we  establish the existence of a maximal equilibrium of system \eqref{log_allee}.

\begin{theorem}[Properties of the maximal equilibrium]
    \label{0only}
    For any $a \in \mathbb{R}^n_+$, the system \eqref{log_allee} admits a non-negative equilibrium   $x_+(a)$  that is maximal,  in the sense that any other equilibrium $x^*_a$ satisfies  $x^*_a \le x_+(a)$. Moreover, $x_+(a)$ is GAS on the set $\{x \geq x_+(a) \}$.

    Last,
    if $D$ is irreducible,  then any two non-negative equilibria $x^*_a$ and  $\tilde{x}^*_a$ satisfy
    \begin{equation}
        \label{positive_eq_irred_struct}
        x^*_a < \tilde{x}^*_a  \Rightarrow  x^*_a \ll \tilde{x}^*_a.
    \end{equation}
\end{theorem}

\begin{proof}
    See \Cref{prooof_0only}.
\end{proof}

Notice that when $D$ is irreducible and since the origin is an equilibrium, it follows from \eqref{positive_eq_irred_struct} that either the population is eliminated in all patches or it persists in all of them.

\begin{remark}
    \label{remark_0only}
    An immediate consequence of \Cref{0only} is that if $x_+(a) = 0_{n}$, then it is GAS on $\mathbb{R}^n_+$.
\end{remark}

\begin{remark}
    The existence of a unique positive equilibrium  for system \eqref{log_allee} does not necessarily imply its global asymptotic stability on $\mathbb{R}^n_+ \setminus \{ 0_n\}$. Indeed, in such case and when
    the origin is LAS, then by  monotonicity of  system \eqref{log_allee}, the origin attracts every trajectory starting in the set  $\{x: 0_n \le x < x_+(a) \}$ (\cite[Theorem 2.2, Chapter 2]{smith1995monotone}). This represents a significant difference compared to the case of
    perfect mating conditions ($a = 0_n$), for which the existence of a  positive equilibrium implies both its uniqueness and its global asymptotic stability on  $\mathbb{R}^n_+\setminus \{0_n\}$.
\end{remark}

\subsubsection{A sufficient condition for  elimination}

\label{global_elimination}

In this section, we seek a condition guaranteeing that the extinction equilibrium of system \eqref{log_allee} is GAS on $\mathbb{R}^n_+$, i.e., that the maximal equilibrium of \eqref{log_allee} is null. To this aim,
for all $a \in \mathbb{R}^n_+$,   define the matrix
\begin{equation}
    \label{A_a}
    A(a) := \diag\left(- \mu_{1,i} - \rho_i+ \left(\left( \sqrt{b_i} - \sqrt{a_i \mu_{2,i}} \right)^+\right)^2\right) + D.
\end{equation}
As noted in the  discussion immediately preceding \Cref{max}, each diagonal entry in $\diag()$ corresponds to
\begin{equation}
    \label{max_g_i}
    \max_{x_i \geq 0} g_i(x_i),
\end{equation}
where $g_i$ is the  per capita growth rate  in Patch $i$, defined in \eqref{g_i}. 
Thus, the matrix $A(a)$ is the Jacobian matrix of system \eqref{log_allee} evaluated at the unique point in the state space where  the per capita growth rate is maximal in each patch.

Building on these notations,  the following theorem is a direct application of  \cite[Theorem 1]{takeushi93}.

\begin{theorem} [Sufficient condition for global asymptotic stability of the origin]
    \label{convergence_allee}
    For any $a \in \mathbb{R}^n_+$, if $s(A(a)) < 0$, then $x_+(a) = 0_n$ and the origin of \eqref{log_allee}  is GAS on $\mathbb{R}^n_+$.
\end{theorem}

In the one-patch case ($n=1$), this theorem recovers the sufficient condition stated in \Cref{max}. Notice that when $\mathcal{N} >1$,  for $\mathcal{N}$ defined in \eqref{N_chap3}, this condition is also  necessary by \Cref{max}.

In \Cref{convergence_allee}, only a sufficient condition for population elimination is  established. When the mating conditions are perfect, i.e., when $a=0_n$,  a stronger result can be obtained. In such case, the matrix $J_a$,  defined in \eqref{J_a}, satisfies
$$J_a = J_0 = \diag(b_i - \mu_{1,i}-\rho_i)+D,$$
and we have
$$A(0) = J_0.$$
 This case has been studied in \cite{bliman2024feasibility}. In particular, when $D$ is irreducible, a necessary and sufficient condition for population elimination has been established in Theorem 2.2 of that paper, based on the sign of the stability modulus $s(J_0)$ of the matrix $J_0$. We extend this result to the case where $D$ is not required to be irreducible.

\begin{theorem}
    \label{elimination_model_uncontrolled_chapter3}
     If $a=0_n$, then the following holds:
    \begin{enumerate}
        \item  If $s(J_0) \le 0$, then $x_+(0) = 0_n$.
        \item  If $s(J_0) > 0$, then $x_+(0) > 0_n$. Moreover, $x_+(0) \gg 0_n$ if $D$ is irreducible.
    \end{enumerate}
    \label{convergence2}
\end{theorem}

This theorem is a generalisation of the result given in \Cref{logistique_a=0} for the one-patch model \eqref{simple_allee}. As a matter of fact, when $n=1$, the condition $s(J_0) \le 0$ is equivalent to $\mathcal{N} \le 1$.

\begin{proof}[Proof of \Cref{elimination_model_uncontrolled_chapter3}]
    When $D$ is irreducible, \Cref{elimination_model_uncontrolled_chapter3} is a direct application of  \cite[Theorem 2.2]{bliman2024feasibility}.  When $D$ is not irreducible, one can decompose the underlying directed graph into its \emph{strongly connected components} (\emph{SCCs}), which form irreducible subgraphs  (see \Cref{SCC_chap3} introduced later in \Cref{The n-patch sterile insect model}).
    Since $D$ (and thus $J_0$)  is reducible,  there exists a permutation matrix such that $J_0$ can be brought into block  triangular form where each block on the diagonal corresponds to an SCC, and thus is irreducible (see \cite{varga1999matrix} p.~50). Therefore, if $s(J_0) \le 0$, then all the diagonal blocks have a non-positive stability modulus, and we deduce as a consequence of  \cite[Theorem 4.2]{bliman2024feasibility} that $x_+(0) = 0_n$. Conversely if $s(J_0) > 0$, then there exists at least one diagonal block with a strictly positive stability modulus, and we deduce, still as a consequence of \cite[Theorem 4.2]{bliman2024feasibility}, that $x_+(0) >0_n$.
\end{proof}

\subsubsection{Monotonicity properties}

\label{section_monotonicity_properties}

In this section, we establish some important monotonicity properties with respect to the Allee parameter $a$.
These results will be applied to the SIT model in \Cref{introduction_ control_chap3}, since, as will be seen, the
release of  sterile insects   artificially induces or amplifies a  strong Allee effect in some patches, through an increase in the Allee parameter $a$.

\begin{theorem}
    \label{max_eq_decreasing}
    For any $a \in \mathbb{R}^n_+$, the maximal equilibrium $x_+(a)$ of system \eqref{log_allee} is non-increasing with $a$, that is, for any $a,a'\in \mathbb{R}^n_+$,
    \begin{equation}
        \label{decrease_eq}
        a < a' \Rightarrow x_+(a') \le x_+(a).
    \end{equation}

    Moreover, if $D$ is irreducible and $x_+(a) > 0_n$, then the inequality $\le$ in \eqref{decrease_eq} may be replaced by the inequality $\ll$.

    Last,    the basin of attraction  $\mathcal{B}_a$ of the origin of system \eqref{log_allee} is non-decreasing with $a$, that is, for any $a,a' \in \mathbb{R}^n_+$,
    $$a \le a' \Rightarrow \mathcal{B}_a \subset \mathcal{B}_{a'}.$$
\end{theorem}

\begin{proof}
    See \Cref{proof_of_max_eq_decreasing}.
\end{proof}

 In particular, \Cref{max_eq_decreasing} states that when $D$ is irreducible, then the increase of $a_i$ in \emph{any} patch reduces the maximal equilibrium value in \emph{every} patch of the network, which constitutes a  strong result from a population reduction perspective.

From the inequality \eqref{decrease_eq} in \Cref{max_eq_decreasing} and the fact that $x_+(a)$ is GAS on the set $\{x\geq x_+(a)\}$, one derives the  following result.

\begin{corollary}
    \label{corollary_max_eq}
    Let $a \in \mathbb{R}^n_+$. If the origin of \eqref{log_allee} is GAS on $\mathbb{R}^n_+$, then the origin of system \eqref{log_allee} with parameter $a'$ is GAS on $\mathbb{R}^n_+$ for all $a' \geq a$.
\end{corollary}

Finally, the following result provides monotonicty  and convexity properties of the stability modulus of $A(a)$. It will be
useful later in \Cref{minimization_section} to study an optimisation problem related to the elimination of the wild
population through the release of sterile insects.

\begin{theorem}
    \label{s(A(a))_decreasing}
    The function   $a \in \mathbb{R}^n_+ \mapsto s(A(a))$ is convex and  non-increasing in the following sense:
    \begin{equation}
        \label{non-decreasing A(a)}
        a < a' \Rightarrow s(A(a')) \le s(A(a)).
    \end{equation}

    Last, if $D$ is irreducible,   then  this function is  continuous on $\mathbb{R}^n_+$, and
    continuously differentiable and piecewise twice continuously differentiable on the set  $\{ a \gg  0_n \}$. If in addition, $a_i < \frac{b_i}{\mu_{2,i}}$ for all $i =1,\ldots,n$, then
    the inequality $ \le $ in  \eqref{non-decreasing A(a)}
    may be replaced by the inequality $<$.
\end{theorem}

\begin{proof}
    See \Cref{proof_s(A)_decreasing}.
\end{proof}

\begin{remark}
    \label{remarkyves}
    By \Cref{s(A(a))_decreasing}, we have $s(A(a))\leq s(A(0))$, for all $a\in \mathbb{R}^n_+$.
\end{remark}

 We now summarise the main results of this section. We have  shown that a strong Allee effect in each patch ($a \gg 0_n$) implies the local asymptotic stability of the origin of \eqref{log_allee} (\Cref{0LAS}) and  an estimate of the basin of attraction of the extinction equilibrium  is provided in \Cref{estimation basin of attraction}. We have also proved the existence of a maximal equilibrium (\Cref{0only}). Then, in \Cref{convergence_allee}, we have derived a sufficient condition for this maximal equilibrium to be null, thereby guaranteeing elimination of the wild population in all patches. This condition depends on the stability modulus of the Jacobian matrix of system \eqref{log_allee} evaluated at the unique point in the state space where  the per capita growth rate is maximal in each patch. Finally, we have established monotonicity properties of the maximal equilibrium (\Cref{max_eq_decreasing}) and of the stability modulus (\Cref{s(A(a))_decreasing}) with respect to the Allee parameter $a$, in order to  anticipate the effect of sterile insect releases, which is studied later in \Cref{introduction_ control_chap3}.

\section{The $n$-patch sterile insect model}

\label{The n-patch sterile insect model}

We now consider the release of sterile insects.
In \Cref{Presentation of the SIT model}, the  model describing the evolution of the sterile population is presented, along with some relevant boundedness assumptions on the release rate.
Then in \Cref{Asymptotic behaviour of the sterile insects}, the asymptotic behaviour of the sterile insects is studied. Finally, in \Cref{Supremal sterile insect equilibrium}, the supremal sterile insect equilibrium  for large admissible release rates is determined.

\subsection{Sterile insect model}
\label{Presentation of the SIT model}

We assume that sterile insects are released in each patch at a constant  rate $\Lambda_i \in \mathbb{R}_+$ for $i = 1,\ldots,n$, with a mortality rate of $\mu_{s,i}>0$. The population of sterile individuals in patch $i$, denoted as $x_{s,i}$, is then  modelled by:

\begin{equation}
    \label{sterile_model}
    \dot{x}_{s,i} = \Lambda_i - \left(\mu_{s,i} + \rho_{s,i}\right)x_{s,i} + \sum\limits_{\substack{j = 1 \\ i \neq j}}^{n}d^s_{ij} x_{s,j} -  \sum\limits_{\substack{j = 1 \\ i \neq j}}^{n}d^s_{ji} x_{s,i}, \qquad i = 1,\ldots,n.
\end{equation}
 The scalar $\rho_{s,i} \in \mathbb{R}_+$ represents the additional mortality rate induced by  mass trapping in Patch $i$. Indeed, in the field,  traps targeting wild insects  may also capture sterile insects. We do not impose $\rho_i = \rho_{s,i}$, as the relative susceptibility to trapping of sterile  insects  depends on the types of traps used.
 Since sterile insects are released as adults, there is no inter- or intraspecific competition for resources, such that there is no density-dependent mortality rate as in the wild population.

Similarly to the coefficients $d_{ij}$ introduced in~\Cref{n_patch} to describe the dispersal of  wild insects, $d^s_{ij}$ is a non-negative coefficient describing the flow of   sterile insects from  Patch $j$ to  Patch $i$ ($i \neq j$). We then introduce the matrix $D^s$, where
\begin{equation*}
    D^s_{ij} := \begin{cases}
        d^s_{ij} & \text{ if } i \ne j \\
        -\sum\limits_{\substack{k=1    \\ k \neq i}}^{n} d^s_{ki} & \text{ if } i=j.
    \end{cases}
\end{equation*}
Note that system \eqref{sterile_model} writes in vector form as
$$\dot{x}_s = \Lambda + J^s x_s,$$
where $x_s := \left(x_{s,i} \right)_{i = 1,\ldots,n}$, $\Lambda = \left(\Lambda_i \right)_{i = 1,\ldots,n}$, two vectors, and
\begin{equation}
    \label{def cal A}
     J^s := D^s - \diag(\mu_s+\rho_s),
\end{equation}
with $\mu_s = \left(\mu_{s,i} \right)_{i = 1,\ldots,n}$ and $\rho_s = \left(\rho_{s,i} \right)_{i = 1,\ldots,n}$.

We consider that in the field,  the release of sterile insects may be limited and restricted to a specific subset of plots, due, for example,
to environmental constraints that prevent control implementation or to the reluctance of local stakeholders to allow certain forms
of intervention.   Accordingly, we  introduce a boundedness assumption on the sterile insect release rate $\Lambda$  in~\eqref{sterile_model},   recalling that the extended positive real line $\overline{\mathbb{R}}_+$ has been defined in \Cref{positive_real_line}.

\begin{definition}
    \label{def C^s}
    Let   $\overline{\Lambda} \in \overline{\mathbb{R}}_+^n$.  Define the sets $\mathcal{C}^s_B$, $\mathcal{C}^s_I$ and $\mathcal{C}^s$ as follows:
    \begin{subequations}
        \begin{equation}
            \mathcal{C}^s_B := \{ i \in \{1,\ldots,n\}\ :\ 0 < \overline{\Lambda}_i < + \infty   \},\qquad
            \mathcal{C}^s_I := \{ i \in \{1,\ldots,n\}\ :\ \overline{\Lambda}_i = + \infty   \}
        \end{equation}
        and
        \begin{equation}
            \mathcal{C}^s := \mathcal{C}^s_B \cup \mathcal{C}^s_I.
        \end{equation}
    \end{subequations}
\end{definition}

\begin{definition}
    \label{Lambda-admissible}     Let $\overline{\Lambda} \in \overline{\mathbb{R}}^n_+$. Any vector $\Lambda\in\mathbb{R}^n_+$ such that $\Lambda \leq \overline{\Lambda}$ is called $\overline{\Lambda}$\emph{-admissible}.
\end{definition}

\subsection{Asymptotic behaviour of the sterile insect population}
\label{Asymptotic behaviour of the sterile insects}

In this section, we study the asymptotic behaviour of the solutions of system \eqref{sterile_model}. First, the fact that the matrix $J^s$ defined in \eqref{def cal A} is Metzler and Hurwitz (it is strictly column diagonally dominant) implies that  $J^s$ is invertible, with $-\left(J^s\right)^{-1} > 0_{n \times n}$ (\cite[Theorem 2.5.3]{horn1994topics}).   Thus,
        for any  sterile insect initial condition $x_{s,0}\in \mathbb{R}^n_+$, the solution  $x_s(t, x_{s,0})$ of \eqref{sterile_model}  satisfies
        \begin{equation}
            \label{solution_sterile_control}
            x_s(t, x_{s,0}) = e^{J^s t}x_{s,0}+ \left(J^s\right)^{-1}  \left( e^{J^s t}-I \right)\Lambda.
        \end{equation}
When $D^s$ is irreducible, we more precisely have $-\left(J^s\right)^{-1} \gg 0_{n \times n}$ (\cite[Theorem 10.3]{bullo2018lectures}). These observations yield the following result.

\begin{theorem}
    \label{eq_steriles}
    For any $\overline{\Lambda}$-admissible $\Lambda \in \mathbb{R}^n_+$ the system \eqref{sterile_model} admits a unique  equilibrium  $x^*_s(\Lambda) = -\left(J^s\right)^{-1} \Lambda \in \mathbb{R}^n_+$ which is GAS.
    In particular,  it is increasing with $\Lambda \in \mathbb{R}^n_+$, and the following holds:
    \begin{itemize}
        \item  If $\Lambda > 0_n$, then $x^*_s(\Lambda)$ is positive $(x^*_s(\Lambda) > 0_n)$.
        \item  If $D^s$ is irreducible and $\Lambda > 0_n$, then $x^*_s(\Lambda)$ is strictly positive $(x^*_s(\Lambda) \gg 0_n)$.
    \end{itemize}
\end{theorem}

In  what follows,   the value of $x^*_s(\Lambda)$  on the patches is characterised more precisely, depending on the connectivity matrix $D^s$  of sterile insects and  the release vector $\Lambda$.  When the   network of the sterile insects is irreducible, the sterile insect equilibrium is strictly positive. When it is not irreducible, a  detailed examination of the network is required.  Introducing $\Gamma_s$ as  the underlying graph of the matrix $D^s$, we revert to the irreducible case by decomposing $\Gamma_s$ into its {\em strongly connected components} (SCCs), which are irreducible subgraphs (see \Cref{SCC_chap3}). Such a diakoptic approach has been used to determine asymptotic stability of linear cooperative systems \cite{greulich2019stability} and the GAS equilibrium of system \eqref{log_allee} when $a = 0_n$  \cite{bliman2024feasibility}.

Let us  introduce some adequate notations and definitions from graph theory.

\begin{definition}
    For any subgraph or set of subgraphs $\mathcal{G}$ of  $\Gamma_{s}$,
    the vertex set of $\mathcal{G}$ is noted  $V_{\mathcal{G}}$, and the number of vertices in $V_{\mathcal{G}}$ is noted  $n_{V_{\mathcal{G}}}$.
\end{definition}

We now introduce the definition of a \emph{strongly connected component} \cite{bang2008digraphs}.

\begin{definition}
    \label{SCC_chap3}

    A \emph{strongly connected component (SCC)} of a  graph $\mathcal{G}$ is  a maximal  subgraph of $\mathcal{G}$ which is irreducible, i.e.,  no additional arcs or vertices from $\mathcal{G}$ can be included in the subgraph without breaking its property of being irreducible.
\end{definition}

Assume  $\Gamma_s$ admits $N$  SCCs, where $1 \le N \le n$.
The case $N = 1$ corresponds to the situation where $D^s$  is irreducible, and $N = n$  corresponds to the case where $\Gamma_s$   is composed of $n$  SCCs, each composed of a single patch (in other words, the system \eqref{sterile_model}  consists of $n$ independent subsystems).  It is a classical result that one may arrange the SCCs in an {\em acyclic ordering}, that is, into a sequence  $\mathcal{G}_1, \ldots, \mathcal{G}_N$, so that there is no arc (no direct path)  from $\mathcal{G}_j$ to $\mathcal{G}_i$ unless  possibly if $j <i$ \cite[p.~17]{bang2008digraphs}.

\begin{definition}\label{indegree}For any   SCC $\mathcal{G}$ of  $\Gamma_s$, the set of its \emph{in-neighbouring} (resp. \emph{out-neighbouring}) \emph{SCCs}, noted $\mathcal{G}^{-}$ (resp. $\mathcal{G}^{+}$), is the set of other SCCs from which there exists a direct path leading to $\mathcal{G}$ (resp. that are reachable from $\mathcal{G}$ via a direct path),
    that is, for any SCC $\mathcal{H} \neq \mathcal{G}$:
    $$ \mathcal{H}  \in  \mathcal{G}^{-} \Leftrightarrow \mbox{ there exists } i \in V_{\mathcal{G}} \mbox{ and }  j \in V_{\mathcal{H}} \mbox{ such that }  d^s_{ij} >0. $$
    $$( \mbox{resp. }\mathcal{H}  \in  \mathcal{G}^{+} \Leftrightarrow \mbox{ there exists } j \in V_{\mathcal{G}} \mbox{ and }  i \in V_{\mathcal{H}} \mbox{ such that }  d^s_{ij} >0. )$$
    The \emph{in-degree}  of  $\mathcal{G}$ is the number of vertices in $\mathcal{G}^-$. In particular, $\mathcal{G}$ is of in-degree zero if $\mathcal{G}^- = \emptyset$.
\end{definition}
Notice that  an SCC $\mathcal{H}$ cannot belong simultaneously to   $\mathcal{G}^-$ and $\mathcal{G}^+$,  and that $\mathcal{H} \in \mathcal{G}^-$ if and only if $\mathcal{G} \in \mathcal{H}^+$.

Let us introduce the concept of upstream (resp.~downstream) subgraphs.

\begin{definition}  For any SCC $\mathcal{G}$ of  $\Gamma_s$,
    the set of its \emph{upstream} (resp.~\emph{downstream}) \emph{SCCs}, noted $\mathcal{G}^{up}$ (resp. $\mathcal{G}^{down}$), is the set of other SCCs from which there exists a path leading to $\mathcal{G}$ (resp.~that are reachable from $\mathcal{G}$ via a path).
\end{definition}

  Intuitively, sterile insects can only move from an upstream SCC to a downstream SCC, or remain within the same SCC.

With these considerations, one may now introduce the following theorem, which refines \Cref{eq_steriles} by providing a more precise characterisation of the positivity of the equilibrium point, depending on the patches where sterile insects are released and the network connectivity.

\begin{theorem}
    \label{sterile_equilibrium_reducible}
    The equilibrium $x^*_s(\Lambda)$ satisfies, for any SCC $\mathcal{G}$ of $\Gamma_s$:
    \begin{equation}
        \label{eq_SCC}
        x^*_s(\Lambda)|_{V_{\mathcal{G}}} = -  \left(J^s|_{V_{\mathcal{G}}}\right)^{-1} \left( \Lambda|_{V_{\mathcal{G}}} + \left(\sum\limits_{\substack{j \in V_{\mathcal{G}^-}}}d^s_{ij} x^*_{s,j}(\Lambda) \right)_{i \in V_{\mathcal{G}}}  \right).
    \end{equation}
    In particular, for any $\overline{\Lambda}$-admissible $\Lambda \in \mathbb{R}^n_+$ such that $\Lambda_j > 0$ for every $j \in \mathcal{C}^s$,
    \begin{equation}
        \label{eq_sterile_red}
        x^*_s(\Lambda)|_{V_{\mathcal{G}}} =
        \begin{cases}
            \gg 0_{n_{V_{\mathcal{G}}}} & \text{if } \{V_{\mathcal{G}} \cup  V_{\mathcal{G}^{\text{up}}}\} \cap \mathcal{C}^s \neq \emptyset, \\
            0_{n_{V_{\mathcal{G}}}}     & \text{otherwise}.
        \end{cases}
    \end{equation}
\end{theorem}

\begin{proof}
    See \Cref{proof_sterile_reducible}.
\end{proof}

\Cref{sterile_equilibrium_reducible} is  a corollary of the second point in  \Cref{eq_steriles}.  Indeed, \Cref{sterile_equilibrium_reducible} can be  applied to the particular case where
$D^s$ is irreducible. In such case,  $\Gamma_s$ is composed of a single SCC, which is $\Gamma_s$ itself.

\begin{remark}
    \label{Remark4_1}
     As a consequence of \Cref{sterile_equilibrium_reducible} and its proof, sterile insects disperse (and persist) throughout the entire network if and only if they are released in each SCC with in-degree zero (see \Cref{indegree} for the definition of an in-degree). 
\end{remark}

\begin{figure}[H]
    \centering
    \includegraphics[scale = 0.68]{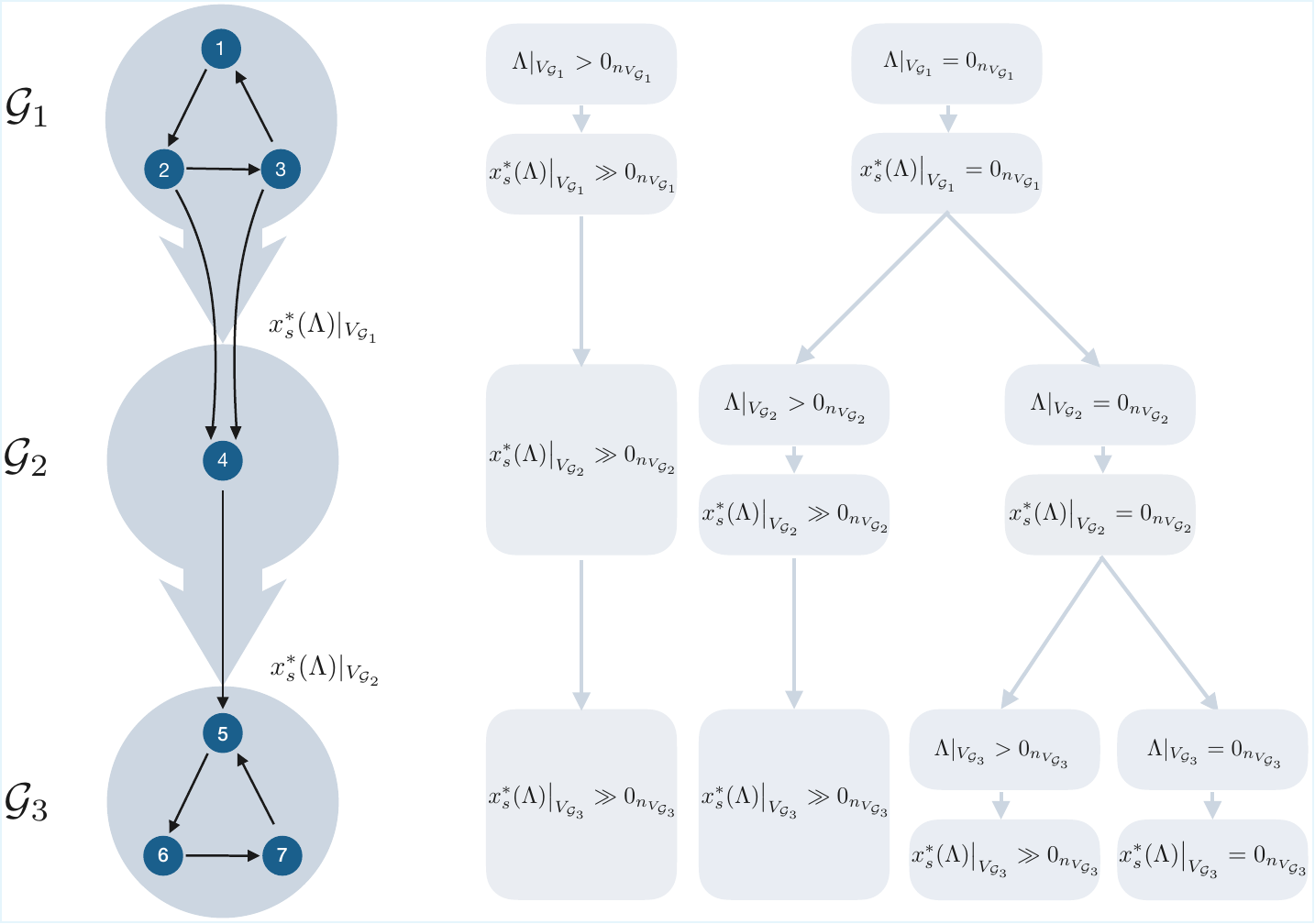}
    \caption{ \textbf{Determination of the sterile insect  equilibrium in the case of reducible  graph.}
        The  graph (left) is composed of three SCCs, namely $\mathcal{G}_1$, $\mathcal{G}_2$ and $\mathcal{G}_3$ in acyclic ordering.  Upstream and downstream relations between SCCs are  given by this ordering.
        The right figure depicts the various scenarios possible for the sterile insect equilibrium value.}
    \label{fig_reducible}
\end{figure}

An example of reducible  graph is illustrated in \Cref{fig_reducible}.

\subsection{Supremal sterile insect equilibrium}
\label{Supremal sterile insect equilibrium}

 We now evaluate the supremal sterile insect equilibrium  $x^\infty_s$.  For each admissible release vector $\Lambda$, we have shown in \Cref{eq_steriles} that   system \eqref{sterile_model} admits a unique equilibrium, which increases as $\Lambda$ increases.  The supremal sterile insect equilibrium is the  upper bound of these equilibria over all $\overline{\Lambda}$-admissible release vectors $\Lambda$, and is obtained when each component of $\Lambda$ is set to its maximal admissible value:
\begin{equation}
    \label{x_s_res}
    x^\infty_s := \sup\{x^*_s(\Lambda): \Lambda ~ \overline{\Lambda}\text{-admissible} \} = \lim_{\Lambda \rightarrow \overline{\Lambda}} x^*_s(\Lambda).
\end{equation}
 The supremum is defined componentwise.   The quantity $x^\infty_s$ will be useful for assessing the feasibility of wild insect elimination through SIT in \Cref{introduction_ control_chap3}, under constraints on the release rate $\Lambda$.

When $\mathcal{C}^s_I = \emptyset$, every component of the vector $\overline{\Lambda}$ is finite,  which implies as a consequence of \Cref{eq_steriles} that $x^\infty_s$ is simply equal to $- \left(J^s\right)^{-1} \overline{\Lambda}$. When $\mathcal{C}^s_I \neq \emptyset$ and following the same approach as in \Cref{Asymptotic behaviour of the sterile insects}, we now characterise $x_s^\infty$ depending on the connectivity matrix $D^s$.  In  what follows, we first consider  in \Cref{Irreducible sterile insects network}  the case where $D^s$ is irreducible. Building on this result, we then address in \Cref{Reducible sterile insects network} the general case,  without any assumption on the irreducibility of $D^s$.

\subsubsection{Irreducible sterile insect network}
\label{Irreducible sterile insects network}

\begin{lemma}
    \label{lemma_supremal_sterile_eq_irred}
    Assume $\mathcal{C}^s_I \neq \emptyset$ and  $D^s$ is irreducible. Then, the supremal sterile insect equilibrium $x^\infty_s$ satisfies
    $$  x^\infty_{s,i} =  +  \infty, \quad i = 1,\ldots,n.$$
\end{lemma}

\begin{proof}
    By  \Cref{eq_steriles}, for any $\overline{\Lambda}$-admissible $\Lambda$, one has
    $$x^*_s(\Lambda) = - \left(J^s\right)^{-1} \Lambda.$$
    As seen before in \Cref{Asymptotic behaviour of the sterile insects}, the matrix $\left(J^s\right)^{-1}$ is strictly negative when $D^s$ is irreducible.  Therefore, if  $\Lambda_k \to +\infty$ for some  $k \in \{1,\ldots,n\}$, then $x^*_{s,i}(\Lambda) \to \infty$ for all $i =1,\ldots,n$. Since $\overline{\Lambda}_i = + \infty$ for any $i \in \mathcal{C}^s_I$, one has
    $$ \lim_{\Lambda \rightarrow \overline{\Lambda}} x^*_{s,i}(\Lambda) = +\infty, \quad i = 1,\ldots,n.$$
    This completes the proof of \Cref{lemma_supremal_sterile_eq_irred}.
\end{proof}

\subsubsection{Reducible sterile insect network}
\label{Reducible sterile insects network}

\begin{lemma}
    \label{lemma_supremal_sterile_eq}
    For any SCC $\mathcal{G}$ of $\Gamma_s$ and for any $k \in V_{\mathcal{G}}$, the supremal sterile insect equilibrium $x^\infty_s$ satisfies
    $$x^\infty_{s,k} = \left\{
        \begin{array}{ll}
            +\infty                                                                                                                                                                                                                                           & \mbox{if } \{ V_{\mathcal{G}} \cup  V_{\mathcal{G}^\text{up}}  \} \cap \mathcal{C}^s_I \neq \emptyset \\
            0                                                                                                                                                                                                                                                 & \mbox{if }  \{ V_{\mathcal{G}} \cup  V_{\mathcal{G}^\text{up}}  \} \cap \mathcal{C}^s = \emptyset     \\
            -\left( \left(J^s|_{V_{\mathcal{G}}}\right)^{-1}\left[\overline{\Lambda}|_{V_{\mathcal{G}}}  + \left(\sum\limits_{\substack{j \in V_{\mathcal{G}^-}}}{d^s_{ij} x^\infty_{s,j} } \right)_{i \in V_{\mathcal{G}}} \right]\right)_k \in (0, +\infty) & \mbox{otherwise}.
        \end{array}
        \right.$$
\end{lemma}

\begin{proof}
    See \Cref{proof_lemma_sup}.
\end{proof}

\begin{remark} 
    As a consequence of \Cref{lemma_supremal_sterile_eq} and \Cref{Remark4_1}, the sterile insect equilibrium can be made arbitrarily large in all patches if and only if the releases can be made arbitrarily large in each SCC with in-degree zero (see \Cref{indegree} for the definition of an in-degree).
\end{remark}

Let us summarise the main results of this section. We showed in \Cref{eq_steriles} that the solutions of the sterile insect system \eqref{sterile_model} converge to a unique equilibrium. The result in  \Cref{sterile_equilibrium_reducible} characterises in which patches the components of this equilibrium are positive, depending on the network structure and the release patches. Last, we also assessed in \Cref{lemma_supremal_sterile_eq} how large this equilibrium can be under constraints on the release rates. These results will be useful to assess the feasibility of wild insect elimination through SIT in Section 5.

\section{General properties of the  SIT model}

\label{introduction_ control_chap3}

The aim of this section is to adjust the release rate $\Lambda$ in system \eqref{sterile_model}  to drive the wild population towards elimination.  This section is essentially an application of \Cref{n_patch} since, as will be seen,
     the release of sterile insects artificially induces or amplifies a  strong Allee effect in some patches, through an increase in the Allee parameter $a$.

In \Cref{Controlled model and assumptions},  the SIT model which consists in the coupling of the wild and sterile populations is presented. Next, in \Cref{controlproperties}, we establish  important  properties of monotonicity with respect to  the release rate. A  sufficient  condition for  the feasibility of elimination of the wild population under some constraints on the release rates is then  derived in \Cref{limitbehaviourforlargeaction_chap3}.
Finally, in \Cref{Section time entrance basin}, assuming this condition holds  and that the origin of the uncontrolled system \eqref{log_allee} is LAS, we provide an upper bound  on the minimal release duration necessary   to drive the wild population into the basin of attraction of the origin of the uncontrolled system \eqref{log_allee}.

\subsection{SIT model}
\label{Controlled model and assumptions}

When sterile insects are released, they primarily disrupt the mating process.  This implies that the birth
rate of the wild population has to be modified to take into account the releases of sterile insects, as
done in many other SIT models like \cite{bliman2019implementation,cai2014dynamics,yang2020dynamics}. Consequently, in each patch $i$, the term $\frac{x_i}{x_i + a_i}$ is replaced by $\frac{x_i}{x_i + a_i + \gamma  x_{s,i}}$. This term represents the probability of a fertile mating occurring in Patch $i$, and $\gamma$ is a parameter reflecting the average relative competitiveness of the sterile insects.

By coupling the wild population system \eqref{log_allee} with the sterile insect system \eqref{sterile_model}, we obtain the following SIT model  of $2n$ ODEs:
\begin{subequations}
    \label{coupled_model_control}
    \begin{align}
        \dot{x}_i     & =  x_i \big(   p(x_i, a_i + \gamma x_{s,i})  b_i   - \mu_{1,i} - \rho_i -\mu_{2,i} x_i\big) + \sum\limits_{\substack{j = 1 \\ j \neq i}}^{n} d_{ij} x_j -  \sum\limits_{\substack{j = 1\\ j \neq i}}^{n} d_{ji} x_i, \label{sterile_a}
        \\
        \dot{x}_{s,i} & = \Lambda_i - \left(\mu_{s,i} +  \rho_{s,i}\right) x_{s,i}+ \sum\limits_{\substack{j = 1                                         \\ i \neq j}}^{n}d^s_{ij} x_{s,j} -  \sum\limits_{\substack{j = 1 \\ i \neq j}}^{n}d^s_{ji} x_{s,i}, \qquad \qquad \qquad i = 1,\ldots,n, \label{sterile_b}
    \end{align}
\end{subequations}
where  the function $p$ has been defined in  \eqref{chap3_eq1}.

     Using the formula in \eqref{solution_sterile_control} for the solution $x(t,x_{s,0})$ of \eqref{sterile_b} with sterile insect initial condition $x_{s,0}$, the wild insect population  satisfies the $n$-dimensional system of time-varying ODEs:
\begin{equation}
    \label{log_allee_control_time}
    \dot{x}_i =  x_i \Big(  p\big(x_i, a_i + \gamma  x_{s,i}(t,x_{s,0})\big)  b_i  - \mu_{1,i} - \rho_i -\mu_{2,i} x_i\Big) + \sum\limits_{\substack{j = 1\\ j \ne i}}^{n} d_{ij} x_j -  \sum\limits_{\substack{j = 1\\ j \ne i}}^{n} d_{ji} x_i, \qquad i = 1,\ldots,n.
\end{equation}

 Since $x_s$ converges exponentially towards the equilibrium $x^*_s(\Lambda)$ by \Cref{eq_steriles}, one has  by \cite[Theorem 3.1]{vidyasagar2003decomposition} that the long-term behaviour of  \eqref{log_allee_control_time} is equivalent to that of the  following $n$-dimensional autonomous system:
\begin{equation}
    \label{log_allee_control}
    \dot{x}_i =  x_i \Big(  p\big(x_i, a_i + \gamma x^*_{s,i}(\Lambda)\big)  b_i  - \mu_{1,i} - \rho_i -\mu_{2,i} x_i\Big) + \sum\limits_{\substack{j = 1\\ j \ne i}}^{n} d_{ij} x_j -  \sum\limits_{\substack{j = 1\\ j \ne i}}^{n} d_{ji} x_i, \qquad i = 1,\ldots,n.
\end{equation}

 A central point here is that, by applying  \Cref{0only} with the substitution $a \mapsto a+ \gamma x^*_s(\Lambda)$, it follows that   \eqref{log_allee_control} admits a maximal equilibrium, noted $x_+(a+\gamma x^*_s(\Lambda))$, which is GAS on the set $\{x \geq x_+(a+\gamma x^*_s(\Lambda))\}.$

\Cref{convergence_allee} also applies to system \eqref{log_allee_control}, yielding  the following  result, for the matrix function $A$ defined in \eqref{A_a}.

\begin{theorem}
    \label{convergence_allee_control}
    If $s\big(A(a + \gamma x^*_s(\Lambda))\big)  < 0$, then $x_+(a+\gamma x^*_s(\Lambda)) = 0_n$, and the origin  of  \eqref{log_allee_control} is GAS on $\mathbb{R}^n_+$.
\end{theorem}

\subsection{Monotonicity properties}
\label{controlproperties}

We now introduce some  monotonicity properties  relative to the release rate $\Lambda$.

First, since   $x^*_s(\Lambda)$  is a linear and increasing function of $\Lambda$  by \Cref{eq_steriles}, the following result is derived from \Cref{max_eq_decreasing}.

\begin{theorem}
    \label{x(A)decreases}
    The maximal equilibrium $x_+(a+\gamma x^*_s(\Lambda))$ is a non-increasing  function of $\Lambda$, that is, for any $\overline{\Lambda}$-admissible $\Lambda, \Lambda' \in \mathbb{R}^n_+$,
    \begin{equation}
        \label{u<u'}
        \Lambda < \Lambda' \Rightarrow x_+(a+\gamma x^*_s(\Lambda'))\le x_+(a+\gamma x^*_s(\Lambda)).
    \end{equation}

    Last, if $D$ is irreducible and $x_+(a+\gamma x^*_s(\Lambda))> 0_n$,  the inequality  $\le$ in \eqref{u<u'} may be replaced  by the inequality  $\ll$.
\end{theorem}

  \Cref{x(A)decreases} implies that when $D$ is irreducible,  increasing the release rate in any single patch reduces the maximal wild population in every patch of the network.

Similarly, the next result    follows directly from \Cref{corollary_max_eq} and \Cref{s(A(a))_decreasing}.

\begin{theorem}
    \label{s(A)_decreasing}
    If the origin of \eqref{log_allee_control} is GAS on $\mathbb{R}^n_+$, then the origin of system \eqref{log_allee_control} with release rate  $\Lambda'$ is GAS on $\mathbb{R}^n_+$ for all  $\Lambda' \geq \Lambda$.

    Moreover, the function  $\Lambda \in \mathbb{R}^n_+ \mapsto s\big(A(a + \gamma x^*_s(\Lambda))\big) $ is convex and  non-increasing in the following sense:
    \begin{equation}
        \label{non-decreasing A_a(u)}
        \Lambda < \Lambda' \Rightarrow s\big(A(a + \gamma x^*_s(\Lambda'))\big)  \le s\big(A(a + \gamma x^*_s(\Lambda))\big) .
    \end{equation}

    Last, if $D$ is irreducible,   then  this function is  continuous on $\mathbb{R}^n_+$, and
    continuously differentiable and piecewise twice continuously differentiable on the set  $\{\Lambda \in \mathbb{R}^n_+ :  a + \gamma x^*_s(\Lambda) \gg  0_n \}$. If in addition $a_i + \gamma x^*_{s,i}(\Lambda) < \frac{b_i}{\mu_{2,i}}$ for all $i = 1,\ldots,n$, the inequality $\le$ in  \eqref{non-decreasing A_a(u)}
    may be  replaced by the inequality $<$.
\end{theorem}

\begin{proof} 
    The first part of the proof is a direct application of \Cref{corollary_max_eq}.

    We only prove the convexity of the function $\Lambda \mapsto s\big(A(a + \gamma x^*_s(\Lambda))\big) $, since the remainder of \Cref{s(A)_decreasing} follows by composing  $a \mapsto s(A(a))$ with the increasing function  $\Lambda \mapsto  a + \gamma x^*_s(\Lambda)$ and applying  \Cref{s(A(a))_decreasing}.

    The function $\Lambda \mapsto a+ \gamma x^*_s(\Lambda)$ is affine and thus convex. Since the function $a \mapsto s(A(a))$ is convex and non-decreasing  by \Cref{s(A(a))_decreasing}, we deduce by composition that   the function $\Lambda \in \mathbb{R}^n_+ \mapsto s\big(A(a + \gamma x^*_s(\Lambda))\big) $ is convex.
\end{proof}

\subsection{A sufficient condition for wild population elimination feasibility}
\label{limitbehaviourforlargeaction_chap3}

In this section, we provide a sufficient   condition under which, for given $\overline{\Lambda} \in \overline{\mathbb{R}}^n_+$,  the elimination of the wild population may be achieved for some $\overline{\Lambda}$-admissible  release rates $\Lambda \in \mathbb{R}^n_+$.

     Since the stability modulus is a non-decreasing function of the  release rate (\Cref{s(A)_decreasing}), our analysis focuses on evaluating  the stability modulus of  $A(a + \gamma x^*_s(\Lambda))$ for large $\overline{\Lambda}$-admissible  $\Lambda$. To this aim, denote
        $$
            A(a + \gamma x^\infty_s)  = \lim_{\Lambda \to \overline{\Lambda}} A(a + \gamma x^*_s(\Lambda)),$$
        where  $x^\infty_s$  is the limit of $x^*_s(\Lambda)$ as $\Lambda \to \overline{\Lambda}$. Recall that an explicit formula for   $x^\infty_s$  is provided in \Cref{lemma_supremal_sterile_eq}.

        The quantity  $A(a + \gamma x^\infty_s)$ is obtained from the matrix $A(a)$ in \eqref{A_a}, by modifying only the local diagonal terms. For indices $i$ such that $x^\infty_{s,i} < + \infty$, the substitution $a_i \mapsto a_i + \gamma x^\infty_{s,i}$ is applied in the corresponding diagonal terms. For indices $i$ such that $x^\infty_{s,i} = + \infty$, the diagonal local terms reduce to $-(\mu_{1,i}+\rho_i)$, reflecting the absence of local growth.

\begin{lemma}
    \label{s(A)_chap3}
    For the matrix $A$ defined in \eqref{A_a},
    one has
    \begin{align*}
        \inf \left\{s\big(A(a + \gamma x^*_s(\Lambda))\big): \Lambda~ \overline{\Lambda}\mbox{-admissible} \right\} & = s\big(A(a + \gamma x^\infty_s)\big).
    \end{align*}
\end{lemma}

\begin{proof}
    Since $s\big(A(a + \gamma x^*_s(\Lambda))\big) $ is a non-increasing function of $\overline{\Lambda}$-admissible $\Lambda$ by \Cref{s(A)_decreasing}, it follows that
    $$\inf \left\{s\big(A(a + \gamma x^*_s(\Lambda))\big) : \Lambda~ \overline{\Lambda}\mbox{-admissible} \right\}  = \lim_{\Lambda \rightarrow \overline{\Lambda}} s\big(A(a + \gamma x^*_s(\Lambda))\big) .$$
    By \eqref{x_s_res}, one has that
    $\lim_{\Lambda \rightarrow \overline{\Lambda}} x^*_s(\Lambda) = x^\infty_s$, which implies
    $$\lim_{\Lambda \rightarrow \overline{\Lambda}} A(a+\gamma x^*_s(\Lambda)) = A(a + \gamma x^\infty_s).$$
    The equality in \Cref{s(A)_chap3} then follows from the continuity of the stability modulus of a Metzler matrix.
\end{proof}

Building on  \Cref{s(A)_chap3}, we derive a sufficient condition for  the feasibility of population elimination. In addition, recall that $s(J_a) < 0$, for $J_a$ defined in \eqref{J_a}, implies the local asymptotic stability of the origin of  the uncontrolled model \eqref{log_allee} (see \Cref{0LAS}). Thus, we also introduce the time required for a trajectory of the  model  with prescribed releases \eqref{log_allee_control_time}  to enter  $\mathcal{B}_a$, the basin of attraction of the origin of the uncontrolled system \eqref{log_allee}.

\begin{definition}
    \label{def temps inf}
    For any wild and sterile insect initial conditions $x_0, x_{s,0}\in \mathbb{R}^n_+$, let  $\tau(x_0,x_{s,0})$ denote the \emph{release duration}, that is,
    $$
        \tau(x_0,x_{s,0})
        := \inf\{t \geq 0 : x(t,x_0, x_{s,0}) \in \mathcal{B}_a \},$$
    where   $x(t,x_0,x_{s,0})$  is the solution of the  system   with prescribed releases \eqref{log_allee_control_time}, with wild and sterile insect initial conditions respectively equal to $x_0$ and $x_{s,0}$, and $\mathcal{B}_a$ is the basin of attraction of the uncontrolled system \eqref{log_allee}.

     In addition, denote by $\Sigma_s(x_0, x_{s,0})$ the \emph{total number of sterile insects released until entering the uncontrolled basin of attraction of the origin}, that is, for prescribed release rate vector $\Lambda$,
    $$\Sigma_s(x_0, x_{s,0}) := \mathbf{1}\t_n \Lambda ~ \tau(x_0,x_{s,0}).$$
\end{definition}
Departing from $x_0$, stopping the release of sterile insects for $t\geq \tau(x_0,x_{s,0})$ slows down the elimination, but does not suppress it.

The next result provides a sufficient condition for the feasibility  of the elimination  of the wild insect population within a finite time. 

\begin{theorem}
    \label{s(A)><0_chap3}
     If  $s\big(A(a + \gamma x^\infty_s)\big) < 0$, then there exists $\Lambda^*$ $\overline{\Lambda}$-admissible such that $x_+(a+\gamma x^*_s(\Lambda))= 0_n $ for any  $\overline{\Lambda}$-admissible $\Lambda \geq \Lambda^*$.

    Moreover if $s(J_a) < 0$, then $  \tau(x_0,x_{s,0})  < + \infty$  and $\Sigma_s(x_0,x_{s,0}) < + \infty$ for any $\Lambda \geq \Lambda^*$   and any $x_0,x_{s,0}\in \mathbb{R}^n_+$.

\end{theorem}

\begin{proof}
    If $s\big(A(a + \gamma x^\infty_s)\big) < 0$, then by continuity of the stability modulus, 
    there exists a  $\overline{\Lambda}$-admissible $\Lambda^*$ such that $s\big(A(a + \gamma x^*_s(\Lambda))\big) < 0$ for any  $\overline{\Lambda}$-admissible  $\Lambda \geq \Lambda^*$. By applying \Cref{convergence_allee_control}, the origin of \eqref{log_allee_control} is then GAS for
    for any $\Lambda \geq \Lambda^*$, and  $x_+(a+\gamma x^*_s(\Lambda))= 0_n$.

    If moreover $s(J_a)<0$, then the origin
    of the uncontrolled system \eqref{log_allee} is LAS by \Cref{0LAS}.  Since $x_+(a+\gamma x^*_s(\Lambda))= 0_n$, there exists a time $t$ such that $x(t,x_0,x_{s,0})$ enters the basin of attraction $\mathcal{B}_a$ of the uncontrolled system \eqref{log_allee},  implying $ \tau(x_0,x_{s,0}) < + \infty$,  and thus $\Sigma_s(x_0, x_{s,0}) < + \infty$.
\end{proof}

\Cref{s(A)><0_chap3} states in particular that if elimination is feasible and  $s(J_a)<0$,  the   release duration  $\tau(x_0,x_{s,0})$ is  finite:
the releases
can be stopped as soon as the solution of \eqref{log_allee_control} enters the basin of attraction of the origin of the uncontrolled system \eqref{log_allee}. In contrast,  as a consequence of \Cref{elimination_model_uncontrolled_chapter3}, when $a = 0_n$  and $s(J_a)>0$, that is, the extinction equilibrium of \eqref{log_allee} is unstable, the release must be maintained indefinitely  to reach elimination  in all patches,  causing the total number $\Sigma_s(x_0,x_{s,0})$ of sterile insects released to increase linearly with time.  However, notice that using non constant release values, e.g., values elaborated via a feedback law based on measurement of the wild population, may lead to stabilisation with a finite total release effort (see, e.g., \cite{bhaya2025feedback}).

Moreover, notice that when $\mathcal{C}^s_I \neq \emptyset$ and $D^s$ is irreducible, elimination is always feasible. Indeed in such case, we have by \Cref{lemma_supremal_sterile_eq_irred} that  all components of $ x^\infty_s$ are infinite. This implies that
$$A(a + \gamma x^\infty_s) = -\diag(\mu_{1,i}+\rho_i) + D.$$
This matrix is Metzler and strictly column diagonally dominant,  and is therefore Hurwitz.

\subsection{Estimate of the release duration}

\label{Section time entrance basin}

In the previous section, we showed that when $s(J_a) < 0$ and $s(A( a + \gamma x^\infty_s) ) < 0$, elimination of  wild insects may be achieved in finite time by  large enough  $\overline{\Lambda}$-admissible release rates $\Lambda$. In the following, building on the estimation of the basin of attraction of the origin of the uncontrolled system \eqref{log_allee} in \Cref{estimation basin of attraction}, we give an upper bound for the release  duration $\tau(x_0,x_{s,0})$.

\begin{theorem}
    \label{time entrance basin of attraction}
    Let  $a > 0_n$. Assume $D$ is irreducible and that $s(J_a) < 0$, and let $c \gg 0_n$ be a corresponding left  (Perron) eigenvector of $J_a$,  normalised so that $\min_i c_i = 1$.
    Also assume that $s\big(A(a + \gamma x^\infty_s)\big)< 0$, and let $\Lambda$ $\overline{\Lambda}$-admissible such that  $\Lambda \geq \Lambda^*$, for $\Lambda^*$ defined in \Cref{s(A)><0_chap3}.   Then, for any $x_0, x_{s,0}\in \mathbb{R}^n_+$,
    \begin{equation*}
        \tau(x_0, x_{s,0}) \le \inf\left\{t \geq 0 : \Psi\left(c^T x(t,x_0,x_{s,0})\right) < -s(J_a)  \right\}  < + \infty,
    \end{equation*}
    where   $x(t,x_0,x_{s,0})$  is the solution of the  system  with prescribed releases \eqref{log_allee_control_time}, and the function $\Psi$ has been defined in \Cref{estimation basin of attraction}.

    In particular, when $x_{s,0}= 0_n$,
    \begin{equation}
        \label{time entrance x_s,0 = 0_n }
        \tau(x_0, 0_n) \le \inf\left\{t \geq 0 : \Psi\left(c^T x(t,x_0,0_n)\right) < -s(J_a)  \right\}  < + \infty.
    \end{equation}

\end{theorem}

To apply this result, one must simulate the non-autonomous differential equation  \eqref{log_allee_control_time}, and compute the evolution of the scalar $\Psi(c\t x(t))$ along the resulting trajectory.

\begin{proof}[Proof of \Cref{time entrance basin of attraction}]
    We first recall the result of \Cref{estimation basin of attraction}:
    \begin{equation*}
        \{ x \in \mathbb{R}^n_+ : \Psi( c\t x) < -s(J_a) \} \subset \mathcal{B}_a,
    \end{equation*}
    where $\mathcal{B}_a$ is the basin of attraction in $\mathbb{R}^n_+$ of the origin of the uncontrolled system \eqref{log_allee}. Therefore, one has
    $$\tau(x_0,x_{s,0}) \le \inf\left\{t \geq 0 : \Psi\left(c^T x(t,x_0,x_{s,0})\right) < -s(J_a)  \right\}.$$

    Finally, we  prove that

    \begin{equation}
        \label{inf finite}
        \inf\left\{t \geq 0 : \Psi\left(c^T x(t,x_0,x_{s,0})\right) < -s(J_a)  \right\} < + \infty.
    \end{equation}
    As $\Lambda \geq \Lambda^*$ for $\Lambda^*$ defined in \Cref{s(A)><0_chap3}, one has thanks to that theorem that
    $$\lim_{t \to +\infty }c^T x(t,x_0,x_{s,0}) = 0.$$
    Since $-s(J_a) > 0$ and $\Psi(0) = 0$, we deduce by continuity of $\Psi$ that \eqref{inf finite} holds.
    This concludes the proof of the statement.
\end{proof}

Notice that, by cooperativeness of  system \eqref{log_allee_control_time}, and since   the solution $x_s(t,x_{s,0})$ of \eqref{sterile_b} satisfies  $x_s(t,x_{s,0}) \geq x_s(t,0_n)$, the
upper bound provided in \eqref{time entrance x_s,0 = 0_n } when the sterile insect initial condition $x_{s,0}$ is null is valid for all $x_{s,0} \in \mathbb{R}^n_+$.

Moreover, still by monotonicity of the  system \eqref{log_allee_control_time}, for all viable states, that is, for all wild insect initial condition $x_0 \in \{x: 0_n \le x <x_+(a) \}$, where $x_+(a)$ is the maximal equilibrium of the uncontrolled system \eqref{log_allee}, one has in particular
$$\tau(x_0,x_{s,0}) \le \inf\left\{t \geq 0 : \Psi\left(c^T x(t,x_+(a),0_n)\right) < -s(J_a)  \right\}  < + \infty.$$

Also notice that the maximal equilibrium $x_+(a)$ can be easily  computed numerically.
When $a=0_n$, by \Cref{elimination_model_uncontrolled_chapter3}, the maximal equilibrium
$x_+(0)$ is GAS on $\mathbb{R}^n_+\setminus\{0_n\}$. Thus, for any initial condition $x_0>0_n$ the solution of \eqref{log_allee} with $a=0_n$ converges to $x_+(0)$.
For $a>0_n$, \Cref{0only} shows that $x_+(a)$ is GAS on the set $\{x\ge x_+(a)\}$. Moreover, \eqref{decrease_eq} implies
$$
    x_+(a)\le x_+(0)\quad\text{for all }a>0_n.
$$
Therefore, initialising the numerical solution of \eqref{log_allee} at $x_+(0)$ yields convergence to $x_+(a)$.

\vspace{0.5cm}

    We now summarise the main results of this section. We derived a sufficient condition for the elimination of the wild population through SIT in \Cref{convergence_allee_control}, and established some monotonicity properties with respect to the release rate in \Cref{x(A)decreases} and \Cref{s(A)_decreasing}. This led to a sufficient condition for the feasibility of wild population elimination under constraints on the release rates in \Cref{s(A)><0_chap3}. Last, when   the origin of the uncontrolled system \eqref{log_allee} is LAS, we provided a finite upper bound for the release duration $\tau(x_0,x_{s,0})$ in \Cref{time entrance basin of attraction}.

\section{A cost-effective approach for achievable elimination scenarios}

\label{minimization_section}

In this section, it is  assumed that $s\big(A(a + \gamma x^\infty_s)\big)<0$, so that the wild population can be driven to elimination   for sufficiently large $\overline{\Lambda}$-admissible release rates
        (see \Cref{s(A)><0_chap3}). We build a cost-efficient strategy  to achieve   elimination,  while minimising a certain control cost. Here, this cost is taken as a weighted sum of the  sterile insect  release rates $\Lambda_i$,  but  the results  presented below apply equally to the minimisation of any other convex, increasing and coercive function of $\Lambda$. Notice that  a related optimisation problem consisting  of optimising over  additional mortality rates  has  been studied in \cite{bliman2024feasibility}.

We first introduce the minimisation problem and its solvability properties (\Cref{minsol}) and address a monotonicity property of the solutions with respect to the Allee parameter $a$ (\Cref{f_a_decrease}).  Then, in order to solve numerically the problem, we provide  in   \Cref{gradprobleme} a formula for the gradient of the  problem’s constraint, which relies on an explicit expression for  the gradient of the stability modulus of an irreducible Metzler matrix with respect to its  diagonal elements.
\newline

Let us  introduce a  fixed vector $\pi$, where for each \( i \in \mathcal{C}^s \), \( \pi_i > 0 \) denotes the relative intervention price for the release of sterile insects in Patch \( i \).

For any  $\alpha$ such that $0 < \alpha < - s\big(A(a + \gamma x^\infty_s)\big) $, let us introduce the following minimisation problem:
\begin{equation}
    \label{minimizationproblem_chap3}
    \begin{tabular}{llll}
        \text{minimise}   & $\pi\t \Lambda$                          \\
        \text{subject to} & $\Lambda ~\overline{\Lambda}$-admissible \\
                          & $h(\Lambda) \le -\alpha$,
    \end{tabular}
\end{equation}
where  the constraint in \eqref{minimizationproblem_chap3} is defined by the function:
\begin{equation}
    \label{h}
    h(\Lambda) := s\big(A(a + \gamma x^*_s(\Lambda))\big) .
\end{equation}
The quantity $\pi\t \Lambda$ represents the total daily  cost corresponding to the $\overline{\Lambda}$-admissible release rate $\Lambda \in \mathbb{R}^n_+$. Since $\alpha>0$,  the constraint $h(\Lambda) \le - \alpha$ ensures  global asymptotic stability of the origin of \eqref{log_allee_control} (see \Cref{convergence_allee_control}).

The solvability of problem~\eqref{minimizationproblem_chap3}  is addressed in the next result. The proof is exactly the same as for Theorem 5.1 in \cite{bliman2024feasibility}  and thus is not provided here. It
relies on the monotonicity and convexity properties of the function $h$ (\Cref{s(A)_decreasing}).

\begin{theorem}
    \label{minsol}
    Let $\alpha$ satisfy $ s\big(A(a + \gamma x^\infty_s)\big)< -\alpha  < 0 $.
    Then Problem \eqref{minimizationproblem_chap3} admits local minimisers. The latter  are also global and their set is convex. Moreover, the following holds:
    \begin{enumerate}
        \item If $s(A(a)) > -\alpha$, then any minimiser $\Lambda^*$ satisfies $\Lambda^* > 0_n$ and $h(\Lambda^*)=-\alpha$.
        \item If $s(A(a)) \le -\alpha$, then the minimiser is unique and equal to $0_n$.
    \end{enumerate}
\end{theorem}

\begin{theorem}
    \label{f_a_decrease}
    Let $a,a' \in \mathbb{R}^n_+$ and let $\alpha$ satisfy $ s\big(A(a + \gamma x^\infty_s)\big)< -\alpha  < 0 $. Let  $\Lambda_{a'}^*, \Lambda_a^*$ be  corresponding minimisers of Problem \eqref{minimizationproblem_chap3}.
    Then,
    $$
        a \leq a' \Rightarrow \pi\t \Lambda^*_{a'} \le \pi\t \Lambda^*_{a}.
    $$
\end{theorem}

\begin{proof}
    For any $a' \geq a$, we have by \Cref{s(A(a))_decreasing} that $s\big(A(a' + \gamma x^\infty_s)\big) \le s\big(A(a + \gamma x^\infty_s)\big) <0$. Hence, by \Cref{minsol}, Problem~\eqref{minimizationproblem_chap3} with parameters $a$ and $a'$ admits global minimisers.

    Let $\Lambda^*_a$ and $\Lambda^*_{a'}$ be minimisers of Problem \eqref{minimizationproblem_chap3} with parameters $a$ and $a'$, respectively. Since $a' \geq a$, one has,
    still by \Cref{s(A(a))_decreasing}, that
    $$s\big(A(a' + \gamma x^*_s(\Lambda^*_a))\big)  \le s( A(a + \gamma x^*_s(\Lambda^*_a)))= -\alpha,$$
    where the last equality is a consequence of \Cref{minsol}. Therefore, $\Lambda^*_a$ satisfies the constraints of Problem \eqref{minimizationproblem_chap3} with parameter $a'$. By definition of a minimiser, one deduces that
    $ \pi\t \Lambda^*_{a'} \le \pi\t \Lambda^*_a.$
\end{proof}

By \Cref{minsol}, \eqref{minimizationproblem_chap3} is a convex minimisation problem.  It can be solved efficiently by an interior-point method \cite[Chapter 11]{boyd2004convex}, slightly adapted from \cite{bliman2024feasibility}, where the optimisation was performed with respect to the additional mortality rates  $\rho_i$.  Interior-point algorithms are particularly well-suited for convex problems with constraints, as they approach the optimal solution from the interior of the feasible region while ensuring constraint satisfaction at each iteration.

To apply this algorithm, it is necessary to compute the gradient of the function $h$. For this, one takes advantage of the concept of \textit{group inverse}. The reader is referred to \cite{kirkland2012group} for a thorough presentation of this notion.

\begin{definition}
    Assume that $B \in \mathbb{C}^{n \times n}$ is a complex singular matrix, and that its eigenvalue $0$ is semisimple, that is, its algebraic and geometric multiplicities coincide. Then, the \emph{group inverse} of $B$, denoted $B^{\#}$, is the unique matrix $X \in \mathbb{C}^{n \times n}$ such that
    $$(i)~ BXB = B;\qquad
        (ii)~ XBX = X;\qquad
        (iii)~ BX = XB.
    $$
\end{definition}

The following technical result shows how to compute the gradient $\nabla h(\Lambda)$  of the function $h(\Lambda)$.
\begin{lemma}
    \label{gradprobleme}
    Assume  $D$ is irreducible. Moreover, for any $\Lambda \in \mathbb{R}^n_+$, let $Q(\Lambda) := h(\Lambda)I - A(a+\gamma x^*_s(\Lambda))$, and note $Q^\#(\Lambda)$ its group inverse.  Then, for all $\Lambda \gg 0_n$ and $i=1,\ldots,n$,
     \begin{equation}
        \label{grad_h}
        \big(\nabla h(\Lambda)\big)_i = \gamma  \sum\limits_{p = 1}^n  \mu_{2,p}  \left(J^s \right)^{-1}_{pi}  \left( \sqrt{\frac{b_p }{\big(a_p + \gamma x^*_{s,p}(\Lambda)\big) \mu_{2,p}}} - 1\right)^+
        \left(I- Q(\Lambda) Q^\#(\Lambda) \right)_{pp}.
    \end{equation}
\end{lemma}

\begin{proof}
    See  \Cref{Gradient evaluation}.
\end{proof}

To compute the group inverse $Q^{\#}(\Lambda)$  of the matrix $Q(\Lambda)$, one uses the formula provided in \cite[Remark 2.5.3]{kirkland2012group}. This formula is stated therein for matrices  written in the form $s(A)I-A$,  for $A$ a non-negative and irreducible matrix.
One shows easily that it remains valid for an irreducible Metzler matrix $A$.

In \Cref{gradprobleme}, the gradient of $ h(\Lambda)$ is given for $\Lambda \gg 0_n$,  which ensures  $x^*_s(\Lambda) \gg 0_n$ by \Cref{sterile_equilibrium_reducible}. Since $a + \gamma x^*_s(\Lambda) \gg 0_n$,  $h$ is then differentiable with respect to $\Lambda$ by  \Cref{s(A)_decreasing}.  We are only interested in strictly positive  values of $\Lambda_i$ for $i \in \mathcal{C}^s$, since  the variables corresponding to  indices in $\{1,\ldots,n\}\setminus \mathcal{C}^s$
are of course discarded and admissible points for the interior-point algorithm must strictly satisfy the inequality constraints.

\section{Numerical exploration of different scenarios}
\label{simulations}

    In this section,
    we consider the SIT control of the oriental fruit fly, \textit{Bactrocera dorsalis}, by exploring different scenarios. Typical values for local parameters are provided in  \Cref{Table_parametres_general}, page \pageref{Table_parametres_general}.

In \Cref{Table_parametres_general}, the competitiveness parameter $\gamma$ is  set to $0.39$, reflecting that sterile insects are generally less competitive than wild insects in mating \cite{dyck2021sterile,dumont2025improvement}. Notice that a higher $\gamma$ reduces the number of insects that need to be released. Some treatments can  enhance the competitiveness of the  sterile insects. For example, when releasing sterile males of \textit{Bactrocera dorsalis}, the male  competitiveness parameter can be increased up to $1$, notably through treatment with methyl eugenol before the releases \cite{ji2013effect}.

Regarding the dispersal parameters $d_{ij}$, no data appear to be available in the literature. However, it is likely that mark–release–recapture experiments could provide such estimates \cite{jang2017mark,shelly2010mark}. We   assume here that the dispersal rates of sterile and wild insects are identical, i.e.,
$$D^s =  D.$$

\begin{table}[H]
    \centering
    \begin{tabular}{lll ll}
        \toprule
        \textbf{Symbol} & \textbf{Description }                         & \textbf{Value} & \textbf{Unit}        & \textbf{Reference}           \\
        \midrule
        $b_i$           & Adult birth rate                              & $6.60$         & day$^{-1}$           & \cite{ekesi2006field}        \\
        $\mu_{1,i}$     & Wild  insect death rate                       & 1/80.75        & day$^{-1}$           & \cite{ekesi2006field}        \\
        $\rho_i$        & Wild  insect additional death rate with MT    & 0.5/80.75      & day$^{-1}$           & Chosen                       \\
        $\mu_{2,i}$     & Wild insect density death rate                & 0.001          & Ind$^{-1}$day$^{-1}$ & Chosen                       \\
        $\mu_{s,i}$     & Sterile insect death rate                     & 2/86.4         & day$^{-1}$           & \cite{dumont2025improvement} \\
        $\rho_{s,i}$    & Sterile  insect additional death rate with MT & 1/86.4         & day$^{-1}$           & Chosen                       \\
        $\gamma$        & Sterile insect competitive parameter          &  0.39    & --                   & \cite{dumont2025improvement} \\
        \bottomrule
    \end{tabular}
    \caption{Parameter values for \textit{Bactrocera dorsalis}.  When mass trapping (MT) is applied, it is assumed to increase the natural mortality rates of both wild and sterile insects by 50\%.}
    \label{Table_parametres_general}
\end{table}

In this section, we display  solutions of Problem \eqref{minimizationproblem_chap3} in different configurations,      assuming that the release rates can be arbitrarily large in the release patches, that is, $\mathcal{C}^s_I = \mathcal{C}^s$, with the notations introduced in \Cref{def C^s}. Throughout this exploration, we take  the  price vector to be $\pi= 1_n$, i.e., the total number of sterile insects to be released daily to eliminate the wild population is minimised.  We  also set  $\alpha = 10^{-5}$ days$^{-1}$. The optimal release rate in each release patch is rounded up to the nearest integer.

We also  assume guaranteed natural mating success in each patch ($a = 0_n$), as no empirical data are currently available.  Under this worst-case scenario,   by  \Cref{remarkyves}, the resulting optimal release rates ensure elimination regardless of the actual value of $a$, providing robustness to the results with respect to uncertainties on the value of $a$. The impact of the Allee  parameter $a$ on the  release duration  $\tau(x_0,x_{s,0})$, defined in \Cref{def temps inf}, is studied in \Cref{section allee}.

A three-patch example is first studied (\Cref{Three-patch model}), followed by a more complex seven-patch network (\Cref{seven patch model}).

\subsection{Three-patch model}

\label{Three-patch model}

In this section, we assume that the network is composed of three patches ($n=3$), as illustrated in \Cref{Fig three-patch system}.

\begin{figure}[htbp]
    \centering
    \begin{tikzpicture}[auto,
            cloud/.style={
                    circle,
                    draw,
                    thick,
                    minimum size=1.5cm, 
                    inner sep=0pt, 
                    fill=bleu_fonce,
                    text=white,           
                    font=\bfseries\large
                },
            line/.style={draw, -latex', thick}
        ]

        \def\L{3.4} 
        \node[cloud]  (1) at (0,0)                        {$1$};
        \node[cloud] (2) at (\L,0)                      {$3$};
        \node[cloud]   (3) at ({\L/2},{\L*sqrt(3)/2})     {$2$};

        \path[line] (1) -- node[midway, above] {$d_{31}$} (2);
        \path[line] (2) -- node[midway, left]  {$d_{23}$} (3);
        \path[line] (3) -- node[midway, right] {$d_{12}$} (1);
        \path[line, bend left] (2) to node[midway, below] {$d_{13}$} (1);
        \path[line, bend left] (1) to node[midway, left]  {$d_{21}$} (3);
        \path[line, bend left] (3) to node[midway, right] {$d_{32}$} (2);

        \path[line, looseness=2, out=210, in=150, min distance=1cm]
        (1) to node[left] {$g_1(x_1)$} (1);

        \path[line, looseness=2, out=330, in=30, min distance=1cm]
        (2) to node[right] {$g_3(x_3)$} (2);

        \path[line, looseness=2, out=60, in=120, min distance=1cm]
        (3) to node[above] {$g_2(x_2)$} (3);

    \end{tikzpicture}
    \caption{3-patch flow diagram.}
    \label{Fig three-patch system}
\end{figure}

For parameters identical in every patch and symmetric dispersal, 
the natural persistence equilibrium  ${x_+(0)}_i$  remains the same regardless of the dispersal coefficients \cite{elbetch2021multi}. Thus, its value coincides with the persistence equilibrium obtained in the absence of dispersal, given by \Cref{logistique_a=0}:
\begin{equation}
    \label{natural_equilibrium}
    {x_+(0)}_i = 6587.62, \quad i = 1,\ldots,n.
\end{equation}

\subsubsection{Effects of dispersal rates and density-dependent mortality rates}

In this section, the impact of dispersal parameters   and  density-dependent mortality rates on the SIT strategy is investigated, using the parameter values provided in \Cref{Table_parametres_general}, page \pageref{Table_parametres_general}, and identical dispersal rates between connected  patches.

\paragraph{Interconnected patches}

The SIT strategies associated with different release patches $\mathcal{C}^s$ are presented in \Cref{Table1} and \Cref{Table2}, page \pageref{Table1}, for dispersal coefficients all  equal  to 0.02 and 0.01, respectively.
In these tables, the last row represents the optimal total number of sterile insects to release daily to ensure elimination of the wild population.

Since identical dispersal rates is assumed  between all patches and the local parameters are identical in every patch, the outcome depends only on the number of release patches.

Even if the wild insect persistence equilibrium is the same regardless of the dispersal coefficients,  by comparing \Cref{Table1} and \Cref{Table2},  the  critical  release rate  increases as the dispersal parameters decrease, except when the releases occur in all patches. In that case, both the allocation and the total number of insects released per day remain unchanged, regardless of the dispersal coefficients.  As expected, the best strategy is to release the sterile insects in all three patches.

However, while our simulations show all optimal releases strategies, it is important to discuss their usefulness in the field. For instance, for practical purposes, it may be better to release in one patch only: the amount of sterile insects to release is larger but it requires only one team in the field to carry out the releases, except e.g. if the three patches are sufficiently close to be visited the same day.

\begin{table}[H]
    \centering
    \small

    \begin{minipage}{\textwidth}
        \centering
        \begin{tabular}{|>{\centering\arraybackslash}m{2cm}||c|c|c|c|c|c|c|}
            \hline
            $\mathcal{C}^s$                              & \{1\}        & \{2\}        & \{3\}        & \{1,2\}      & \{1,3\}      & \{2,3\}      & \{1,2,3\}     \\
            \hline
            $\left(\Lambda_i^*\right)$                   & (1468,--,--) & (--,1468,--) & (--,--,1468) & (713,713,--) & (713,--,713) & (--,713,713) & (359,359,359) \\
            \hline
            $\rule{0pt}{2.3ex} \sum_{i=1}^3 \Lambda_i^*$ & 1468         & 1468         & 1468         & 1426         & 1426         & 1426         & 1077          \\
            \hline
        \end{tabular}
        \caption{Estimates of the critical (total) release rate and the optimal release strategy for the 3-patch SIT control with $d_{ij}= d^s_{ij}= 0.02, i,j = 1,2,3, i \neq j$. }
        \label{Table1}
    \end{minipage}

    \vspace{1em}


    \begin{minipage}{\textwidth}
        \centering
        \begin{tabular}{|>{\centering\arraybackslash}m{2cm}||c|c|c|c|c|c|c|}
            \hline
            $\mathcal{C}^s$                              & \{1\}        & \{2\}        & \{3\}        & \{1,2\}      & \{1,3\}      & \{2,3\}      & \{1,2,3\}     \\
            \hline
            $\left(\Lambda_i^*\right)$                   & (1881,--,--) & (--,1881,--) & (--,--,1881) & (921,921,--) & (921,--,921) & (--,921,921) & (359,359,359) \\
            \hline
            $\rule{0pt}{2.3ex} \sum_{i=1}^3 \Lambda_i^*$ & 1881         & 1881         & 1881         & 1842         & 1842         & 1842         & 1077          \\
            \hline
        \end{tabular}
        \caption{Estimates of the (total) critical release rate and the optimal release strategy for the 3-patch SIT control with $d_{ij}=d^s_{ij}=0.01, i,j = 1,2,3, i \neq j$.}
        \label{Table2}

    \end{minipage}
\end{table}

We now consider the case where one patch   (say, Patch 2) has a lower density-dependent mortality rate ($\mu_{2,2} = 0.0005$), corresponding to a higher carrying capacity. The  natural persistence equilibria are  reported in \Cref{Table_nat_eq_1}, page \pageref{Table_nat_eq_1}  and the SIT strategies in \Cref{Table3} and \Cref{Table4}, page \pageref{Table4}, for dispersal coefficients equal to $0.02$ and $0.01$, respectively.
As expected, the natural equilibrium as well as the number of sterile insects to be released daily are larger than in the case of patches with
equal density-dependent mortality rates (compare with \Cref{Table1} and \Cref{Table2}),  since the total carrying capacity is larger.  However, unlike the  equal density-dependent mortality rate scenario---where the persistence  equilibrium is independent of dispersal---the equilibrium here
varies with dispersal rates.

A  patch  with a lower density-dependent mortality rate also changes the allocation strategy. For any number of release patches, the best strategy is always to release more sterile insects in the  patch with lower  density-dependent mortality rate  (Patch 2). This is due to the fact that Patch $2$ has  the largest pest population at equilibrium.

If releases occur only in  Patch 2, the corresponding release rate remains the same as in the equal   density-dependent mortality rate case,
resulting in the  same amount of sterile insects to be released daily (compare the second columns of Tables \ref{Table3} and \ref{Table4}  with Tables \ref{Table1} and \ref{Table2}).
Moreover, according to \Cref{Table3} and \Cref{Table4}, releasing in two patches, say $\{1,2\}$ or $\{2,3\}$,  yields nearly the same outcome as releasing in Patch 2  alone. Therefore when
a patch  has a lower density-dependent mortality rate  and the two others are identical, the best strategy —if the releases cannot occur in all  patches——is to release only in the  patch  with lower density-dependent mortality rate, which  avoids deploying two teams without affecting  significantly the  release rate.

As in the  equal density-dependent mortality rate scenario, the required  release rate increases as dispersal decreases, except when all patches are treated. In that case, the total  release rate remains  nearly unchanged, as in the case of   patches with equal  density-dependent mortality rate; however, the allocation now varies with dispersal.

\begin{table}[H]
    \centering
    \begin{tabular}{|c|c|}
        \hline
        Dispersal rate between all patches & Natural persistence equilibrium \\
        \hline
        0.02                               & (6607.38, 13135.47,  6607.38)   \\
        \hline
        0.01                               & (6597.56, 13155.29,  6597.56)   \\
        \hline
    \end{tabular}
    \caption{ 3-patch model with $\mu_{2,2} = 0.0005$: persistence equilibrium of the  wild  population versus  dispersal rates.}
    \label{Table_nat_eq_1}
\end{table}

\begin{table}[H]
    \centering
    \small

    \begin{minipage}{\textwidth}
        \centering

        \begin{tabular}{|>{\centering\arraybackslash}m{2cm}||c|c|c|c|c|c|c|}
            \hline
            $\mathcal{C}^s$                               & \{1\}        & \{2\}        & \{3\}        & \{1,2\}      & \{1,3\}        & \{2,3\}      & \{1,2,3\}    \\
            \hline
            $\left(\Lambda^*_i\right)$                    & (2850,--,--) & (--,1469,--) & (--,--,2850) & (93,1356,--) & (1425,--,1425) & (--,1356,93) & (70,1291,70) \\
            \hline
            $ \rule{0pt}{2.3ex} \sum_{i=1}^3 \Lambda_i^*$ & 2850         & 1469         & 2850         & 1449         & 2850           & 1449         & 1431         \\
            \hline
        \end{tabular}
        \caption{Estimates of the critical (total) release rate and the optimal release strategy for the 3-patch  SIT control with  $d_{ij}= d^s_{ij}=0.02, i,j = 1,2,3, i \neq j$ and $\mu_{2,2} = 0.0005$.}
        \label{Table3}

    \end{minipage}

    \vspace{1em}

    \begin{minipage}{\textwidth}
        \centering
        \begin{tabular}{|>{\centering\arraybackslash}m{1.5cm}||c|c|c|c|c|c|c|}
            \hline
            $\mathcal{C}^s$                              & \{1\}        & \{2\}        & \{3\}        & \{1,2\}       & \{1,3\}        & \{2,3\}       & \{1,2,3\}      \\
            \hline
            $\left(\Lambda^*_i\right)$                   & (3683,--,--) & (--,1881,--) & (--,--,3683) & (838,1004,--) & (1842,--,1842) & (--,1004,838) & (213,1007,213) \\
            \hline
            $\rule{0pt}{2.3ex} \sum_{i=1}^3 \Lambda_i^*$ & 3683         & 1881         & 3683         & 1842          & 3684           & 1842          & 1433           \\
            \hline
        \end{tabular}
        \caption{Estimates of the critical (total) release rate and the optimal release strategy for the 3-patch  SIT control with  $d_{ij}=d^s_{ij}=0.01, i,j = 1,2,3, i \neq j$ and $\mu_{2,2} = 0.0005$.}
        \label{Table4}
    \end{minipage}
\end{table}

\begin{remark}
    The carrying capacity affects the release strategy here, whereas in our previous work  \cite{bliman2024feasibility},  where a logistic-type local growth function was considered, it did not influence the optimal additional mortality rate induced by mass trapping (MT).  This is because the release of sterile insects induces a strong Allee effect in the model.
\end{remark}

\paragraph{Chain structure}

In this section, we assume that $d_{13} = d_{31} = 0$, so that the network has a chain structure, as illustrated in \Cref{fig:three_patch_chain}.

\begin{figure}[htbp]
    \centering
    \begin{tikzpicture}[auto,
            cloud/.style={
                    circle,
                    draw,
                    thick,
                    minimum size=1.5cm, 
                    inner sep=0pt,      
                    fill=bleu_fonce,
                    text=white,
                    font=\bfseries\large
                },
            line/.style={draw, -latex', thick}
        ]

        \def\L{3.4}
        \node[cloud] (1) at (0,0)     {$1$};
        \node[cloud] (2) at (\L,0)    {$2$};
        \node[cloud] (3) at (2*\L,0)  {$3$};

        \path[line, bend left=25] (1) edge node[above] {$d_{21}$} (2);
        \path[line, bend left=25] (2) edge node[below] {$d_{12}$} (1);

        \path[line, bend left=25] (2) edge node[above] {$d_{32}$} (3);
        \path[line, bend left=25] (3) edge node[below] {$d_{23}$} (2);

        \path[line, looseness=2, out=60, in=120, min distance=1cm]
        (3) to node[above] {$g_3(x_3)$} (3);

        \path[line, looseness=2, out=60, in=120, min distance=1cm]
        (2) to node[above] {$g_2(x_2)$} (2);

        \path[line, looseness=2, out=60, in=120, min distance=1cm]
        (1) to node[above] {$g_1(x_1)$} (1);

    \end{tikzpicture}
    \caption{ 3-patch flow diagram with a chain structure.}
    \label{fig:three_patch_chain}
\end{figure}

When the patches have equal density-dependent mortality rates ($\mu_{2,i}=0.001$, $i = 1,2,3$)  and the dispersal is symmetric, the persistence equilibrium of the wild insect population coincides  with  that  in the absence of dispersal (see \eqref{natural_equilibrium}).
The SIT strategies are listed in \Cref{Table5} and \Cref{Table6}, page \pageref{Table6}, for the dispersal coefficients equal to $0.02$ and $0.01$, respectively.  If releases are allowed in all three patches, in Patches $\{1,3\}$, or in Patch $2$---cases that reflect the model symmetry---the SIT strategy remains unchanged compared to the fully connected scenario  (see \Cref{Table1} and \Cref{Table2}, page \pageref{Table2}). In the remaining configurations, the strategy is altered and larger releases are required.

Notice that releasing in two patches, or in Patch 2 only, does not substantially affect the number of sterile insects to be released daily.
However,  in the two-patch release case, only the releases in Patches $\{1,3\}$ are realistic. In particular, releasing in Patch 2 alone is nearly equivalent to releasing in Patches $\{1,2\}$ or $\{2,3\}$,  since Patch 2 is directly connected to the two other patches.
     Given the limited gain achieved by releasing in two patches instead of one, it may be sufficient to adopt the slightly suboptimal   strategy  of releasing in Patch 2 only, thereby avoiding the need to deploy two teams in the field.

\begin{table}[H]
    \centering
    \small

    \begin{minipage}{\textwidth}

        \begin{tabular}{|>{\centering\arraybackslash}m{2cm}||c|c|c|c|c|c|c|}
            \hline
            $\mathcal{C}^s$                              & \{1\}        & \{2\}        & \{3\}        & \{1,2\}      & \{1,3\}      & \{2,3\}      & \{1,2,3\}     \\
            \hline
            $\left(\Lambda^*_i\right)$                   & (3117,--,--) & (--,1468,--) & (--,--,3117) & (16,1448,--) & (713,--,713) & (--,1448,16) & (359,359,359) \\
            \hline
            $\rule{0pt}{2.3ex} \sum_{i=1}^3 \Lambda_i^*$ & 3117         & 1468         & 3117         & 1464         & 1426         & 1464         & 1077          \\
            \hline
        \end{tabular}
        \caption{Estimates of the critical (total) release rate and the optimal release strategy for the 3-patch SIT control when the network has a chain structure, with
            $d_{12} = d_{21} = d_{23} = d_{32} =0.02$.}
        \label{Table5}
    \end{minipage}

    \vspace{1em}

    \begin{minipage}{\textwidth}
        \begin{tabular}{|>{\centering\arraybackslash}m{2cm}||c|c|c|c|c|c|c|}
            \hline
            $\mathcal{C}^s$                              & \{1\}        & \{2\}        & \{3\}        & \{1,2\}     & \{1,3\}      & \{2,3\}     & \{1,2,3\}     \\
            \hline
            $\left(\Lambda^*_i\right)$                   & (6185,--,--) & (--,1881,--) & (--,--,6185) & (8,1870,--) & (921,--,921) & (--,1870,8) & (359,359,359) \\
            \hline
            $\rule{0pt}{2.3ex} \sum_{i=1}^3 \Lambda_i^*$ & 6185         & 1881         & 6185         & 1878        & 1842         & 1878        & 1077          \\
            \hline
        \end{tabular}
        \caption{Estimates of the critical (total) release rate and the optimal release strategy for the 3-patch SIT control when the network has a chain structure, with
            $d_{12} = d_{21} = d_{23} = d_{32} =0.01$.}
        \label{Table6}
    \end{minipage}
\end{table}

When Patch 2  has a lower density-dependent mortality rate $(\mu_{2,2} = 0.0005)$, the natural persistence equilibrium of the wild population is still the same as in the fully connected case (see \Cref{Table_nat_eq_1},  page \pageref{Table_nat_eq_1}). The SIT strategies for dispersal coefficients of 0.02 and 0.01 between connected patches are reported in \Cref{Table7} and \Cref{Table8},  page \pageref{Table8}, respectively. A comparison of \Cref{Table7} (\Cref{Table8}) with \Cref{Table5} (\Cref{Table6}) shows that, unlike in the fully connected case,   reducing the  density-dependent mortality rate in  Patch 2 only alters the SIT strategy when releases  occur in  both Patches $1$ and $3$.

Notice that when releases can occur in at most two patches,  as in the  case with identical patches,  the best strategy is  to release only in the  patch  with lower density-dependent mortality rate (Patch 2), especially when dispersal increases.

\begin{table}[H]
    \small
    \centering

    \begin{minipage}{\textwidth}
        \begin{tabular}{|>{\centering\arraybackslash}m{2cm}||c|c|c|c|c|c|c|}
            \hline
            $\mathcal{C}^s$                              & \{1\}        & \{2\}        & \{3\}        & \{1,2\}      & \{1,3\}        & \{2,3\}      & \{1,2,3\}    \\
            \hline
            $\left(\Lambda^*_i\right)$                   & (3123,--,--) & (--,1469,--) & (--,--,3123) & (14,1452,--) & (1425,--,1425) & (--,1452,14) & (70,1291,70) \\
            \hline
            $\rule{0pt}{2.3ex} \sum_{i=1}^3 \Lambda_i^*$ & 3123         & 1469         & 3123         & 1466         & 2850           & 1466         & 1431         \\
            \hline
        \end{tabular}
        \caption{Estimates of the critical (total) release rate and the optimal release strategy for the 3-patch SIT control when the network has a chain structure, with
            $d_{12} = d_{21} = d_{23} = d_{32} =0.02$ and $\mu_{2,2} = 0.0005$.}
        \label{Table7}
    \end{minipage}

    \vspace{1em}

    \begin{minipage}{\textwidth}
        \begin{tabular}{|>{\centering\arraybackslash}m{2cm}||c|c|c|c|c|c|c|}
            \hline
            $\mathcal{C}^s$                              & \{1\}        & \{2\}        & \{3\}        & \{1,2\}     & \{1,3\}        & \{2,3\}     & \{1,2,3\}      \\
            \hline
            $\left(\Lambda^*_i\right)$                   & (6185,--,--) & (--,1881,--) & (--,--,6185) & (8,1870,--) & (1842,--,1842) & (--,1870,8) & (213,1007,213) \\
            \hline
            $\rule{0pt}{2.3ex} \sum_{i=1}^3 \Lambda_i^*$ & 6185         & 1881         & 6185         & 1878        & 3684           & 1878        & 1433           \\
            \hline
        \end{tabular}
        \caption{Estimates of the critical (total) release rate and the optimal release strategy for the 3-patch SIT control when the network has a chain structure, with
            $d_{12} = d_{21} = d_{23} = d_{32} =0.01$ and $\mu_{2,2} = 0.0005$.}
        \label{Table8}
    \end{minipage}
\end{table}

To conclude this section, 
when releasing in all patches is  not feasible, the   best strategy is to focus mainly on the  patch  with lower density-dependent mortality rate.  In particular, when this  patch is well connected to the others, then releasing  only in this single patch may be the best strategy. Indeed, releasing in two patches does not substantially decrease the total  critical release rate, 
but requires the deployment of two teams on the ground, increasing the cost of the SIT  programme.

\subsubsection{Combining SIT with mass trapping (MT)}

\label{combining SIT with mass trapping}

In this section, the effect of  mass trapping (MT) on the SIT strategy is investigated for a fully connected network, as in \Cref{Fig three-patch system}, page \pageref{Fig three-patch system}. Recall that traps  also capture sterile insects,  and the corresponding additional mortalities are provided in \Cref{Table_parametres_general}. Each dispersal coefficient is set to 0.02.

In practice, producers use MT against \textit{Bactrocera dorsalis}. In the following numerical simulations, traps are assumed to be removed in the patches where sterile insects are released, because MT can negatively impact the sterile insects. However, trap removal can be costly, especially when large areas are involved.

The release strategies in \Cref{Table1}, page \pageref{Table1}, obtained without MT, will serve as reference scenario. The corresponding strategies with MT in the patches where sterile insects are not released  are shown in \Cref{Table_SIT_MT_patch2}, page \pageref{Table_SIT_MT_patch2}.  Comparing this table with \Cref{Table1} shows that trap deployment  requires higher  release rates to achieve elimination.    When releasing in a single patch, applying MT increases the  critical release rate by   35.5\% (1989 vs.~1468). When releasing in two patches, applying MT increases it by  24.0\% (1768 vs 1426).  Such  increases  may pose challenges when sterile insect production is limited.

\begin{table}[H]
    \small
    \centering
    \begin{tabular}{|>{\centering\arraybackslash}m{1.7cm}||c|c|c|c|c|c|c|}
        \hline
        $\mathcal{C}^s$                              & \{1\}        & \{2\}        & \{3\}        & \{1,2\}      & \{1,3\}      & \{2,3\}      \\
        \hline
        $\left(\Lambda^*_i\right)$                   & (1989,--,--) & (--,1989,--) & (--,--,1989) & (884,884,--) & (884,--,884) & (--,884,884) \\
        \hline
        $\rule{0pt}{2.3ex} \sum_{i=1}^3 \Lambda_i^*$ & 1989         & 1989         & 1989         & 1768         & 1768         & 1768         \\
        \hline
    \end{tabular}
    \caption{Estimates of the critical (total) release rate and the optimal release strategy for the 3-patch SIT control depending on the release patches $\mathcal{C}^s$ when   MT  is applied in the remaining patches.
    }
    \label{Table_SIT_MT_patch2}
\end{table}

To conclude this section, regarding release rates, the additional sterile mortality induced by MT outweighs the gain from increased wild mortality. However,  we  show in the next section that,
when a  strong Allee effect is naturally present  in every patch and for fixed release rates, MT may reduce  the release duration  $\timetot$  and  thus the total number  $\nbtot$ of sterile insects released until entering the uncontrolled basin of attraction of the origin.

\subsubsection{Strong Allee effect and sterile insect release duration}

\label{section allee}

In the previous simulations, mating success in each patch was assumed to be guaranteed, i.e., $a = 0_n$, which implies that sterile insects must be released indefinitely \cite{anguelov2020sustainable}; otherwise, the wild population would  return to its persistence equilibrium. When mating conditions are imperfect in each patch, the SIT can  be stopped after a finite  time needed to enter the basin of attraction of the origin of the uncontrolled model \eqref{log_allee}, such that the population will decay till elimination (\Cref{s(A)><0_chap3}).

     In the following simulations, we use as a baseline the optimal critical release rate obtained in the case $a = 0_n$, denoted by $\Lambda_0^*$.  Recall that this release rate guarantees elimination for any value of $a \in \mathbb{R}^n_+$.  Compared to the previous sections, the index $0$ in $\Lambda_0^*$ emphasises that the critical release rates considered here correspond to the case where no  strong Allee effect is present  in any patch. 

        In field applications, mass releases are performed.  By \Cref{s(A)_decreasing} and \Cref{corollary_max_eq}, any release rate $\Lambda$ such that $\Lambda \geq \Lambda^*_0$ drives the wild population to elimination, regardless of $a\in \mathbb{R}^n_+$.
        We express here the release values as multiples of $\Lambda_0^*$, namely  $\Lambda = p \times \Lambda_0^*$ for  $p = 2,5,10$ and $30$. For each value of $\Lambda$, we examine how the Allee parameter $a$ impacts the release duration  $\timetot$ and the total number $\nbtot$ of sterile insects released until entering the uncontrolled basin of attraction of the origin (see \Cref{def temps inf}). We also analyse how these two indicators vary with different release rates $\Lambda$.

As said in the beginning of \Cref{combining SIT with mass trapping},
numerical simulations are performed with dispersal coefficients set to $0.02$.

\begin{figure}[H]
    \centering
    \begin{subfigure}[b]{0.4965\textwidth}
        \centering
        \includegraphics[width=\textwidth]{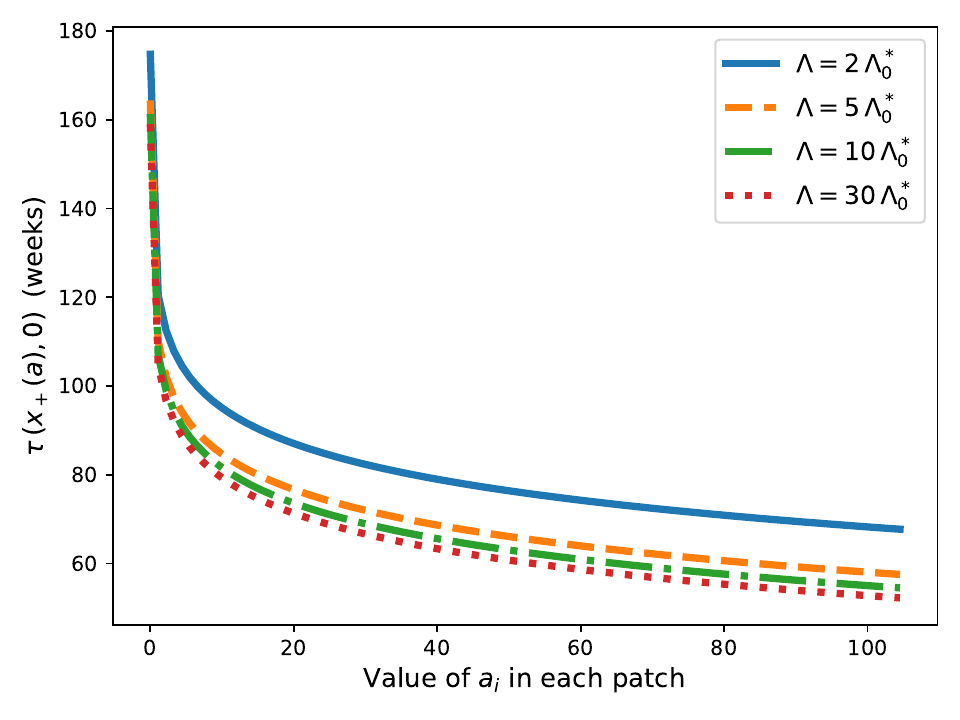}
    \end{subfigure}
    \begin{subfigure}[b]{0.4965\textwidth}
        \centering
        \includegraphics[width=\textwidth]{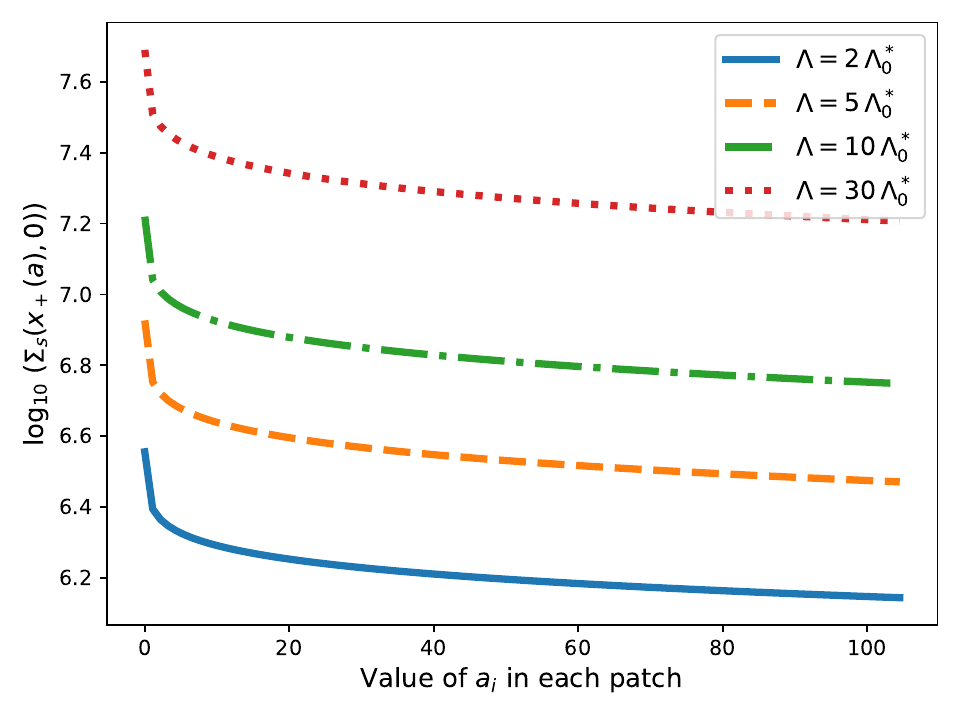}
    \end{subfigure}
    \caption{$3$-patch SIT control with $\mathcal{C}^s=\{2\}$: \textbf{(a)} $\tau(x_+(a),0)$  and  \textbf{(b)} $\Sigma_s(x_+(a),0)$  as functions of the Allee  parameter.  Both $\timetot$ and $\nbtot$ decrease with respect to both the Allee parameter and the release rates.}
    \label{Fig_a_i_1}
\end{figure}

\begin{figure}[H]
    \centering
    \begin{subfigure}[b]{0.4965\textwidth}
        \centering
        \includegraphics[width=\textwidth]{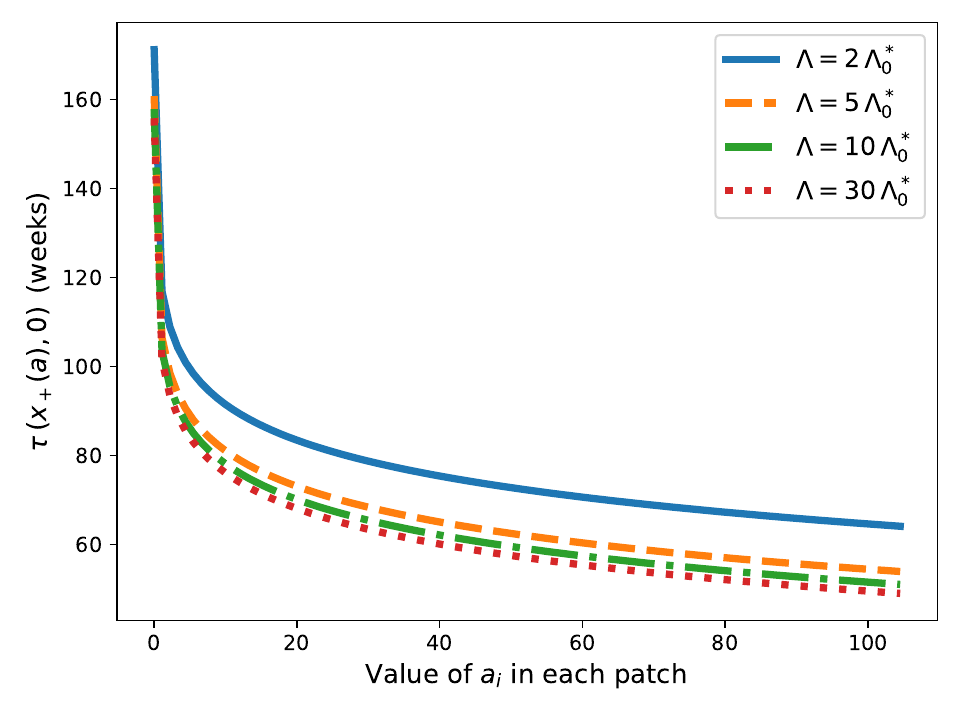}
    \end{subfigure}
    \begin{subfigure}[b]{0.4965\textwidth}
        \centering
        \includegraphics[width=\textwidth]{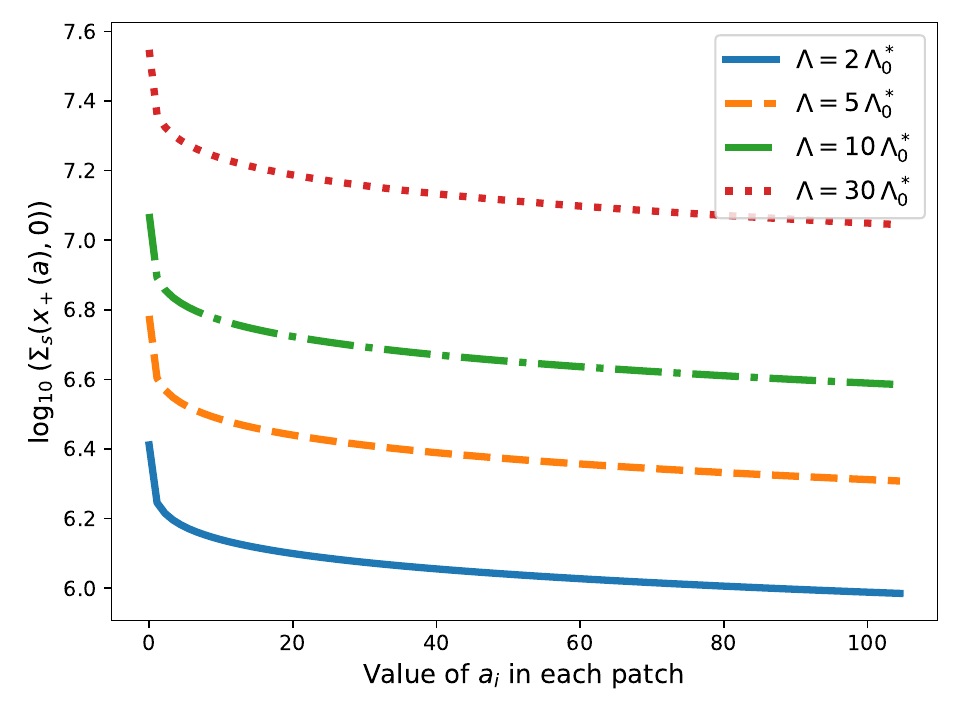}
    \end{subfigure}
    \caption{$3$-patch SIT control with $\mathcal{C}^s=\{1,2,3\}$:  \textbf{(a)} $\tau(x_+(a),0)$  and  \textbf{(b)} $\Sigma_s(x_+(a),0)$  as functions of the Allee  parameter.  Compared with \Cref{Fig_a_i_1}, releasing in all three patches reduces  $\nbtot$ relative to the one-patch release case.}
    \label{Fig_a_i_2}
\end{figure}

        Figures \ref{Fig_a_i_1}(a) and \ref{Fig_a_i_2}(a), page \pageref{Fig_a_i_1}, show that  $\timetot$ decreases rapidly as $a_i$ increases in all patches, and also decreases with respect to $p=2$,  $5$  $10$ and $30$.  However, the release duration curves for
                $p = 5$, $p=10$ and $p=30$ are very  close, and thus  $\nbtot$
                grows dramatically as $p$ increases (see \Cref{Fig_a_i_1}(b) and \Cref{Fig_a_i_2}(b)).
                These observations suggest the existence of a threshold release rate, beyond which further increases provide no practical advantage. The best result for $\nbtot$ is obtained for $\Lambda = 2  \Lambda_0^*$.

        As shown earlier in Tables \ref{Table5}--\ref{Table8}, the critical release rate decreases as the number of release patches increases. Moreover, comparing \Cref{Fig_a_i_1}(a) and \Cref{Fig_a_i_2}(a) shows that the release durations in the one-patch and three-patch release cases are similar. This implies that, due to the lower critical release rates in the three-patch release case,  $\nbtot$
        is significantly reduced when releases are carried out in three patches (compare \Cref{Fig_a_i_1}(b) and \Cref{Fig_a_i_2}(b)).

\begin{figure}[H]
    \centering

    \begin{minipage}{\textwidth}

        \begin{subfigure}[b]{0.495\textwidth}
            \centering
            \includegraphics[width=\textwidth]{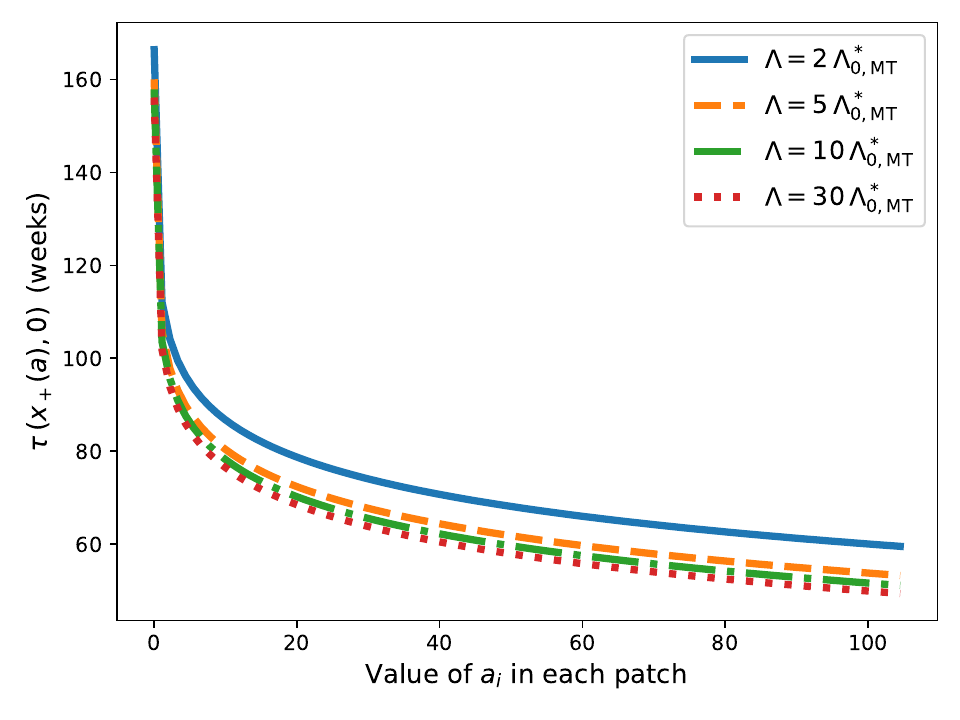}
        \end{subfigure}
        \begin{subfigure}[b]{0.495\textwidth}
            \centering
            \includegraphics[width=\textwidth]{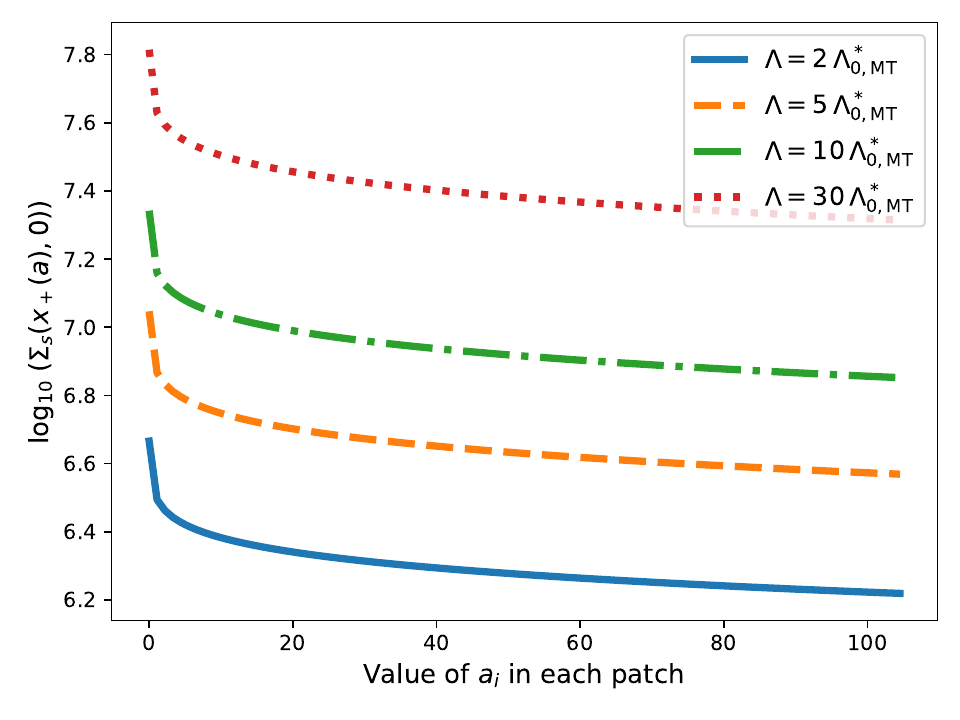}
        \end{subfigure}
        \caption{$3$-patch SIT control with $\mathcal{C}^s=\{2\}$: \textbf{(a)} $\tau(x_+(a),0)$  and  \textbf{(b)} $\Sigma_s(x_+(a),0)$  as functions of the Allee  parameter.   In contrast to \Cref{Fig_a_i_1}, the release rates $\Lambda$ are chosen as multiples of the critical release rate obtained with MT in Patches~$\{1,3\}$.}
        \label{Fig_a_i_3}

    \end{minipage}

    \vspace{1em}

    \begin{minipage}{\textwidth}

        \begin{subfigure}[b]{0.495\textwidth}
            \centering
            \includegraphics[width=\textwidth]{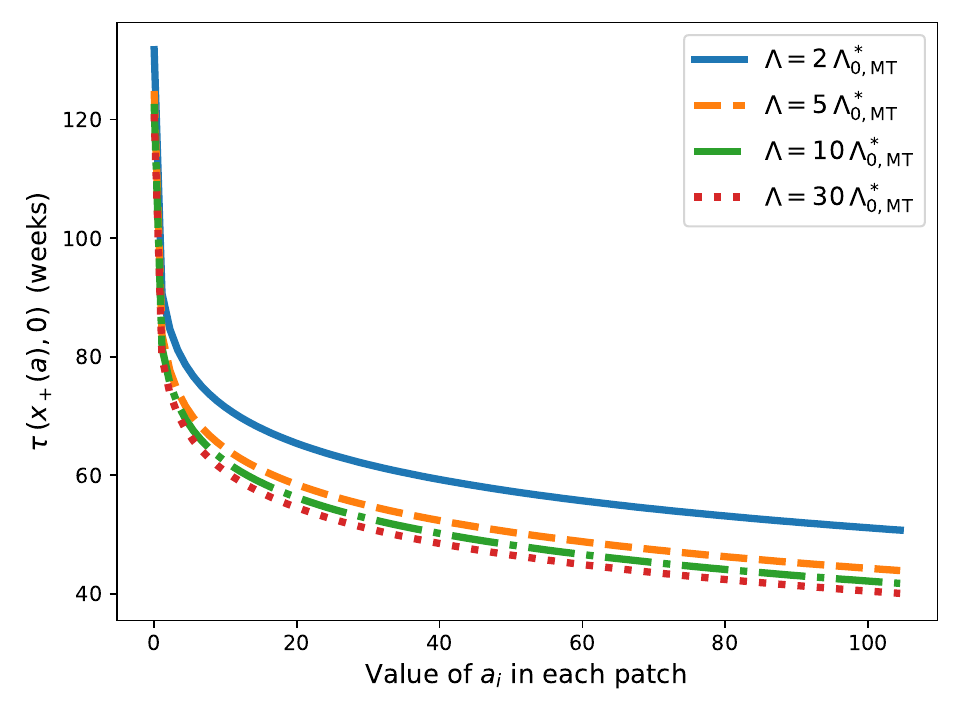}
        \end{subfigure}
        \begin{subfigure}[b]{0.495\textwidth}
            \centering
            \includegraphics[width=\textwidth]{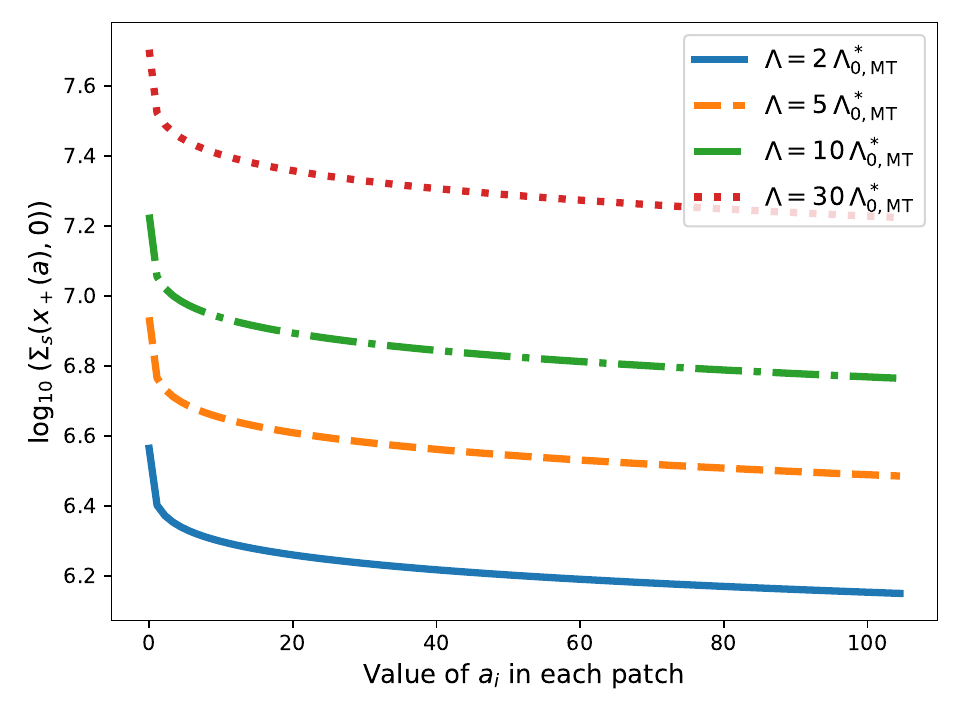}
        \end{subfigure}
        \caption{$3$-patch SIT control with $\mathcal{C}^s=\{2\}$ and MT in Patches \{1,3\}: \textbf{(a)} $\tau(x_+(a),0)$  and  \textbf{(b)} $\Sigma_s(x_+(a),0)$  as functions of the Allee  parameter.  By comparing with \Cref{Fig_a_i_3}, applying MT reduces both $\timetot$ and $\nbtot$ relative to the SIT-only case.}
        \label{fig_estimated_basin}
    \end{minipage}
\end{figure}

 We now investigate the impact of mass trapping (MT) on the release duration. As shown in \Cref{combining SIT with mass trapping}, the use of MT induces a higher critical release rate required for elimination compared to the SIT-only case.
To enable a meaningful comparison of release durations between the SIT-only case and the SIT-MT case, we consider identical release rates, using as a reference the critical release rates obtained with MT, which are larger and denoted by $\Lambda^*_{0,\mathrm{MT}}$ (see \Cref{Table_SIT_MT_patch2}).
By comparing \Cref{fig_estimated_basin}(a) and \Cref{Fig_a_i_3}(a) (see page \pageref{Fig_a_i_3}),  applying MT  reduces $\timetot$,  and consequently $\nbtot$ (compare \Cref{fig_estimated_basin}(b) and \Cref{Fig_a_i_3}(b)).

In all simulations, MT is stopped at the same time as the releases, i.e.,  the entrance time into the basin of attraction of the extinction equilibrium of the uncontrolled model \eqref{log_allee} is evaluated  by setting $\rho_i = \rho_{s,i} = 0$ in \eqref{log_allee}, although trap removal can be costly. In practice, traps may be maintained beyond the  release period. Here, they are removed simultaneously with the end of the releases to consider a conservative scenario; maintaining them would in fact reduce the required release duration.

\subsection{Seven-patch model and impact of the network configuration}

\label{seven patch model}

We consider in this section examples of larger networks, with seven patches ($n=7$). The same configurations were examined in \cite{bliman2024feasibility} with mass trapping (MT) control only.

In the field, with numerous areas, it may be  impossible or too costly to deploy teams to release sterile insects in all patches. Therefore, for each network configuration, three scenarios are examined, with the number of sterile-release patches set consecutively to 1, 2, and 3.  This amounts to solving the minimisation problem \eqref{minimizationproblem_chap3} with an additional constraint on the number of positive entries of the release vector $\Lambda$, i.e.,  for $k = 1,2,3$, solving the following problem:
        \begin{equation}
            \label{minimizationproblem_chap3_card}
            \begin{tabular}{llll}
                \text{minimise}   & $\mathbf{1}_n^\top \Lambda$              \\
                \text{subject to} & $\Lambda ~\overline{\Lambda}$-admissible \\
                                  & $h(\Lambda) \le -10^{-5}$                \\
                                  & $\#\{ i \;:\; \Lambda_i > 0 \} \le k$,
            \end{tabular}
        \end{equation}
        where the function $h$ is defined in \eqref{h}, and $\#\{i: \Lambda_i > 0 \}$ denotes the number of positive components of $\Lambda$. 
        The combinations of $k$ patches that minimise the total  release rate are then identified as the optimal combinations of  release patches.

     		The objective of this section is to assess the impact of the network’s spatial structure on the optimal selection of sterile-release patches. To this aim, we consider two different networks. For each of them, we first consider the unconstrained scenario, in which releases are allowed in all patches $(\mathcal{C}^s = \mathcal{C}^s_I = \{1,\ldots,7\})$. We then examine how prohibiting releases in certain patches, as may occur when landowners do not allow them, can significantly increase the number of sterile insects required. Finally, we study the application of MT in specific areas and its impact on the optimal critical release rate and the release duration.

To see the net effect of the network structure, we choose uniform  local parameters across all patches, with parameters in each patch as in \Cref{Table_parametres_general}, page \pageref{Table_parametres_general}. The dispersal parameters for the wild and sterile insects are assumed identical ($D^s = D$), and the dispersal coefficients between connected patches are  set to $0.02$.

\subsubsection{Configuration 1}

\begin{figure}[H]
    \centering
    \includegraphics[scale=0.8]{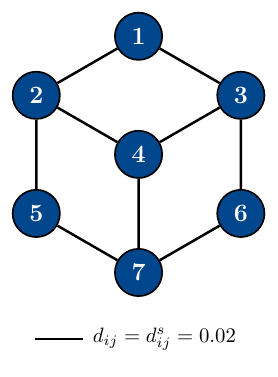}
    \caption{A network configuration with seven patches}
    \label{graph_init}
\end{figure}

Let us begin with the example illustrated in  \Cref{graph_init}, where  lines connecting two  patches  illustrate bidirectional and symmetric dispersal.

 The optimal combinations of $k$ release patches are presented in \Cref{T1_opti1}, page \pageref{T1_opti1}, in the absence of MT, for $k=1, 2$ and $3$.  The corresponding optimal release rates are also reported in the table: first without MT  (SIT-only case), and then with MT in the remaining patches  (SIT-MT case).

Without MT, releasing insects in a single patch requires a total daily release rate of $4259$ individuals, which is $29.4\%$ more than when two patches are available for release ($3291$ individuals). Releasing in two patches results in only a $7.0\%$ higher total daily release compared with three patches ($3075$ individuals).
      Applying SIT in 3 patches rather than 2 thus results in limited savings in terms of sterile insect production, to be compared with the greater organisational complexity.

\begin{figure}[H]
    \centering
    \begin{minipage}{\textwidth}
        \begin{subfigure}[b]{0.4965\textwidth}
            \centering
            \includegraphics[width=\textwidth]{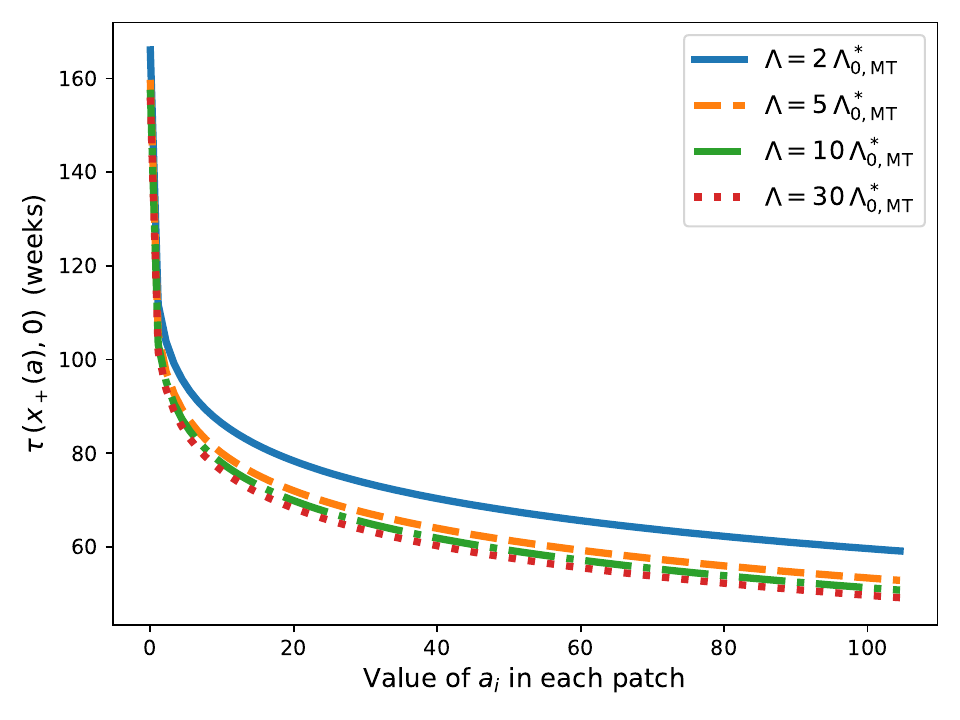}
        \end{subfigure}
        \begin{subfigure}[b]{0.4965\textwidth}
            \centering
            \includegraphics[width=\textwidth]{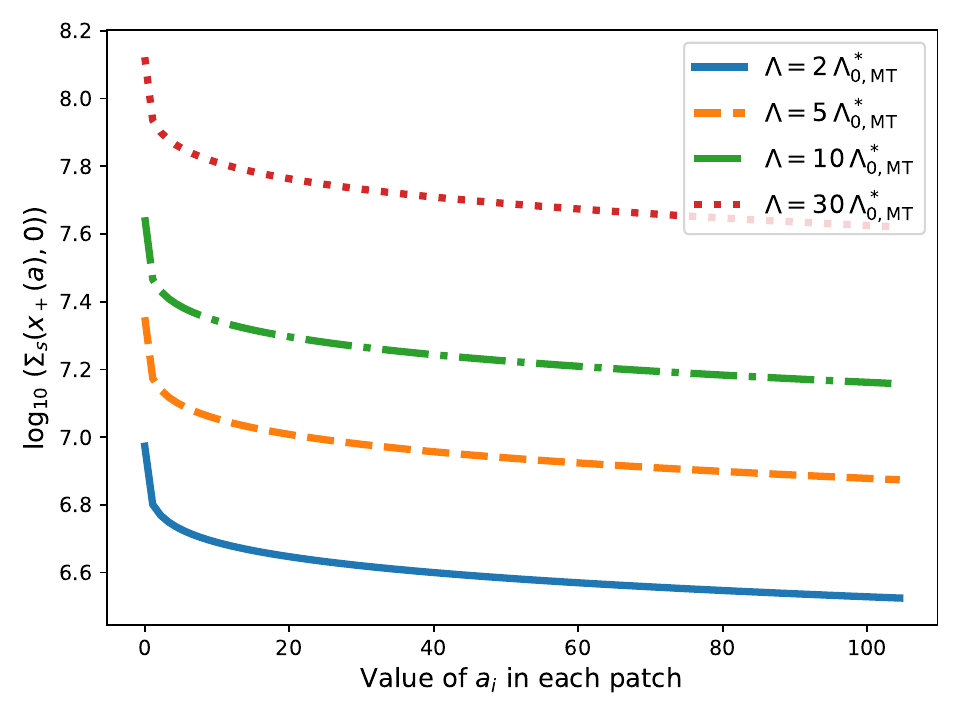}
        \end{subfigure}
    \end{minipage}
    \caption{$7$-patch SIT control  for Configuration 1 with $\mathcal{C}^s=\{2,3,7\}$: \textbf{(a)} $\tau(x_+(a),0)$  and  \textbf{(b)} $\Sigma_s(x_+(a),0)$  as functions of the Allee  parameter.}
    \label{F3}

    \vspace{1em}
    \begin{minipage}{\textwidth}
        \begin{subfigure}[b]{0.4965\textwidth}
            \centering
            \includegraphics[width=\textwidth]{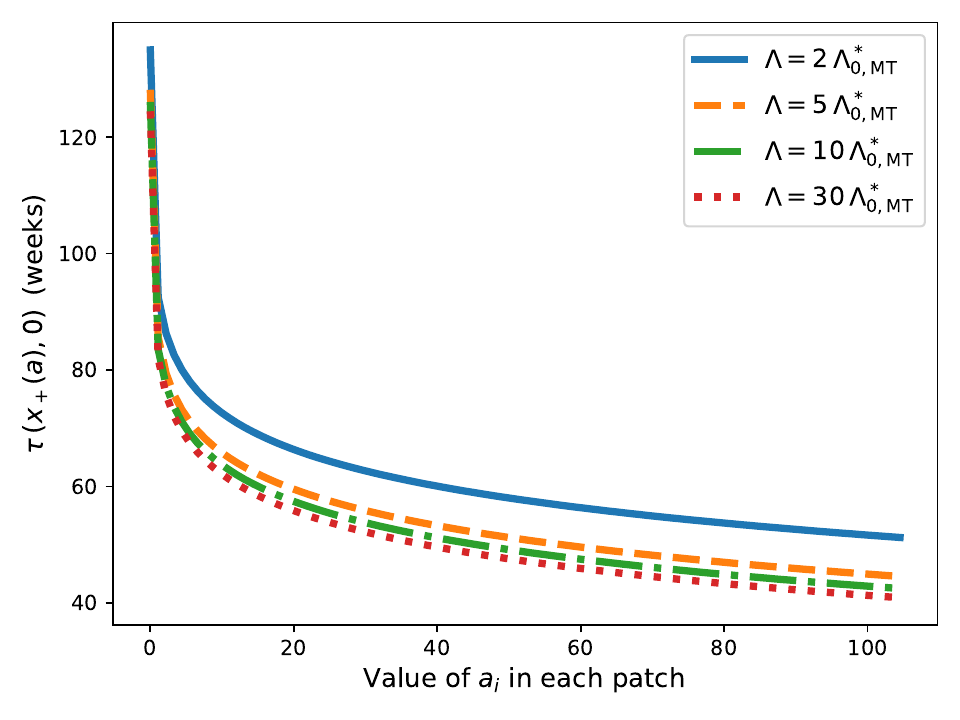}
        \end{subfigure}
        \begin{subfigure}[b]{0.4965\textwidth}
            \centering
            \includegraphics[width=\textwidth]{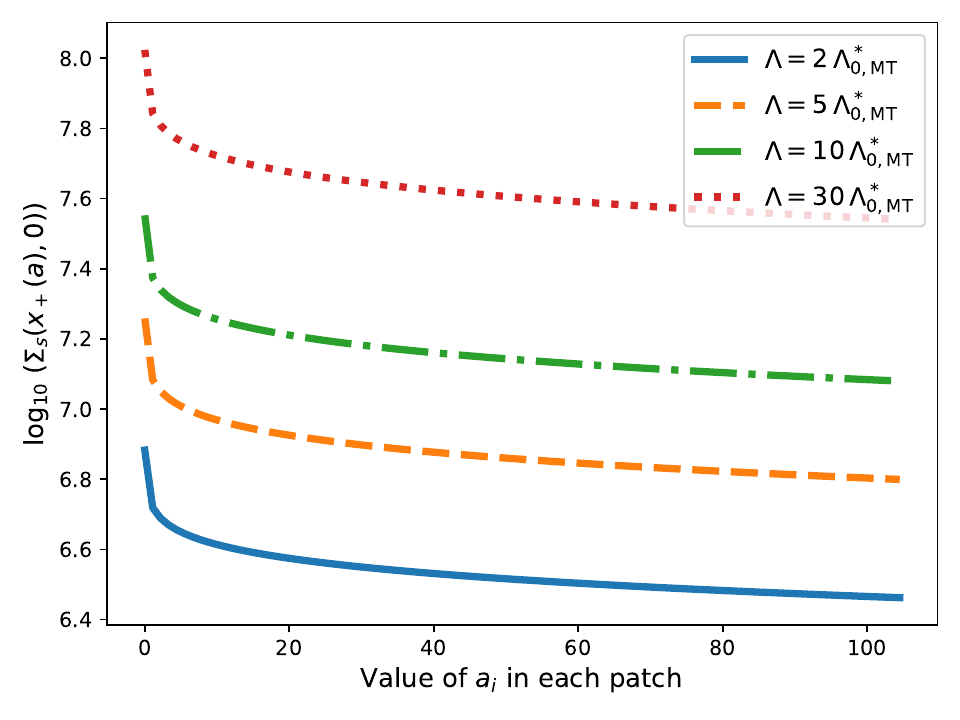}
        \end{subfigure}
    \end{minipage}
    \caption{$7$-patch SIT control  for Configuration 1 with $\mathcal{C}^s=\{2,3,7\}$ and MT in the remaining patches: \textbf{(a)} $\tau(x_+(a),0)$  and  \textbf{(b)} $\Sigma_s(x_+(a),0)$  as functions of the Allee  parameter.  By comparing with \Cref{F3}, applying MT reduces both $\timetot$ and $\nbtot$ relative to the SIT-only case.}
    \label{F4}
\end{figure}

When MT is applied outside the release patches, the daily release rates required for elimination increase according to \Cref{T1_opti1} (see page \pageref{T1_opti1}).  However, for fixed release rates, $\timetot$ decreases (compare  Fig.~\ref{F4}(a) with Fig.~\ref{F3}(a), page \pageref{F3}, for $\mathcal{C}^s=\{2,3,7\}$), leading to a reduction in    $\nbtot$
(compare  Fig.~\ref{F4}(b) with Fig.~\ref{F3}(b)).   Further numerical simulations, not reported here, show that increasing the number of patches where MT is applied leads to higher optimal critical release rates but, for fixed release rates, to shorter release durations and thus lower total quantities $\nbtot$.

\begin{table}[H]
    \centering
    \renewcommand{\arraystretch}{1.2}
    \begin{tabular}{|>{\centering\arraybackslash}m{0.6cm}|
        >{\centering\arraybackslash}m{10.5cm}|
        >{\centering\arraybackslash}m{1.7cm}|
        >{\centering\arraybackslash}m{1.9cm}|}
        \hline
        \multirow{2}{*}{$k$}                                &
        \multirow{2}{*}{Optimal release patch combinations} &
        \multicolumn{2}{c|}{$\sum_{i=1}^7 \Lambda_i^*$}                                          \\ \cline{3-4}
                                                            &                & SIT-only & SIT-MT \\
        \hline
        $1$                                                 & \vspace{0.1cm}
        \begin{minipage}[t]{0.3\linewidth}
            \centering
            \includegraphics[width=0.9\linewidth]{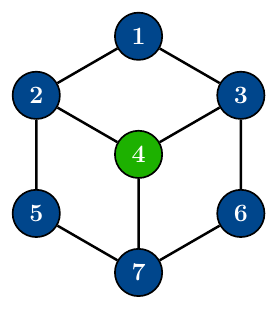}\\
            $\Lambda^*_4 = 4259$\\[2pt]
            $\Lambda^*_{\mathrm{MT},4} = 6687$
        \end{minipage}
                                                            & 4259           & 6687              \\
        \hline
        $2$                                                 & \vspace{0.1cm}
        \begin{minipage}[t]{0.3\linewidth}\centering
            \includegraphics[width=0.9\linewidth]{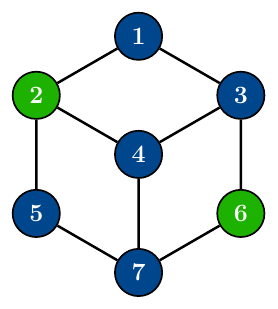}\\
            $(\Lambda^*_2,\Lambda^*_6) = (1882,1409)$\\[2pt]
            $(\Lambda^*_{\mathrm{MT},2},\Lambda^*_{\mathrm{MT},6})=(2566, 1937)$
        \end{minipage}%
        \hspace{0.5em}%
        \begin{minipage}[t]{0.3\linewidth}\centering
            \includegraphics[width=0.9\linewidth]{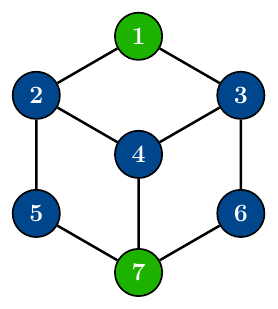}\\
            $(\Lambda^*_1,\Lambda^*_7) = (1409,1882)$\\[2pt]
            $(\Lambda^*_{\mathrm{MT},1},\Lambda^*_{\mathrm{MT},7}) = (1937,2566)$
        \end{minipage}%
        \hspace{0.5em}%
        \begin{minipage}[t]{0.3\linewidth}\centering
            \includegraphics[width=0.9\linewidth]{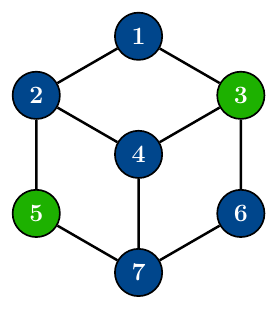}\\
            $(\Lambda^*_3,\Lambda^*_5) = (1882,1409)$\\[2pt]
            $(\Lambda^*_{\mathrm{MT},3},\Lambda^*_{\mathrm{MT},5}) = (2566, 1937)$
        \end{minipage}
                                                            & 3291           & 4503              \\
        \hline
        $3$                                                 & \vspace{0.1cm}
        \begin{minipage}[t]{0.9\linewidth}\centering
            \includegraphics[width=0.3\linewidth]{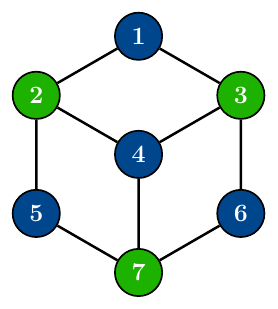}\\
            $(\Lambda^*_2,\Lambda^*_3,\Lambda^*_7)=(1025,1025,1025)$\\[2pt]
            $(\Lambda^*_{\mathrm{MT},2},\Lambda^*_{\mathrm{MT},3},\Lambda^*_{\mathrm{MT},7}) =  (1348,1348,1348)$
        \end{minipage}
                                                            & 3075           & 4044              \\
        \hline
    \end{tabular}
    \caption{Optimal combinations of release patches (in green) depending on the number $k$ of  release sites. The amounts of sterile insects released daily are shown with and without MT in the remaining patches (in blue).}
    \label{T1_opti1}
\end{table}

  Restricting the number of patches on which  SIT is applied is not the only logistical preoccupation: for various reasons, some patches may be unavailable. This constraint  can alter the optimal intervention strategy. That is why  numerical simulations can be useful to  explore  all possible intervention strategies in order to help the field experts and  increase the SIT efficacy.

The (optimal) release strategies when Patches \{4,5,6\} are unavailable for  releases  ($\mathcal{C}^s= \mathcal{C}^s_I = \{1,2,3,7\}$)  are presented in \Cref{T1_opti2}, page \pageref{T1_opti2}. Three categories of patches are considered, with the corresponding colours shown in the table.   Patches where releases cannot be achieved are shown in red, while the set of  optimal patches -- among those where releases are permitted -- is shown in green. In the remaining patches, shown in blue, MT may or may not be applied.

When the releases  are limited to a single patch,  the  best strategy is to release insects in Patch $1$. This  placement necessitates a daily release of $6757$ sterile insects, representing a $58.7\%$ increase compared to the unconstrained scenario, in which the optimal patch is the central  Patch 4 (see \Cref{T1_opti1}).  

If intervention is restricted to at most two patches ($k=2$), among the  three optimal strategies identified in \Cref{T1_opti1} for the unconstrained case, only one remains feasible: Patches $\{1,7\}$. When releases are limited to a maximum of three patches ($k=3$), and since the optimal combination for the unconstrained case does not involve the forbidden Patches \{4,5,6\} (see \Cref{T1_opti1}, page \pageref{T1_opti1}),  the best strategy  remains unchanged.

Naturally, in the absence of MT, the daily release rates for $k = 2$ or $k = 3$ remain identical to those of the unconstrained scenario since the optimal release patches are unchanged. Differences arise only when MT is applied in patches that are neither release patches nor forbidden. In that case, the required daily release decreases relative to the unconstrained setting, since  traps are deployed in fewer patches.

% ...

\begin{table}[H]
    \centering
    \begin{tabular}{|>{\centering\arraybackslash}m{0.6cm}|
        >{\centering\arraybackslash}m{10.5cm}|
        >{\centering\arraybackslash}m{1.7cm}|
        >{\centering\arraybackslash}m{1.9cm}|}
        \hline
        \multirow{2}{*}{$k$}                                &
        \multirow{2}{*}{Optimal release patch combinations} &
        \multicolumn{2}{c|}{$\sum_{i=1}^7 \Lambda_i^*$}                                                                                                                            \\ \cline{3-4}
                                                            &                                                                                                  & SIT-only & SIT-MT \\
        \hline
        $1$                                                 & \vspace{2pt}
        \begin{minipage}[t]{0.9\linewidth}\centering
            \includegraphics[width=0.3\linewidth]{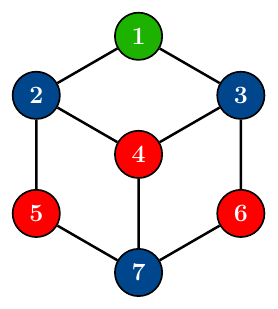}\\
            $\Lambda^*_1 = 6757$\\[2pt]
            $\Lambda^*_{\mathrm{MT},1} = 9705$
        \end{minipage}

                                                            & 6757                                                                                             & 9705              \\
        \hline
        $2$                                                 & \vspace{2pt}
        \begin{minipage}[t]{0.9\linewidth}\centering
            \includegraphics[width=0.3\linewidth]{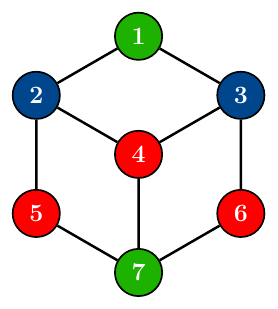}\\
            $(\Lambda^*_1,\Lambda^*_7) = (1409, 1882)$\\[2pt]
            $(\Lambda^*_{\mathrm{MT},1},\Lambda^*_{\mathrm{MT},7}) = (1880, 1926)$
        \end{minipage}

                                                            & 3291                                                                                             & 3806              \\
        \hline
        $3$                                                 & \vspace{2pt}  \begin{minipage}[t]{0.9\linewidth}\centering
                                                                                \includegraphics[width=0.3\linewidth]{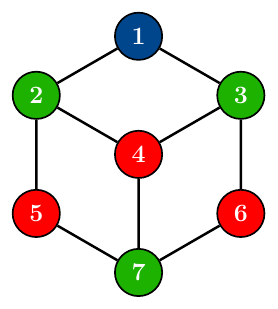}\\
                                                                                $(\Lambda^*_2,\Lambda^*_3,\Lambda^*_7)=(1025,1025,1025)$\\[2pt]
                                                                                $(\Lambda^*_{\mathrm{MT},2},\Lambda^*_{\mathrm{MT},3},\Lambda^*_{\mathrm{MT},7}) =  (1330, 1330, 691)$
                                                                            \end{minipage}
                                                            & 3075                                                                                             & 3351              \\
        \hline
    \end{tabular}
    \caption{Optimal combinations of release patches (in green) depending on the number $k$ of release sites, when some patches (in red) are excluded. The amounts of sterile insects released daily are shown with and without MT in the remaining patches (in blue).}
    \label{T1_opti2}
\end{table}

Also notice that when $k=3$ in \Cref{T1_opti2}, the release rates without MT are similar in Patches $2$, $3$, and $7$, whereas this is no longer the case when MT (in Patch~$1$) is considered. In that case, the release rate in Patch~$7$ is significantly lower than those in Patches~$2$ and $3$, since Patch~$7$ is more distant from Patch~$1$. However, it is important to keep in mind that these values correspond to minimal release rates ensuring elimination in all patches. If sterile insect production is sufficiently large,  field experts may opt for higher release rates in Patch~$7$.

\subsubsection{Configuration 2 -- chain structure}

Let us now consider the scenario where the  graph exhibits a chain structure. This type of structure occurs when the orchards, or plots, are next to each other.  In such cases, it is likely that insects will migrate from one orchard to another sequentially, without skipping any.

The optimal combinations are shown in \Cref{T5_opt}, page \pageref{T5_opt}. In particular, when sterile insects can be released in only one patch, the optimal one is the same as in Configuration 1 (Patch $4$), and requires a daily release of $11\,847$ individuals.  This is $2.78$ times more sterile insects to release than in Configuration $1$: see \Cref{T1_opti1}, page \pageref{T1_opti1}. This shows the importance of the network structure in the optimal size of the releases.

When MT is applied in the remaining patches, $24\,026$ sterile insects per day are required, which is $3.59$ times more than in Configuration~$1$. Thus, the simultaneous use of MT and SIT may render the elimination of the wild population impossible in practice   if the daily production  of sterile insects is insufficient.

\begin{table}[H]
    \centering
    \begin{tabular}{|>{\centering\arraybackslash}m{0.6cm}|
        >{\centering\arraybackslash}m{8cm}|
        >{\centering\arraybackslash}m{3.2cm}|
        >{\centering\arraybackslash}m{2.9cm}|}
        \hline
        \multirow{2}{*}{$k$}                                &
        \multirow{2}{*}{Optimal release patch combinations} &
        \multicolumn{2}{c|}{$\sum_{i=1}^7 \Lambda_i^*$}                                         \\ \cline{3-4}
                                                            &              & SIT-only  & SIT-MT \\
        \hline
        $1$                                                 & \vspace{2pt}
        \begin{minipage}[c]{0.8\linewidth}\centering
            \includegraphics[width=\linewidth]{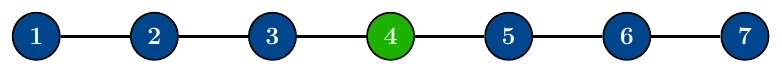}\\
            $\Lambda^*_4 = 11\,847$\\[2pt]
            $\Lambda^*_{\mathrm{MT},4} = 24\,026$
        \end{minipage}

                                                            & $11\,847$    & $24\,026$          \\

        \hline
        $2$                                                 & \vspace{2pt}
        \begin{minipage}[c]{0.8\linewidth}\centering
            \includegraphics[width=\linewidth]{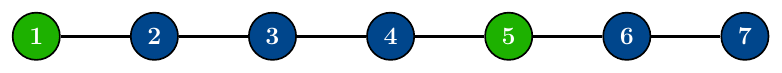}\\
            $(\Lambda^*_1,\Lambda^*_5) = (1110, 4216)$\\[2pt]
            $(\Lambda^*_{\mathrm{MT},1},\Lambda^*_{\mathrm{MT},5}) = (1463, 6963)$\\[2pt]

            \includegraphics[width=\linewidth]{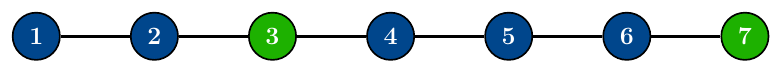}\\
            $(\Lambda^*_3,\Lambda^*_7)= (4216, 1110) $\\[2pt]
            $(\Lambda^*_{\mathrm{MT},3},\Lambda^*_{\mathrm{MT},7}) = (6963,1463)$
        \end{minipage}

                                                            & 5326         & 8426               \\
        \hline
        $3$                                                 & \vspace{2pt}
        \begin{minipage}{0.8\linewidth}\centering

            \includegraphics[width=\linewidth]{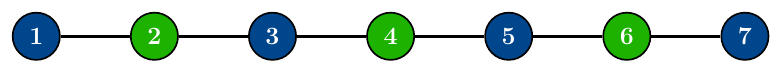}\\
            $(\Lambda^*_2,\Lambda^*_4,\Lambda^*_6) = (1454,380,1454)$\\[2pt]
            $(\Lambda^*_{\mathrm{MT},2},\Lambda^*_{\mathrm{MT},4},\Lambda^*_{\mathrm{MT},6}) = (1971, 386, 1971)$
        \end{minipage}

                                                            & 3288         & 4328               \\
        \hline
    \end{tabular}
    \caption{Optimal combinations  of release patches (in green) depending on the number $k$ of  release sites. The amounts of sterile insects released daily are shown with and without  MT in the remaining patches (in blue).}
    \label{T5_opt}
\end{table}

The release strategies when sterile insects cannot be released in Patches $\{3,4,5\}$ are presented in \Cref{T5_opt_2}, page \pageref{T5_opt_2}. Notably, when sterile insects are released in only one patch, then  Patch $2$ or Patch $6$ is optimal, but  the daily release rate is extremely high compared with the unconstrained scenario ($88\,835$ vs $11\,847$ individuals), which  may be impossible in practice. When MT is applied in the remaining patches, the release rate  rises to the  value of $148\,577$ individuals per day, which is substantially larger.  This example highlights the important additional production cost that may be induced by the unavailability of some patches.

From \Cref{T5_opt_2}, page \pageref{T5_opt_2}, we  also  see that without MT, releasing sterile insects in two or three patches provides the same result, so that it is sufficient  to release  sterile insects in Patches \{2,6\}.

\begin{table}[H]
    \centering
    \begin{tabular}{|>{\centering\arraybackslash}m{0.6cm}|
        >{\centering\arraybackslash}m{8cm}|
        >{\centering\arraybackslash}m{3.2cm}|
        >{\centering\arraybackslash}m{2.9cm}|}
        \hline
        \multirow{2}{*}{$k$}                                &
        \multirow{2}{*}{Optimal release patch combinations} &
        \multicolumn{2}{c|}{$\sum_{i=1}^7 \Lambda_i^*$}                                          \\ \cline{3-4}
                                                            &              & SIT-only   & SIT-MT \\
        \hline
        $1$                                                 & \vspace{2pt}
        \begin{minipage}[c]{0.8\linewidth}\centering
            $~$\\[2pt]
            \includegraphics[width=\linewidth]{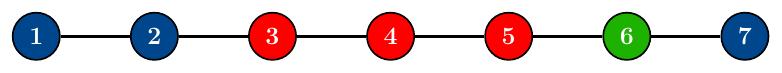}\\[2pt]
            $\Lambda^*_6 = 88\,835$ \\[2pt]
            $\Lambda^*_{\mathrm{MT},6} = 148\,577$\\[4pt] 
            \includegraphics[width=\linewidth]{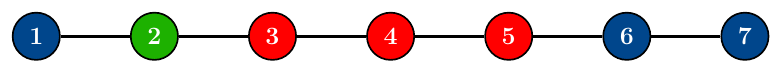}\\[2pt]
            $\Lambda^*_2 = 88\,835$ \\[2pt]
            $\Lambda^*_{\mathrm{MT},2} = 148\,577$
        \end{minipage}

                                                            & $88\,835$    & $148\,577$          \\

        \hline
        $2$                                                 & \vspace{2pt}
        \begin{minipage}[c]{0.8\linewidth}\centering
            $~$\\[2pt]
            \includegraphics[width=\linewidth]{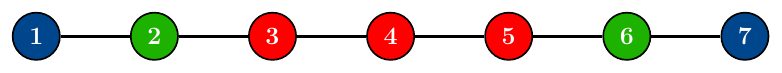}\\[2pt]
            $(\Lambda^*_2,\Lambda^*_6)=(2673,2673)$ \\[2pt]
            $(\Lambda^*_{\mathrm{MT},2},\Lambda^*_{\mathrm{MT},6})=(2787,2787)$
        \end{minipage}

                                                            & 5346         & 5574                \\
        \hline
        $3$                                                 & \vspace{2pt}
        \begin{minipage}[c]{0.8\linewidth}\centering
            $~$\\[2pt]
            \includegraphics[width=\linewidth]{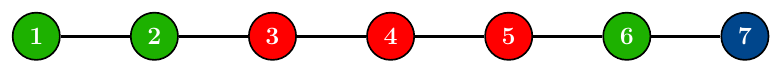}\\[2pt]
            $(\Lambda^*_1,\Lambda^*_2,\Lambda^*_6)=(0,2673, 2673)$\\[2pt]
            $(\Lambda^*_{\mathrm{MT},1},\Lambda^*_{\mathrm{MT},2},\Lambda^*_{\mathrm{MT},6})=(0, 3402, 2030)$

            \includegraphics[width=\linewidth]{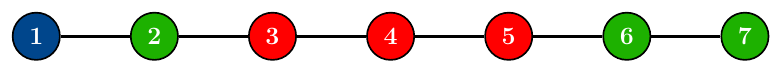}\\[2pt]
            $(\Lambda^*_2,\Lambda^*_6,\Lambda^*_7)=(2673,2673, 0)$\\[2pt]
            $(\Lambda^*_{\mathrm{MT},2},\Lambda^*_{\mathrm{MT},6},\Lambda^*_{\mathrm{MT},7})= (2030, 3402, 0)$
        \end{minipage}

                                                            & 5346         & 5432                \\
        \hline
    \end{tabular}
    \caption{ Optimal combinations of release patches (in green) depending on the number $k$ of release sites, when some patches (in red) are excluded. The amounts of sterile insects  released daily are shown with and without MT in the remaining patches (in blue). }
    \label{T5_opt_2}
\end{table}

\begin{remark}
    In \cite{bliman2024feasibility}, we performed analogous numerical simulations for Configuration 1 and Configuration 2, but we focused on minimising the additional mortality induced by MT and for three controllable patches. Notably, for  these two configurations and when the connectivity matrix is the same for wild and sterile individuals, the optimal combination of three patches for sterile insect release is exactly the same as that for the introduction of additional mortality. Further numerical simulations suggest that  this is still  true when $D^s$ is proportional to $D$, i.e., $D^s = \omega D$  for some $\omega>0$. However, when the networks defined by $D$ and $D^s$ are not correlated, the outcome changes significantly, and the optimal combinations for sterile insect release and for the application of MT have no reason to coincide.
\end{remark}

\subsection{Numerical computations of the release  duration}

 In all the simulations related to release durations that were presented, the exact release durations were determined. In this section, the methodology used for these computations is  presented. An alternative method, based on \Cref{time entrance basin of attraction}, 
    is also provided. Then, we quickly compare the computation times of 
    both methods.

    To estimate the (exact) release duration, we first solve the time-varying system of ODEs \eqref{log_allee_control_time} using the function odeint from the scipy.integrate module of the SciPy package in Python \cite{virtanen2020scipy}. This function relies on the LSODA solver from the ODEPACK library \cite{hindmarsh1983odepack}, which automatically switches between methods for stiff and non-stiff systems \cite{petzold1983automatic}.
    A bisection  procedure is then  applied along the simulated trajectory.  We consider an  interval $[t_{\min}, t_{\max}]$ and iteratively refine it. At each step, we set $t' = (t_{\min} + t_{\max})/2$ and evaluate the state of the  system \eqref{log_allee_control_time} at time $t'$. This state is then used as an initial condition for the uncontrolled system \eqref{log_allee}, and we test whether the corresponding trajectory converges to the origin.   If convergence occurs,  we update $t_{\max} \leftarrow t'$, otherwise, we update $t_{\min} \leftarrow t'$.

    With this simple estimation method, assessing the exact value of the release  duration becomes computationally demanding as the number of patches increases. The formula given in \Cref{time entrance basin of attraction} to approximate the release duration therefore offers  substantial computational savings, at the price of reduced accuracy (but guaranteed overestimation).
    To apply this theorem, and as  with estimating the exact release duration, 
    we solve numerically the time-varying system of  ODEs \eqref{log_allee_control_time}, and then apply a bisection procedure along the simulated trajectory.
    While, at each step of the bisection procedure, the previous approach directly checks convergence of the uncontrolled trajectory, this second method replaces this test with the condition  $ \Psi\left(c\t x^k \right) < -s(J_a)$, where $x^k$ denotes the state at the current iteration, and  $c \gg 0_n$, and $\Psi$ are defined in \Cref{time entrance basin of attraction}.

\begin{figure}[H]
    \centering

    \begin{subfigure}[b]{0.49\textwidth}
        \centering
        \includegraphics[width=\textwidth]{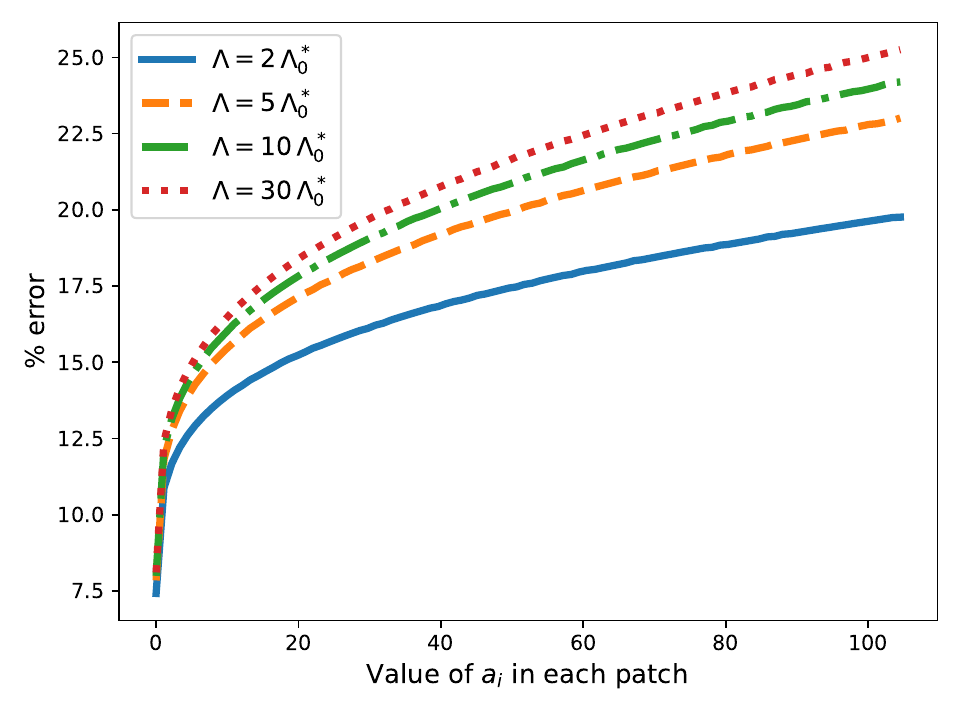}
        \caption{$3$-patch model with $\mathcal{C}^s = \{3\}$}
    \end{subfigure}
    \hfill
    \begin{subfigure}[b]{0.49\textwidth}
        \centering
        \includegraphics[width=\textwidth]{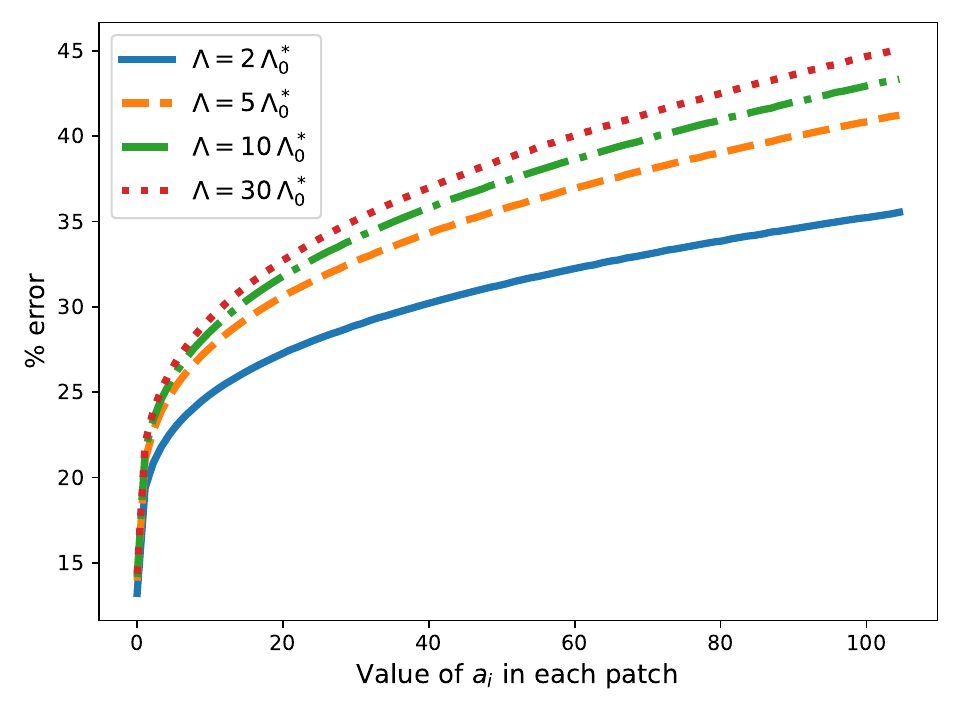}
        \caption{$7$-patch model with $\mathcal{C}^s = \{2,3,7\}$}
    \end{subfigure}

    \vspace{0.3cm}

    \begin{subfigure}[b]{0.49\textwidth}
        \centering
        \includegraphics[width=\textwidth]{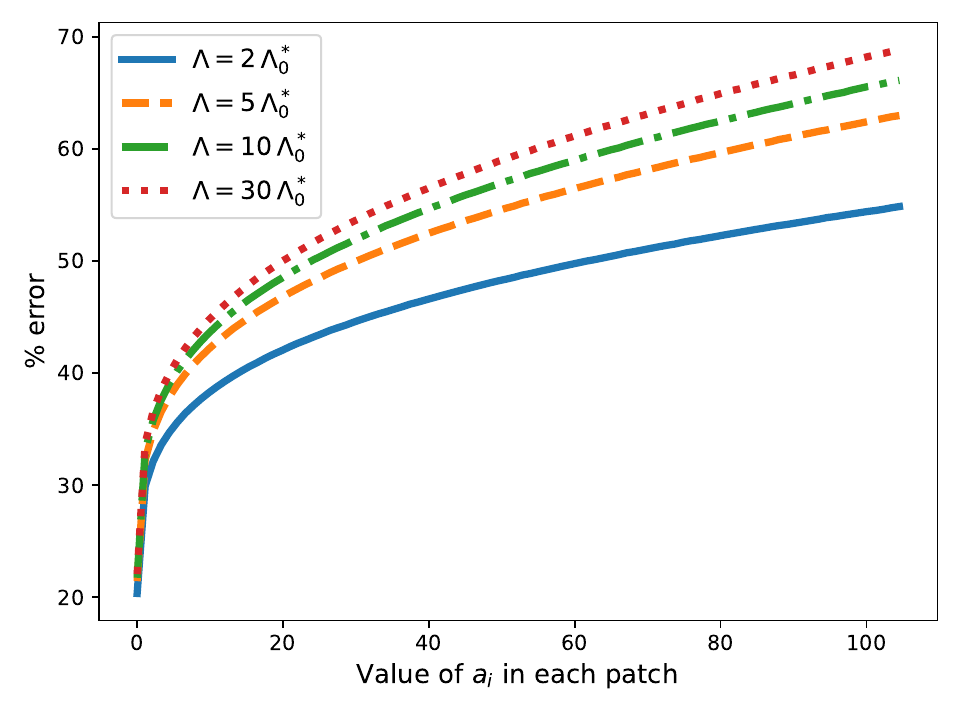}
        \caption{$20$-patch model with $\mathcal{C}^s=\{6, 8, 10, 16, 18\}$}
    \end{subfigure}

    \caption{Overestimates (in percent) of $\timetot$ with \Cref{time entrance basin of attraction},  as a function of the  Allee parameter.  The overestimation increases with  the number of patches in the network, the release rates and the Allee parameter.}
    \label{Comp_error}
\end{figure}

    A summary of the computation times for the exact and estimated release durations
    is provided in \Cref{tab:computation_times}, for the 3-patch and 7-patch (see \Cref{graph_init}) cases, as well as for a network of $20$ patches with a grid structure as in \Cref{F28}.
    It is straightforward to check that the computation of the approximate  release duration is, on average, $7$ to $8$ times faster than the exact estimate. However, as shown in \Cref{Comp_error}, this gain in computational time is offset by a loss of precision in the release duration, which is overestimated in a proportion of up  to 70\% and increases with both the Allee parameters $a_i$ and the release rates. Moreover, a comparison of Figs.~\ref{Comp_error}(a)--(c) also shows that this overestimation increases as the number of patches in the network increases.

\begin{table}[H]
    \small
    \centering
    \begin{tabular}{|c|c|c|}
        \hline
        \multirow{2}{*}{Number of patches}
           & \multicolumn{2}{c|}{Computational time (seconds) to estimate $\timetot$}         \\
        \cline{2-3}
           & Exact estimation
           & Estimation with \Cref{time entrance basin of attraction}                         \\
        \hline
        3  & 17.48                                                                    & 2.46  \\
        7  & 34.64                                                                    & 4.49  \\
        20 & 111.35                                                                   & 14.72 \\
        \hline
    \end{tabular}
    \caption{Numerical computation time for the release duration, for $96$ values of $a_i \in (0,100]$ and all tested values of $\Lambda = p \times \Lambda_0^*$ with $p =  2, 5, 10,30$.   Computations have been made on a laptop with an Apple M3 chip and $16$ Gb of RAM.}
    \label{tab:computation_times}
\end{table}

\begin{figure}[H]
    \centering

    \begin{tikzpicture}[scale=1.2,
            node/.style={
                    circle, draw, minimum size=8mm,
                    font=\small\bfseries, fill=bleu_fonce, text=white
                },
        ]

        \def\rows{4}
        \def\cols{5}

        \foreach \i in {0,1,2,3}{
                \foreach \j in {0,1,2,3,4}{
                        \pgfmathtruncatemacro{\idx}{\i*\cols + \j + 1}
                        \node[node] (n\i\j) at (\j,-\i) {\idx};
                    }
            }

        \foreach \i in {0,1,2,3}{
                \foreach \j in {0,1,2,3}{
                        \pgfmathtruncatemacro{\jp}{\j+1}
                        \draw (n\i\j) -- (n\i\jp);
                    }
            }

        \foreach \i in {0,1,2}{
                \foreach \j in {0,1,2,3,4}{
                        \pgfmathtruncatemacro{\ip}{\i+1}
                        \draw (n\i\j) -- (n\ip\j);
                    }
            }

    \end{tikzpicture}
    \caption{20-patch network.}
    \label{F28}
\end{figure}

\section{Conclusion}

  In order to design spatial SIT strategies, an explicit metapopulation model with an arbitrary number of patches has been proposed and studied, allowing for distinct dispersal networks for wild and sterile populations. Considering constant sterile insect release rates,
  a sufficient condition  has been established to ensure elimination of the pest/vector population under various constraints, including 
  the presence of patches where releases are not feasible. The derived elimination condition  provides a unifying criterion that combines local patch-level conditions and dispersal dynamics into a single framework. Building on this result, a convex constrained optimisation problem has been formulated and numerically solved to minimise a weighted sum of the release rates required for elimination. This framework enables    the identification of efficient  spatial release strategies  to achieve population elimination. Among the admissible strategies, field experts may  choose those that are most suitable for implementation in realistic field settings.

Through numerical simulations, we studied different scenarios, including different network structures.  The results obtained  highlight the importance of taking  dispersal rates and network configuration into account when designing  SIT strategies. In addition, from an integrated control perspective, we  assessed the impact
of combining mass trapping (MT) and SIT on the costs of the latter.
MT also increases the sterile insect mortality, and the numerical simulations show a subsequent increase in the critical release rates required for elimination.
     However, for fixed release rates, MT becomes beneficial in the presence of a natural strong Allee effect  in every patch, as it reduces the  release duration and thus the total number of sterile insects released over the entire SIT  programme. Notice that the negative effect
of MT  on the critical release rates may be limited if the sterile insects are treated before the releases in
order to reduce their susceptibility to the traps.
     For instance, sterile males reared on a methyl eugenol (ME)-based diet exhibit reduced attraction to ME-based traps \cite{wee2025evaluation}.

Overall, our results show that effective strategies must consider not only critical release rates  required to eliminate the population, but also release duration, total releases until elimination, and logistical constraints. For instance,  although releasing in all patches is often optimal, targeting a few key patches may reduce operational costs, requiring a balance between production capacity, time to elimination, and field limitations.

Finally, as our model and results are generic, this approach could be applied to other SIT-related problems involving different pests or disease vectors, such as mosquitoes, provided that appropriate parametrisation is used.


     Some perspectives for future work include the development of a sex- and stage-structured model for $n$ patches, in order to incorporate control methods that target specific sexes and/or developmental stages. For instance, the Male Annihilation Technique is an attract-and-kill strategy that targets adult males, using attractants such as ME \cite{singh2022mass}.
        As another example,  entomopathogenic fungi can be used as biological control agents. These fungi mainly affect pupal mortality, but also reduce adult survival \cite{li2024toxicity, dumont2025improvement}.

\section*{Acknowledgements}
This work was partly supported by the AttracTIS project, funded by ECOPHYTO 2021-2022: ``Construire avec les outre-mer une agro\'ecologie ax\'ee sur la r\'eduction de l’utilisation, des risques et des impacts des produits phytopharmaceutiques''. YD acknowledges the support of the BactroTIS project, grant number ANR-25-ECOM-0031-01.
YD is also (partially) supported by the DST/NRF SARChI Chair in Mathematical Models and Methods in Biosciences and Bioengineering at the University of Pretoria (Grant 82770). YD, MdT and PAB acknowledge the support of the Conseil R\'egional de la R\'eunion, the European Regional Development Fund (ERDF), and the Centre de Coop\'eration Internationale en Recherche Agronomique pour le D\'eveloppement (CIRAD).

\addcontentsline{toc}{section}{References}

% ...

\bibliographystyle{ieeetr}
\bibliography{bibliography.bib}

\section*{Appendix - Proofs}
\label{Appendix}
\addcontentsline{toc}{section}{Appendix -- Proofs}
\appendix

\renewcommand{\theequation}{\thesection.\arabic{equation}}

\section{Theoretical tools}

In this section, we prove results  that are useful for the analysis of model \eqref{log_allee}.

Consider the following system of ODEs:
\begin{equation}
    \label{system_appendix}
    \dot{x} = f(x),
\end{equation}
where $f : \mathcal{D} \rightarrow \mathbb{R}^n$ is differentiable and  $\mathcal{D}$ is an  open  $p$-convex set of  $\mathbb{R}^n$. Assume that all solutions are global.

\begin{theorem}
    \label{positive_equilibrium_irred}
    Assume  system \eqref{system_appendix}  is cooperative and irreducible on $\mathcal{D}$, $\mathcal{D}$ is positively invariant with respect to \eqref{system_appendix}  and let $x^*$ be an equilibrium.  Then, for any $x' \in \mathcal{D}$ such that $f(x') \in \mathbb{R}^n_+$, the following holds:
    $$x^* > x' \Rightarrow x^* \gg x'.$$
\end{theorem}

\begin{proof}
    Let  $x_0 > x'$. Then, as a consequence
    of Theorem 1.1 in \cite[Chapter 4]{smith1995monotone}, any solutions $x(t,x_0)$ and  $x(t,x')$ of system \eqref{system_appendix} with initial condition $x_0$ and $x'$, respectively, satisfy
    \begin{equation}
        \label{kamke_irred}
        x(t,x_0) \gg x(t, x'), \quad t>0.
    \end{equation}
    Since $f(x') \geq 0_n$, then,  by applying  \cite[Proposition 2.1, Chapter 3]{smith1995monotone}, the solution $x(t, x')$ is non-decreasing, which implies $x(t, x') \geq x'$ for all $t \in \mathbb{R}_+$. Therefore, by \eqref{kamke_irred}, it follows $x(t,x_0) \gg x'$ for  all $t>0$.  In particular, any equilibrium $x^*$ such that $x^*> x'$  must satisfy $x^* \gg x'$.
\end{proof}

Now consider the system
\begin{equation}
    \label{system2}
    \dot{y} = g(y),
\end{equation}
where $g : \mathcal{D} \rightarrow \mathbb{R}^n$ is differentiable.

\begin{theorem}
    \label{decrease_equilibrium_point_cooperative_init}
    Assume  system \eqref{system_appendix} or system \eqref{system2} is cooperative and irreducible on $\mathcal{D}$, $\mathcal{D}$ is positively invariant with respect to \eqref{system_appendix} and \eqref{system2} and  $f(x) \le g(x)$ for  every $x \in \mathcal{D}$. Then, for  every $x_0 \in \mathcal{D}$, the following holds:
    $$f(x_0) < g(x_0) \Rightarrow x(t,x_0) \ll y(t,x_0), \quad  t>0, $$
    where  $x(t,x_0)$ and $y(t,x_0)$ denote the respective solutions of \eqref{system_appendix} and \eqref{system2} with initial condition $x_0$.
\end{theorem}

\begin{proof}
    Let $x_0 \in \mathcal{D}$.   Since $f(x_0) < g(x_0)$, we have by continuity  that $x(t, x_0) < y(t,x_0)$ on  an interval $(0, \delta)$ where $\delta>0$. Fix $t> 0$.
    We can write $t = (t- t') + t'$, where  $t' \in (0,\min\{t,\delta\})$. Since  $x(t' , x_0) < y( t',x_0)$,  we obtain by \cite[Theorem 1.1, Chapter 4]{smith1995monotone} that:
    \begin{equation}
        \label{e1}
        x\big(t -  t', x(t',x_0)\big)  \ll x\big(t - t',y( t',x_0)\big).
    \end{equation}
    Since system \eqref{system_appendix} or \eqref{system2} is cooperative on $\mathcal{D}$ with $f \le g$, we also  have  by \cite[Theorem 10]{coppel1965stability} that
    \begin{equation}
        \label{e2}
        x\big(t -  t',y(t',x_0)\big)\le  y\big(t - t',y( t',x_0)\big).
    \end{equation}
    Combining \eqref{e1} and \eqref{e2} implies
    $$x\big(t -  t', x( t',x_0)\big)  \ll y\big(t -  t',y( t',x_0)\big).$$
    On the other hand, since the systems \eqref{system_appendix} and \eqref{system2} are autonomous, we have  $x\big(t -  t', x( t',x_0)\big)  = x(t -  t' +  t',x_0)  = x(t,x_0)$ and similarly, $y\big(t -  t', y( t',x_0)\big)   = y(t,x_0)$. This  yields
    $$x(t,x_0) \ll y(t,x_0).$$
\end{proof}

\section{One-patch model (\Cref{allee_effet_simple})}

\subsection{Proof of \Cref{allee_effet_simple}}

\label{proof_of_allee_effet_simple}

\begin{proof}
    Assume $a > 0$. First, the derivative of the right-hand side of \eqref{simple_allee} when $x=0$ is equal to $-\mu_1 - \rho < 0$, which implies that $0$ is LAS. Moreover, the positive equilibria  of \eqref{simple_allee} are exactly the values $x>0$   such that
    $$b\frac{x}{x + a} - \mu_1 - \rho- \mu_2 x = 0,$$
    that is, the positive roots of the polynomial equation
    \begin{equation}
        \label{eq2_chap3}
        \frac{x^2}{Q} + \left( \frac{a}{Q}+ 1 - \mathcal{N} \right)x + a = 0.
    \end{equation}
    Since $Qa > 0$, the roots of this equation, when real, have the same  sign. Moreover, the sum of the roots is equal to $Q(\mathcal{N}-1)-a$, which implies that the roots, when real, are both negative when $\mathcal{N} \le 1$. In this case, system \eqref{simple_allee} admits no positive equilibrium. Since the polynomial in \eqref{eq2_chap3} is  positive on $\mathbb{R}_+$, the origin of \eqref{simple_allee} is then GAS on $\mathbb{R}_+$.

    Now assume $\mathcal{N} > 1$. The discriminant of   \Cref{eq2_chap3} is
    \begin{align*}
        \Delta & :=\left( \frac{a}{Q}+ 1 - \mathcal{N} \right)^2 - \frac{4a}{Q}          \\
               & = \frac{1}{Q^2} \left(a - a^\text{crit} \right)  \left(a - a_+ \right),
    \end{align*}
    where
    $$a_+ :=  Q\left( \sqrt{\mathcal{N}} + 1\right)^2 > a^\text{crit}.$$

    \noindent $\bullet$ When $0 < a < a^\text{crit}$ or  $a > a_+$, we have $\Delta > 0$, so the polynomial equation \eqref{eq2_chap3} admits two real roots, given by
    \begin{align*}
        x_\pm & = \frac{Q}{2} \left(\mathcal{N}-1 - \frac{a}{Q} \pm \sqrt{\left(\mathcal{N} - 1- \frac{a}{Q} \right)^2 - \frac{4a}{Q}} \right)                \\
              & =  \frac{Q}{2} \left(\mathcal{N} - 1 - \frac{a}{Q}\right)\left(1 \pm \sqrt{1-\frac{4a}{Q\left(\mathcal{N}-1 - \frac{a}{Q}\right)^2}} \right).
    \end{align*}
    If $a > a_+$, then $\mathcal{N} - 1 - \frac{a}{Q} < 0$, which implies that the two roots are negative, and there is no positive equilibrium. The origin is then GAS on  $\mathbb{R}_+$. If, on the other hand, $0 < a < a^{\text{crit}}$,   then $\mathcal{N} - 1 - \frac{a}{Q} > 0$, and the two roots are positive. Since the polynomial  in \eqref{eq2_chap3} is negative between the roots and positive  outside, it follows that
    $x_+$  is LAS,  $x_-$ is unstable, and the basins of attraction of $0$ and $x_+$ are
    $[0,x_-)$ and $(x_-,+\infty)$, respectively.

    \vspace{0.5cm}

    \noindent $\bullet$ When $a = a^{\text{crit}}$, then $\Delta = 0$. Therefore, the equation \eqref{eq2_chap3} admits a unique  root, given by
    $$x^* :=  Q \left(\sqrt{\mathcal{N}}-1 \right).$$
    Since $\mathcal{N} > 1$, this root is positive,  and thus the polynomial  in \eqref{eq2_chap3} is non-negative. It follows that the basins of attraction of $0$ and $x^*$ in $\mathbb{R}_+$ are $[0, x^*)$ and $[x^*,+\infty)$, respectively.

    \vspace{0.5cm}

    \noindent $\bullet$ When $a = a_+$, then $\Delta = 0$, and once again, \eqref{eq2_chap3} admits a unique  root, equal to $-Q \left( \sqrt{\mathcal{N}}+1\right)$.  Since this root is negative, there is no positive equilibrium.

    \vspace{0.5cm}

    \noindent $\bullet$ When $a \in (a^\text{crit}, a_+)$,  then $\Delta < 0$,  which implies that \eqref{eq2_chap3} admits no real root, and the system \eqref{simple_allee} admits no positive equilibrium.
\end{proof}

\subsection{Proof of \Cref{max}}

\label{proof_max}

If \eqref{condition_elim} holds,  then by \eqref{max_equality}, $g(x) < 0$ for all $x \in \mathbb{R}_+$.   Therefore, for any $x_0 \in \mathbb{R}_+$, the solution $x(t,x_0)$ of \eqref{simple_allee} with initial condition $x_0$ is strictly decreasing.
This implies that the origin of \eqref{simple_allee} is GAS on $\mathbb{R}_+$.

Now assume $\mathcal{N}> 1$. To prove that the condition \eqref{condition_elim} is also necessary for the origin of \eqref{simple_allee} to be GAS on $\mathbb{R}_+$, we  simply show that
\begin{equation}
    \label{sens1}
    a > a^\text{crit} \Rightarrow - \mu_1 - \rho +\left(\left( \sqrt{b} - \sqrt{a \mu_2} \right)^+\right)^2
    < 0,
\end{equation}
so that we can apply $\textit{2.(a)}$ in \Cref{allee_effet_simple}.

Assume $a> a^\text{crit}$.

\begin{itemize}
    \item
          If  $a \geq \mathcal{N} Q$, then $- \mu_1-\rho+ \left(\left( \sqrt{b} - \sqrt{a \mu_2} \right)^+\right)^2 = - \mu_1 - \rho
              < 0$, and thus  \eqref{sens1} is satisfied.

    \item   If  $0 \le a < \mathcal{N}Q$, then $- \mu_1-\rho+ \left(\left( \sqrt{b} - \sqrt{a \mu_2} \right)^+\right)^2
              < 0$ is equivalent to
          \begin{equation}
              \label{poly0}
              -\mu_1 - \rho + b + a\mu_2 - 2\sqrt{a \mu_2 b } < 0.
          \end{equation}
          Let $X := \sqrt{a \mu_2}$. The  inequality \eqref{poly0} is then equivalent to
          \begin{equation*}
              X^2-2\sqrt{b}X + b-\mu_1- \rho < 0.
          \end{equation*}
          The discriminant of this polynomial is equal to $4(\mu_1+\rho)$, so it admits two roots
          $$X_\pm := \sqrt{b} \pm ±
              \sqrt{\mu_1+\rho},$$
            which are positive since $\mathcal{N} >1$.
          We deduce that \eqref{poly0} is equivalent to
                $$\sqrt{a} > \frac{\sqrt{b} - \sqrt{\mu_1+ \rho}}{\sqrt{\mu_2}} = \sqrt{Q}(\sqrt{\mathcal{N}}-1) \qquad \text{and} \qquad \sqrt{a} < \frac{\sqrt{b} + \sqrt{\mu_1+ \rho}}{\sqrt{\mu_2}}= \sqrt{Q}(\sqrt{\mathcal{N}}+1),$$
          and thus
          $$a > a^\text{crit} \qquad \text{and} \qquad a < Q(\sqrt{\mathcal{N}}+1)^2.$$
          The second inequality is satisfied  since $a < \mathcal{N}Q$, and thus we deduce that \eqref{sens1} holds.
\end{itemize}

\section{Basin of attraction (\Cref{estimation basin of attraction})}

\label{proof_estimation_basin}

Since $D$ is irreducible, the matrix $J_a$ is also irreducible. Then, by the extension of the Perron-Frobenius theorem, there exists a vector $c \gg 0_n$, such that
$$c\t J_a = -\lambda c\t, \quad \text{with } \lambda := -s(J_a)> 0.$$
Up to multiplying $c$ by $\frac{1}{\min_i c_i}$, we assume that $\min_i c_i = 1$, as done in the statement.

For any $x\in \mathbb{R}^n_+$,  now consider the  Lyapunov function candidate
$$V(x) := c\t x.$$
The   system \eqref{log_allee} rewrites
$$\dot{x} = J_ax + \diag \biggl(\Big(b_i\left[ p(x_i,a_i) - p_{a_i}(0)\right]- \mu_{2,i} x_i\Big)x_i\biggr).$$
Therefore, the  derivative of $V$ along a trajectory of \eqref{log_allee}  satisfies
\begin{align*}
    \dot{V}(x) & = - \lambda c\t x + \sum_{i = 1}^{n} c_i x_i \big(b_i\left[ p(x_i,a_i) - p_{a_i}(0)\right]- \mu_{2,i} x_i\big)                          \\
               & \le- \lambda c\t x + \sum_{i = 1}^{n}  c_i x_i \max_{1 \le k \le n}  \left\{ b_k\left[p_{a_k}(x_k) - p_{a_k}(0)\right]- \mu_{2,k} x_k\right\} \\
               & = c\t x \left(-\lambda+ \max_{1 \le k \le n}  \left\{ b_k\left[p_{a_k}(x_k) - p_{a_k}(0)\right]- \mu_{2,k} x_k\right\} \right).
\end{align*}
   Since for all $k = 1,\ldots,n$, $x_k \le \frac{c\t x}{\min_{i}c_i} =  c\t x$  because $\min_{i}c_i = 1$,  we have
\begin{align*}
    \dot{V}(x) & \le  c\t x \left(-\lambda+   \max_{1 \le k \le n} ~ \max_{0 \le z \le  c\t x}  \{ b_k\left[p_{a_k}(z) - p_{a_k}(0)\right] - \mu_{2,k} z \} \right).
\end{align*}

$\bullet$ If $a_k = 0$, then $b_k\left[p_{a_k}(z) - p_{a_k}(0)\right]- \mu_{2,k} z = -\mu_{2,k}z \le 0$. Thus, the maximum of the function $z \mapsto b_k\left[p_{a_k}(z) - p_{a_k}(0)\right] - \mu_{2,k} z$ is attained at $z = 0$ and is equal to $0$.

\vspace{0.2cm}

$\bullet$ If  $0 < a_k < \frac{b_k}{\mu_{2,k}}$, then $ b_k\left[p_{a_k}(z) - p_{a_k}(0)\right] - \mu_{2,k} z = \frac{z}{z+a_k} - \mu_{2,k} z$. This function
is increasing on
$\left[0, -a_k + \sqrt{\frac{a_k b_k}{\mu_{2,k}}} \right]$
and decreasing on
$
    \left[-a_k + \sqrt{\frac{a_k b_k}{\mu_{2,k}}}, +\infty \right),
        $
            with its maximal value, attained at
        $
        z^*_k = -a_k + \sqrt{\frac{a_k b_k}{\mu_{2,k}}}
        $,
            equal to
        $
        \left(\sqrt{b_k} - \sqrt{a_k \mu_{2,k}}\right)^2
        $.

            \vspace{0.2cm}

        $\bullet$ If  $a_k \geq \frac{b_k}{\mu_{2,k}}$, then the maximum is attained at $z = 0$ and is equal to 0.

            Therefore,
            \begin{align*}
                \max_{1 \le k \le n} ~ \max_{0 \le z \le  c\t x}  \{ b_k\left[p_{a_k}(z) - p_{a_k}(0)\right] - \mu_{2,k} z \} & =     \max_{\substack{1 \le k \le n \\ a_k > 0}} ~ \max_{0 \le z \le  c\t x} ~b_k\frac{z}{z+a_k} - \mu_{2,k} z\\ & =     \max_{\substack{1 \le k \le n \\ a_k > 0}} ~\phi_k( c^\top x).
            \end{align*}
            Thus, this implies that
            $$ \dot{V}(x)   \le   c\t x \left(-\lambda+ \Psi( c\t x) \right).$$
            When  $\Psi( c\t x) <  \lambda = -s(J_a)$, we have $\dot{V}(x) \le 0$ for all $x \in \mathbb{R}^n_+$. Then, since $c \gg 0_n$, $\dot{V}(x) = 0$ if and only if $x = 0_n$.  Therefore, $V(x)$ is a strict Lyapunov function on the set $\{ x \in \mathbb{R}^n_+:  \Psi( c\t x) < - s(J_a)\}$.  This shows that the origin of \eqref{log_allee} is GAS on this set, which proves that  it is included in the basin of attraction $\mathcal{B}_a$ of the origin of \eqref{log_allee}.  Moreover, since $\Psi(0) = 0 < -s(J_a)$, this set contains the origin.
            This completes  the proof of \Cref{estimation basin of attraction}.

            \section{Maximal equilibrium (\Cref{0only})}

            \label{prooof_0only}
            To prove that system \eqref{log_allee} admits a maximal equilibrium $x_+(a)$,   we first show that all the solutions of \eqref{log_allee} are uniformly bounded.  We compute the sum of the $\dot{x}_i$:
            $$ \sum\limits_{\substack{i=1}}^{n}{\dot{x}_i} = \sum\limits_{\substack{i=1}}^{n}{ x_i \big(p(x_i,a_i) b_i  - \mu_{1,i}  - \rho_i- \mu_{2,i} x_i\big)} \le \sum\limits_{\substack{i=1}}^{n}{ x_i \left( b_i  - \mu_{1,i}  - \rho_i- \mu_{2,i} x_i\right)} .$$
            Therefore,  $\sum\limits_{\substack{i=1}}^{n}{\dot{x}_i} < 0$ for sufficiently large $x_i$, $ i = 1,\ldots,n$, which implies that the solutions  of \eqref{log_allee} are uniformly bounded.

            We can  now prove the existence of a maximal equilibrium.  Since  the matrix $D - \diag(\mu_{1,i}+\rho_i)$ is  Metzler and Hurwitz,  the extension of the Perron-Frobenius theorem implies   the existence of a vector $v > 0_n$ such that
            \begin{equation*}
                \big(D - \diag(\mu_{1,i}+\rho_i)\big) v < 0_n.
            \end{equation*}
            Since the solutions of \eqref{log_allee} are uniformly bounded, there exists $\lambda>0$ such that any equilibrium $x^*_a$ satisfies
            \begin{equation}
                \label{x_a < x_0}
                x^*_a \le \lambda v.
            \end{equation}
            Moreover, as the right-hand side $f_a$ of system \eqref{log_allee} satisfies  $f_a \le f_0$ for any $a \in \mathbb{R}^n_+$, one can take $\lambda> 0$ large enough  such that
            \begin{equation}
                \label{ineq}
                f_a(\lambda v) \le f_0(\lambda v) = \lambda \diag(b_i -  \mu_{2,i}\lambda v_i) v+ \lambda \big(D - \diag(\mu_{1,i}+\rho_i)\big)v < 0_n.
            \end{equation}
            Since the solution  of \eqref{log_allee} with initial condition $\lambda v$ is decreasing, it converges towards a non-negative equilibrium   $x_+(a)$. Let $x^*_a$ be another non-negative equilibrium.  By \eqref{x_a < x_0}, it satisfies $ x^*_a \le \lambda v$. By cooperativeness, this implies that, at any time, the solution initiated in  $ \lambda v$ remains above $x^*_a$, and hence $x^*_a \le x_+(a)$. Moreover, by \eqref{ineq} and  by \cite[Theorem 6]{anguelov2012mathematical}, $x_+(a)$ is GAS on the set $\{x \geq x_+(a) \}$.

            Last, when $D$ is irreducible, \eqref{positive_eq_irred_struct} is a direct consequence of \Cref{positive_equilibrium_irred}.

            \section{Monotonicity properties (Theorems \ref{max_eq_decreasing} and \ref{s(A(a))_decreasing})}

            \subsection{Proof of \Cref{max_eq_decreasing}}

            \label{proof_of_max_eq_decreasing}

            Let us  prove that \eqref{decrease_eq} holds. The right-hand side of \eqref{log_allee} satisfies  $f_{a'} \le f_a$ for all $0_n \le a < a'$. Then, for any point $x_0 \in \mathbb{R}^n_+$, the respective solutions $x(t,x_0;a)$ and $x(t,x_0;a')$ of model \eqref{log_allee} with initial condition $x_0$ and corresponding to parameters $a$ and $a'$ satisfy
            \begin{equation}
                \label{xa' <xa}
                x(t,x_0;a') \le x(t,x_0;a), \quad t\geq0.
            \end{equation}
            Assume $x_{0,i} \geq  \max\{{x_+(a)}_i, {x_+(a')}_i\}$ for all $i = 1,\ldots,n$. Then by \Cref{0only}, it follows
            $$\lim_{t \rightarrow + \infty} x(t,x_0;a') = x_+(a'), \quad \lim_{t \rightarrow + \infty} x(t,x_0;a) = x_+(a).$$ By \eqref{xa' <xa}, this implies
            \begin{equation}
                \label{l1}
                x_+(a') \le x_+(a),
            \end{equation}
            so \eqref{decrease_eq} holds.

            Now assume $D$ is irreducible.  We  prove that when $x_+(a) > 0_n$, the inequality $\le$ in \eqref{decrease_eq} is replaced  by  $\ll$.

            If $x_+(a) > 0_n$, then, as the origin is always an equilibrium of \eqref{log_allee},  it follows by  \eqref{positive_eq_irred_struct} in \Cref{0only} that   $x_+(a) \gg  0_n$.  As a consequence, since $f_{a'}(x) < f_a(x)$ for every $x \gg 0_n$  and $a'>a$, we have by  \Cref{decrease_equilibrium_point_cooperative_init}  that
            $$        x\big(t,x_+(a);a'\big) \ll x\big(t,x_+(a);a\big), \quad t>0.  $$
            By \eqref{l1} and by cooperativeness of the system, one has  $x(t,x_+(a');a') \le x(t,x_+(a);a')$. Therefore,
            $$        x\big(t,x_+(a');a'\big) \ll x\big(t,x_+(a);a\big), \quad t>0.  $$
            By definition of an equilibrium, we have  $ x\big(t,x_+(a');a'\big) = x_+(a')$ and $ x\big(t,x_+(a);a\big) = x_+(a)$, and thus
            $$ x_+(a') \ll x_+(a).$$

            Last, we show that the basin of attraction $\mathcal{B}_a$ of $0_n$ is non-decreasing with $a$. Let $a' \geq a$, and let $x_0 \in \mathcal{B}_a \subset \mathbb{R}^n_+$.  By definition of $\mathcal{B}_a$, $x(t,x_0, a) \rightarrow 0_n$ as $t \rightarrow + \infty$. Since $a' \geq a$, $f_{a'}\le f_a$. Therefore, the monotony of  system \eqref{log_allee} implies $x(t,x_0; a') \le x(t,x_0; a)$. By comparison and since it can easily be shown that $\mathbb{R}^n_+$ is positively invariant with respect to system \eqref{log_allee},  $x(t,x_0; a') \rightarrow 0_n$ as $t \rightarrow + \infty$, which proves that $x_0 \in \mathcal{B}_{a'}$.
            This  completes the proof of \Cref{max_eq_decreasing}.

            \subsection{Proof of \Cref{s(A(a))_decreasing}}

            \label{proof_s(A)_decreasing}

            To prove \Cref{s(A(a))_decreasing}, we study the function $a\in \mathbb{R}^n_+ \mapsto s(A(a))$. One has
            $$A(a) = \diag\big(\varphi_i\left(a_i \right)\big)- \diag(\mu_{1,i}+\rho_i) + D,$$
            where, for every $i = 1,\ldots,n$, the function $\varphi_i$ is defined by
            \begin{align*}
                \varphi_i(z) & :=  \left(\left(\sqrt{b_i} - \sqrt{z \mu_{2,i}} \right)^+\right)^2.
            \end{align*}

            Since $\varphi_i$ is non-increasing,
            $$a < a' \Rightarrow A(a') \le A(a).$$
            Extending  \cite[Corollary 1.5, Chapter 2]{nonnegative} to Metzler matrices  then implies the monotonicity statement \eqref{non-decreasing A(a)}.

            Moreover, $\varphi_i$  is  continuous on $\mathbb{R}_+$, which implies that the function  $ a\mapsto A(a)$ is continuous on $\mathbb{R}^n_+$. Its  derivative  on $\left(0, \frac{b_i}{\mu_{2,i}}\right]$ is:
\begin{equation}
    \label{deriv_phi}
    \varphi_i'(z) = \mu_{2,i} -  \sqrt{\dfrac{b_i \mu_{2,i}}{z}}, \quad z \in \left(0, \dfrac{b_i}{\mu_{2,i}}\right].
\end{equation}
This derivative is continuous,  showing that $\varphi_i$ is continuously differentiable on $\left(0, \frac{b_i}{\mu_{2,i}}\right]$. Moreover, $\varphi_i'(\frac{b_i}{\mu_{2,i}}) = 0$ and $\varphi_i$ is constant on $\left(\frac{b_i}{\mu_{2,i}}, +\infty\right)$. Thus $\varphi_i$ is continuously differentiable on $\mathbb{R}^*_+$,  implying that the function $a \mapsto A(a)$ is continuously differentiable on $\{a \gg 0_n\}$. Notice, however, that $\varphi_i$ is not differentiable at $z = 0$, since $\varphi_i'(z) \to -\infty$ as $z \to 0^+$.

For $z>0$,  the second derivative is:

$$\varphi_i''(z)  =  \left\{
    \begin{array}{ll}
        \dfrac{1}{2z} \sqrt{\dfrac{b_i \mu_{2,i}}{z}} & \mbox{if } 0 < z \le \dfrac{b_i}{\mu_{2,i}} \\
        0                                             & \mbox{if } z > \dfrac{b_i}{\mu_{2,i}}.
    \end{array}
    \right.$$
Thus, $\varphi_i$ is piecewise twice continuously differentiable on $\mathbb{R}^*_+$, but not twice continuously differentiable on $\mathbb{R}^*_+$  since  $\varphi_i''(z) \neq 0$ when $z = \frac{b_i}{\mu_{2,i}}$.  This implies that $a \mapsto A(a)$ is piecewise twice continuously differentiable on the set $\{a\gg 0_n\}$.

Finally,   $\varphi_i''(z) =  \frac{1}{2z} \sqrt{\dfrac{b_i \mu_{2,i}}{z}} > 0$  for any $z \in \left(0,\frac{b_i}{\mu_{2,i}}\right]$.  Since in addition,  $\varphi_i'$ is
constant on the interval  $\left(\frac{b_i}{\mu_{2,i}}, + \infty\right)$, and $\varphi_i$ is differentiable on $\mathbb{R}^*_+$, it follows that  $\varphi_i$ is convex on $\mathbb{R}^*_+$.  We now prove that $\varphi_i$ is convex on the whole domain $\mathbb{R}_+$.  For any $z,z' \in \mathbb{R}_+$, there exist sequences $z_n$ and $z'_n \in \left(\mathbb{R}^*_+\right)^\mathbb{N}$ such that $z_n \to z$ and $z_n' \to z'$ as $n \to \infty$. By continuity and by convexity of $\varphi_i$ on $\mathbb{R}^*_+$, for all $t \in [0,1]$,
$$\varphi_i\big(t z +(1-t)z'\big) = \lim_{t \to +\infty} \varphi_i\big(t z_n +(1-t)z'_n\big) \le \lim_{t \to \infty} t \varphi_i( z_n) + (1-t) \varphi_i(z'_n) = t \varphi_i(z) + (1-t) \varphi_i(z'),$$
which proves that $\varphi_i$ is convex on the whole domain $\mathbb{R}_+$.

The convexity of $s(A(a))$  then follows from \cite[Theorem 3.3.5]{kirkland2012group}.
Now assume that $D$ is irreducible.  Then as a consequence of \cite[Theorem 3.2.1]{kirkland2012group}, the stability modulus of $A(a)$ inherits the differentiability properties of $A(a)$.
If moreover, $a_i < \frac{b_i}{\mu_{2,i}}$ for all $i = 1,\ldots,n$, then,  since the function $\varphi_i$ is
strictly decreasing on $\left[0, \frac{b_i}{\mu_{2,i}}\right)$,
$$a < a' \Rightarrow A(a') < A(a).$$
This implies by  \cite[Lemma 1.2]{nussbaum1986convexity} extended to Metzler matrices that the inequality $\le$ in  \eqref{non-decreasing A(a)}
may be replaced by  $<$.

\section{Sterile insect equilibrium for reducible network (\Cref{sterile_equilibrium_reducible} and \Cref{lemma_supremal_sterile_eq})}

\subsection{Proof of \Cref{sterile_equilibrium_reducible}}

\label{proof_sterile_reducible}

Let  $\mathcal{G}$ be an SCC of $\Gamma_s$. By definition of the set $\mathcal{G}^-$ of its  in-neighbouring SCCs, one has $d^s_{ji} = 0$ for every $j \in V_{\mathcal{G}^-}$ and $i \in V_{\mathcal{G}}$. Similarly, one has $d^s_{ij} = 0$ for every $j \in V_{\mathcal{G}^+}$ and $i \in V_{\mathcal{G}}$. It follows that, for any SCC $\mathcal{G}$ of $\Gamma_s$, the system \eqref{sterile_model} writes

$$\dot{x}_{s,i} = \Lambda_i - \left(\mu_{s,i} + \rho_{s,i}  +\sum\limits_{\substack{j \in  V_{\mathcal{G}^+}}} {d^s_{ji}}\right)x_{s,i} + \sum\limits_{\substack{j \in V_{\mathcal{G}} \\ j \neq i}}{d^s_{ij} x_{s,j}} -\sum\limits_{\substack{j \in V_{\mathcal{G}} \\ j \neq i}}{d^s_{ji} x_{s,i}} +  \sum\limits_{\substack{j \in V_{\mathcal{G}^-}}}{d^s_{ij} x_{s,j}}, ~~~~~i \in V_{\mathcal{G}}.$$
To find the value of $x^*_s(\Lambda)|_{V_{\mathcal{G}}}$, we replace $x_{s,j}$  by  $x^*_{s,j}(\Lambda)$  for every  $j \in  V_{\mathcal{G}^-}$. This results in the following system:

\begin{equation}
    \label{N = 2,...,n}
    \dot{x}_{s,i} = \Lambda_i  +  \sum\limits_{\substack{j \in V_{\mathcal{G}^-}}}{d^s_{ij} x^*_{s,j}(\Lambda)}- \left(\mu_{s,i} + \rho_{s,i} +\sum\limits_{\substack{j \in  V_{\mathcal{G}^+}}} {d^s_{ji}}\right)x_{s,i} + \sum\limits_{\substack{j \in V_{\mathcal{G}} \\ j \neq i}}{d^s_{ij} x_{s,j}} -\sum\limits_{\substack{j \in V_{\mathcal{G}} \\ j \neq i}}{d^s_{ji} x_{s,i}}, ~~~~~i \in V_{\mathcal{G}}.
\end{equation}
Therefore, the equilibrium of \eqref{N = 2,...,n} satisfies \eqref{eq_SCC}.  The matrix $D^s|_{V_{\mathcal{G}}}$ is irreducible by definition of an $SCC$, which implies that the matrix $J^s|_{V_{\mathcal{G}}}$ is also irreducible.  Moreover, $J^s|_{V_{\mathcal{G}}}$ is a Metzler and Hurwitz matrix, so $- \left(J^s|_{V_{\mathcal{G}}}\right)^{-1} \gg 0_{n_{V_{\mathcal{G}}} \times n_{V_{\mathcal{G}}}}$ by \cite[Theorem 10.3]{bullo2018lectures}. We can then  prove \eqref{eq_sterile_red} in \Cref{sterile_equilibrium_reducible} by induction, following the acyclic ordering of the SCCs of the sterile insect network.

\textbf{Initialisation:}
First, assume that $\mathcal{G}$ is of in-degree zero. Then $V_{\mathcal{G}_-} = \emptyset$, which implies by  \eqref{eq_SCC} that
\begin{equation*}
    x^*_s(\Lambda)|_{V_{\mathcal{G}}} = - \left(J^s|_{V_{\mathcal{G}}}\right)^{-1} \Lambda|_{V_{\mathcal{G}}}.
\end{equation*}

\begin{itemize}
    \item Assume  $V_{\mathcal{G}}  \cap  \mathcal{C}^s \neq \emptyset$. Then, $\Lambda|_{V_{\mathcal{G}}}   > 0_{n_{V_{\mathcal{G}}}}$, which implies
          $x^*_s(\Lambda)|_{V_{\mathcal{G}}} \gg 0_{n_{V_{\mathcal{G}}}}$.

    \item Assume $ V_{\mathcal{G}} \cap \mathcal{C}^s = \emptyset$. Then  $\Lambda|_{V_{\mathcal{G}}} = 0_{n_{V_{\mathcal{G}}}}$, which implies $x^*_s(\Lambda)|_{V_{\mathcal{G}}} = 0_{n_{V_{\mathcal{G}}}}$.
\end{itemize}
This proves that \eqref{eq_sterile_red} is true for any SCC of in-degree zero.

\vspace{0.5cm}

\textbf{Induction step:} Let $\mathcal{H}$ be an SCC of $\Gamma_s$ and assume \eqref{eq_sterile_red} is true for any SCC {$\mathcal{G}$} in $\mathcal{H}^\text{up}$.
\begin{itemize}
    \item Assume  $V_{\mathcal{H}}  \cap  \mathcal{C}^s \neq \emptyset$. Then, $\Lambda|_{V_{\mathcal{H}}}   > 0_{n_{V_{\mathcal{H}}}}$. Therefore, by the formula in \eqref{eq_SCC}, $x^*_s(\Lambda)|_{V_{\mathcal{H}}} \gg 0_{n_{V_{\mathcal{H}}}}$.

    \item Assume $ V_{\mathcal{H}^\text{up}} \cap \mathcal{C}^s \neq \emptyset$. Then, by definition of $\mathcal{H}^-$ and $\mathcal{H}^\text{up}$,  there exists an SCC $\mathcal{G}$  in $\mathcal{H}^-$ such that $\{V_{\mathcal{G}} \cup  V_{\mathcal{G}^\text{up}}\} \cap  \mathcal{C}^s \neq \emptyset$. Since \eqref{eq_sterile_red} is assumed true for $\mathcal{G}$, this implies that $x^*_s(\Lambda)|_{V_{\mathcal{G}}} \gg 0_{n_{V_{\mathcal{G}}}}$.  Thus, since there exists $i \in V_{\mathcal{H}}$ and $j \in V_{\mathcal{G}}$ such that $d^s_{ij} > 0$, one has $\Lambda|_{V_{\mathcal{H}}} +  \left( \sum\limits_{\substack{j \in V_{\mathcal{H}^-}}}{d^s_{ij} x^*_{s,j}(\Lambda)} \right)_{i \in V_{\mathcal{H}}} > 0_{n_{V_{\mathcal{H}}}}$. This implies  that $x^*_s(\Lambda)|_{V_{\mathcal{H}}} \gg 0_{n_{V_{\mathcal{H}}}}$.

    \item  Assume $ \{V_{\mathcal{H}} \cup V_{\mathcal{H}^{\text{up}}} \}\cap \mathcal{C}^s = \emptyset$. Then,  $\Lambda|_{V_{\mathcal{H}}}   = 0_{n_{V_{\mathcal{H}}}}$, and for any  SCC $\mathcal{G}$  in $\mathcal{H}^-$, one has $\{V_{\mathcal{G}} \cup  V_{\mathcal{G}^\text{up}}\} \cap  \mathcal{C}^s = \emptyset$. Since \eqref{eq_sterile_red} is assumed true for $\mathcal{G}$, this implies that $x^*_s(\Lambda)|_{V_{\mathcal{G}}} = 0_{n_{V_{\mathcal{G}}}}$. Therefore, $\Lambda|_{V_{\mathcal{H}}} +  \left( \sum\limits_{\substack{j \in V_{\mathcal{H}^-}}}{d^s_{ij} x^*_{s,j}(\Lambda)} \right)_{i \in V_{\mathcal{H}}} = 0_{n_{V_{\mathcal{H}}}}$, implying that $x^*_s(\Lambda)|_{V_{\mathcal{H}}} = 0_{n_{V_{\mathcal{H}}}}$.
\end{itemize}
By induction,  this concludes the proof of \Cref{sterile_equilibrium_reducible}.

\subsection{Proof of \Cref{lemma_supremal_sterile_eq}}

\label{proof_lemma_sup}

By \Cref{sterile_equilibrium_reducible}, for any SCC $\mathcal{G}$ of $\Gamma_s$, one has
$$
    x^*_s(\Lambda)|_{V_{\mathcal{G}}} = - \left(J^s|_{V_{\mathcal{G}}}\right)^{-1} \left( \Lambda|_{V_{\mathcal{G}}} + \left(\sum\limits_{\substack{j \in V_{\mathcal{G}^-}}}{d^s_{ij} x^*_{s,j}(\Lambda)} \right)_{i \in V_{\mathcal{G}}}  \right),
$$
with $- \left(J^s|_{V_{\mathcal{G}}}\right)^{-1} \gg 0_{n_{V_{\mathcal{G}}} \times n_{V_{\mathcal{G}}}}$. Then, by reasoning by induction on the SCCs of $\Gamma_s$ as in the proof of \Cref{sterile_equilibrium_reducible}, we obtain that, for all $k \in V_{\mathcal{G}}$  :
\begin{itemize}
    \item  If $\Lambda_i \rightarrow + \infty$ for some $i \in V_{\mathcal{G}} \cup  V_{\mathcal{G}^\text{up}} $, then, since $- \left(J^s|_{V_{\mathcal{G}}}\right)^{-1} \gg 0_{n_{V_{\mathcal{G}}} \times n_{V_{\mathcal{G}}}}$, we have
          $x^*_{s,k}(\Lambda)  \rightarrow    + \infty$.

    \item If $\Lambda_i = 0$ for all $i \in V_{\mathcal{G}} \cup  V_{\mathcal{G}^\text{up}}$, then   $x^*_{s,k}(\Lambda)  =  0$.

    \item  If $\Lambda_i \rightarrow \overline{\Lambda}_i \in [0, + \infty)$ for all $i \in V_{\mathcal{G}} \cup  V_{\mathcal{G}^\text{up}}$ with at least one $i \in V_{\mathcal{G}} \cup  V_{\mathcal{G}^\text{up}}$ such that $\Lambda_i >0$, then  since $- \left(J^s|_{V_{\mathcal{G}}}\right)^{-1} \gg 0_{n_{V_{\mathcal{G}}} \times n_{V_{\mathcal{G}}}}$, it follows that $ 0 < \lim_{\Lambda \rightarrow \overline{\Lambda}} x^*_{s,k}(\Lambda) < +\infty$.
\end{itemize}
Using the definition of $\mathcal{C}^s, \mathcal{C}^s_I$ and $\mathcal{C}^s_B$ in \Cref{def C^s} along with the argument above allows to complete the proof of \Cref{lemma_supremal_sterile_eq}.

\section{Gradient evaluation (\Cref{gradprobleme})}

\label{Gradient evaluation}

Let us compute the gradient  of the function $h(\Lambda) = s\big(A(a + \gamma x^*_s(\Lambda))\big) $.

Since the matrix  $D$ is assumed irreducible, the matrix $Q(\Lambda)$ admits zero as a simple eigenvalue for any $\Lambda \in \mathbb{R}^n_+$.  Consequently, it admits a group inverse, denoted by $Q^\#(\Lambda)$.

For any $\Lambda \gg 0_n$, we have as a consequence of \Cref{sterile_equilibrium_reducible} that $x^*_s(\Lambda) \gg 0_n$.  The function  $h$ is therefore  continuously differentiable with respect to $\Lambda \gg 0_n$ by \Cref{s(A)_decreasing}, and the gradient $\nabla h(\Lambda)$  is then  well-defined. Notice that we are only interested in strictly positive  values of $\Lambda$, since admissible points for the interior-point algorithm must strictly satisfy the inequality constraints.

For any $\Lambda \gg 0_n$, we have  by the chain rule that
\begin{equation}
    \label{grad1}
    \big(\nabla h(\Lambda)\big)_i = \sum\limits_{p = 1}^n  \frac{\partial A_{pp} (a + \gamma x^*_s(\Lambda))}{\partial \Lambda_i} \frac{\partial s}{\partial_{pp}}\big(A(a+\gamma x^*_s(\Lambda))\big), \qquad i = 1,\ldots,n,
\end{equation}
where  $\frac{\partial A_{pp} (a + \gamma x^*_s(\Lambda))}{\partial \Lambda_i}$ is the partial derivative of the $(p,p)$-th entry of $A(a + \gamma x^*_s(\Lambda))$ with respect to $\Lambda_i$, and $\frac{\partial s}{\partial_{pp}}\big(A(a+\gamma x^*_s(\Lambda))\big)$ is the partial derivative of the stability modulus with  respect to the $(p,p)$-th entry of $A(a + \gamma x^*_s(\Lambda))$.

By \cite[Lemma 3.1]{deutsch1984derivatives}, one has
\begin{equation}
    \label{grad2}
    \frac{\partial s}{\partial_{pp}}\big(A(a+\gamma x^*_s(\Lambda))\big) = \left(I- Q(\Lambda) Q^\#(\Lambda) \right)_{pp}, \qquad p = 1,\ldots, n.
\end{equation}
On the other hand,  for  any $\Lambda \gg 0_n$, straightforward computations yields
\begin{equation}
    \label{grad3}
    \frac{\partial A_{pp} (a + \gamma x^*_s(\Lambda))}{\partial \Lambda_i} =  \gamma   \mu_{2,p} \left(J^s\right)^{-1}_{pi}  \left( \sqrt{\frac{b_p }{\big(a_p + \gamma x^*_{s,p}(\Lambda)\big) \mu_{2,p}}} - 1\right)^+.
\end{equation}
Incorporating \eqref{grad2} and \eqref{grad3} into \eqref{grad1} allows to recover \eqref{grad_h} in \Cref{gradprobleme}.

\end{document}